\documentclass{amsart}
\usepackage{amssymb}
\usepackage{amscd}
\usepackage{verbatim}
\begin{document}
\newcommand\Mand{\ \text{and}\ }
\newcommand\Mwith{\ \text{with}\ }
\newcommand\Mfor{\ \text{for}\ }
\newcommand\Mst{\ \text{such that}\ }
\newcommand\Mor{\ \text{or}\ }
\newcommand\Mif{\ \text{if}\ }
\newcommand\Miff{\ \text{iff}\ }
\newcommand\Mthen{\ \text{then}\ }
\newcommand\nin{\notin}
\newcommand\identity{\operatorname{id}}
\newcommand\Id{\operatorname{Id}}
\newcommand\Real{\mathbb{R}}
\newcommand\pos{\Real^+}
\newcommand\Rnp{\Real\setminus\{0\}}
\newcommand\nzero{\setminus\{0\}}
\newcommand\Cx{\mathbb{C}}
\newcommand\Cxp{\Cx^+}
\newcommand\Cxm{\Cx^-}
\newcommand\Nat{\mathbb{N}}
\newcommand\halfNat{{\frac{1}{2}}\mathbb{N}}
\newcommand\intgr{\mathbb{Z}}
\newcommand\im{\operatorname{Im}}
\newcommand\re{\operatorname{Re}}
\newcommand\sign{\operatorname{sign}}
\newcommand\codim{\operatorname{codim}}
\newcommand\End{\operatorname{End}}
\newcommand\Ker{\operatorname{Ker}}
\newcommand\Hom{\operatorname{Hom}}
\newcommand\ideal{{\mathcal I}}
\newcommand\Span{\operatorname{span}}
\newcommand\image{\operatorname{image}}
\newcommand\Range{\operatorname{Ran}}
\newcommand\Graph{\operatorname{graph}}
\newcommand\slim{\operatornamewithlimits{s-lim}}
\newcommand\pa{\partial}
\newcommand\Rn{\Real^n}
\newcommand\Rm{\Real^m}
\newcommand\RN{\Real^N}
\newcommand\RtN{\Real^{2N}}
\newcommand\RM{\Real^M}
\newcommand\sphere{\mathbb{S}}
\newcommand\Sn{\sphere^{n-1}}
\newcommand\Sm{\sphere^{m-1}}
\newcommand\Snp{\sphere^n_+}
\newcommand\Smp{\sphere^m_+}
\newcommand\SN{\sphere^{N-1}}
\newcommand\SNp{\sphere^N_+}
\newcommand\circlep{\sphere^1_+}
\newcommand\Phom{P_{h}}
\newcommand\Shom{S_{h}}
\newcommand\distance{\operatorname{dist}}
\newcommand\cl{\operatorname{cl}}
\newcommand\interior{\operatorname{int}}
\newcommand\Fa{\operatorname{Fa}}
\newcommand\ff{\operatorname{ff}}
\newcommand\mf{\operatorname{mf}}
\newcommand\cf{\operatorname{cf}}
\newcommand\scf{\operatorname{sf}}
\newcommand\lf{\operatorname{lf}}
\newcommand\rf{\operatorname{rf}}
\newcommand\indfam{{\mathcal K}}
\newcommand\calA{{\mathcal A}}
\newcommand\calB{{\mathcal B}}
\newcommand\calR{{\mathcal R}}
\newcommand\calO{{\mathcal O}}
\newcommand\calX{{\mathcal X}}
\newcommand\calF{{\mathcal F}}
\newcommand\calG{{\mathcal G}}
\newcommand\calT{{\mathcal T}}
\newcommand\calC{{\mathcal C}}
\newcommand\calCt{{\tilde \mathcal C}}
\newcommand\calCL{{\mathcal C}_{\text L}}
\newcommand\calCR{{\mathcal C}_{\text R}}
\newcommand\Cinf{{\mathcal C}^{\infty}}
\newcommand\dist{{\mathcal C}^{-\infty}}
\newcommand\dCinf{\dot\Cinf}
\newcommand\ddist{\dot\dist}
\newcommand\Cj{{\mathcal C}^j}
\newcommand\Linf{L^{\infty}}
\newcommand\phg{{\text{phg}}}
\newcommand\bcon{{\mathcal A}}
\newcommand\bconc{{\mathcal A}_{\text{phg}}}
\newcommand\Sch{{\mathcal S}}
\newcommand\temp{\Sch^{\prime}}
\newcommand\Diff{\operatorname{Diff}}
\newcommand\Diffb{\operatorname{Diff}_{\text{b}}}
\newcommand\Diffc{\operatorname{Diff}_{\text{c}}}
\newcommand\Diffsc{\operatorname{Diff}_{\text{sc}}}
\newcommand\DiffI{\operatorname{Diff}_{\text{I}}}
\newcommand\DiffIq{\operatorname{Diff}_{\text{I},q}}
\newcommand\sing{\text{sing}}
\newcommand\supp{\operatorname{supp}}
\newcommand\ssupp{\operatorname{sing\ supp}}
\newcommand\csupp{\operatorname{cone\ supp}}
\newcommand\esupp{\operatorname{ess\ supp}}
\newcommand\Fr{{\mathcal F}}
\newcommand\Frinv{\Fr^{-1}}
\newcommand\bop{{\mathcal B}}
\newcommand\spec{\operatorname{spec}}
\newcommand\pspec{\spec_{pp}}
\newcommand\cspec{\spec_{c}}
\newcommand\FIO{{\mathcal I}}
\newcommand\SP{\operatorname{RC}}
\newcommand\RC{\operatorname{RC}}
\newcommand\Symc{S_c}
\newcommand\Symca{S_c^{\alpha}}
\newcommand\Symczero{S_c^{0,...,0}}
\newcommand\sci{{}^{\text{sc}}}
\newcommand\sct{\sci T^*}
\newcommand\scdt{\sci \dot T^*}
\newcommand\dS{\dot S^*}
\newcommand\dT{\dot T^*}
\newcommand\dSreg{\dot\Sigma_{\text reg}}
\newcommand\scct{\sci\bar{T}^*}
\newcommand\Csc{C_{\text{sc}}}
\newcommand\SNpscd{(\SNp)^2_{\text{sc}}}
\newcommand\scdiag{\Delta_{\text{sc}}}
\newcommand\projscl{\pi^L_{\text{sc}}}
\newcommand\projscr{\pi^R_{\text{sc}}}
\newcommand\scHL{\sci H^{2,0}_{|\zeta|^2-\lambda^2}}
\newcommand\scHrg{\sci H^{2,0}_{\sqrt{g}}}
\newcommand\Hsc{H_{\text{sc}}}
\newcommand\WF{\operatorname{WF}}
\newcommand\WFp{\operatorname{WF^{\prime}}}
\newcommand\WFsc{\operatorname{WF}_{\text{sc}}}
\newcommand\WFscp{\operatorname{WF_{sc}^{\prime}}}
\newcommand\WFC{\operatorname{WF}_C}
\newcommand\WFCi{\operatorname{WF}_{C_i}}
\newcommand\elliptic{\operatorname{ell}}
\newcommand\Psop{\operatorname{\Psi}}
\newcommand\Psiscrs{\operatorname{\Psi_{sc}^{-2,\infty}}}
\newcommand\Psiscr{\operatorname{\Psi_{sc}^{-2,0}}}
\newcommand\Psiscrm{\operatorname{\Psi_{sc}^{0,2}}}
\newcommand\PsiscHam{\operatorname{\Psi_{sc}^{2,0}}}
\newcommand\Psisci{\operatorname{\Psi_{sc}^{*,*}}}
\newcommand\Psiscid{\operatorname{\Psi_{sc}^{0,0}}}
\newcommand\Psiscis{\operatorname{\Psi_{sc}^{0,\infty}}}
\newcommand\Psiscsi{\operatorname{\Psi_{sc}^{-\infty,0}}}
\newcommand\Psiscs{\operatorname{\Psi_{sc}^{-\infty,\infty}}}
\newcommand\Psiscalg{\operatorname{\Psi_{sc}^{\infty,-\infty}}}
\newcommand\nullHam{{\mathcal N}}
\newcommand\charD{\Sigma_{\Delta-\lambda^2}}
\newcommand\charLap{\Sigma_{\Delta-\lambda}}
\newcommand\Snl{\Sn_{\lambda}}
\newcommand\SNl{\SN_{\lambda}}
\newcommand\gammat{\tilde\gamma}
\newcommand\gammasc{\gamma}
\newcommand\Tau{\mathcal{T}}
\newcommand\taut{\tilde\tau}
\newcommand\taub{\bar\tau}
\newcommand\Nout{N^+_{\lambda}}
\newcommand\Nin{N^-_{\lambda}}
\newcommand\Nio{N^{\pm}_{\lambda}}
\newcommand\El{E_{\lambda}}
\newcommand\Elt{\tilde E_{\lambda}}
\newcommand\Eil{E^i_{\lambda}}
\newcommand\Ejl{E^j_{\lambda}}
\newcommand\Eajl{E^{\alpha_j}_{\lambda}}
\newcommand\Eilt{\tilde E^i_{\lambda}}
\newcommand\Np{N^+}
\newcommand\Nm{N^-}
\newcommand\Npm{N^{\pm}}
\newcommand\Fin{F^-(\lambda)}
\newcommand\Fini{F^-_i(\lambda)}
\newcommand\Fout{F^+(\lambda)}
\newcommand\Fouti{F^+_i(\lambda)}
\newcommand\Foutj{F^+_j(\lambda)}
\newcommand\Rout{R^+_{\lambda}}
\newcommand\Routl{R^+_{\lambda^2}}
\newcommand\Routsgnl{R^{\sign\lambda}_{\lambda^2}}
\newcommand\Rin{R^-_{\lambda}}
\newcommand\Rinl{R^-_{\lambda^2}}
\newcommand\Rinsgnl{R^{-\sign\lambda}_{\lambda^2}}
\newcommand\Rio{R^{\pm}_{\lambda}}
\newcommand\Riol{R^{\pm}_{\lambda^2}}
\newcommand\Roi{R^{\mp}_{\lambda}}
\newcommand\Roil{R^{\mp}_{\lambda^2}}
\newcommand\Riob{R^{\pm}}
\newcommand\Roib{R^{\mp}}
\newcommand\Tio{T^{\pm}}
\newcommand\Tiob{T^{\pm}_{\ff}}
\newcommand\Toi{T^{\mp}}
\newcommand\Toib{T^{\mp}_{\ff}}
\newcommand\TIiob{T_I^{\pm}}
\newcommand\Rinb{R^-}
\newcommand\Rinbsgnl{R^{-\sign\lambda}}
\newcommand\Tin{T^-}
\newcommand\Tinb{T^-_{\ff}}
\newcommand\TIinb{T^-_I}
\newcommand\Routb{R^+}
\newcommand\Routbsgnl{R^{\sign\lambda}}
\newcommand\Tout{T^+}
\newcommand\Toutb{T^+_{\ff}}
\newcommand\TIoutb{T^+_I}
\newcommand\Rlkf{(|\xib|^2-(\lambda-i0)^2)^{-1}}
\newcommand\Rlk{\rho_0(\lambda)}
\newcommand\Rmlk{\rho_0(-\lambda)}
\newcommand\Rpmlk{\rho_0(\pm\lambda)}
\newcommand\Rlka{\rho_1(\lambda)}
\newcommand\Rlkb{\rho_2(\lambda)}
\newcommand\Rilk{\rho_i(\lambda)}
\newcommand\reduced{\natural}
\newcommand\Rlf{R_0(\lambda)}
\newcommand\Rla{R_1(\lambda)}
\newcommand\Rlb{R_2(\lambda)}
\newcommand\Ril{R_i(\lambda)}
\newcommand\Rlj{R_j(\lambda)}
\newcommand\Rlft{R_0(\lambda)}
\newcommand\Rflambda{R_0^{\reduced}(\sigma)}
\newcommand\RV{R^{\reduced}_V}
\newcommand\Rfsigma{R_0^{\reduced}(\sigma)}
\newcommand\Rfsigmah{R_0^{\reduced}(\sigma^{1/2})}
\newcommand\Rfzero{R_0^{\reduced}(0)}
\newcommand\RlV{R^{\reduced}_V(\sigma)}
\newcommand\RlVi{R^{\reduced}_{V_i}(\sigma)}
\newcommand\RlVt{R_V(\lambda)}
\newcommand\RlVtL{{R}_V^L(\lambda)}
\newcommand\RlVtR{{R}_V^R(\lambda)}
\newcommand\RlVit{{R}_{V_i}(\lambda)}
\newcommand\RlVta{{R}_V^{(1)}(\lambda)}
\newcommand\RlVtk{{R}_V^{(k)}(\lambda)}
\newcommand\RlVatV{{R}_{V_{\alpha}}(\lambda)V_{\alpha}}
\newcommand\RlVatVa{{R}_{V_{\alpha_1}}(\lambda)V_{\alpha_1}}
\newcommand\RlVatVb{{R}_{V_{\alpha_2}}(\lambda)V_{\alpha_2}}
\newcommand\RlVatVk{{R}_{V_{\alpha_k}}(\lambda)V_{\alpha_k}}
\newcommand\RlVatVkk{{R}_{V_{\alpha_{k+1}}}(\lambda)V_{\alpha_{k+1}}}
\newcommand\RlVaptV{{R}_{V_{\alpha'}}(\lambda)V_{\alpha'}}
\newcommand\RlVapptV{{R}_{V_{\alpha''}}(\lambda)V_{\alpha''}}
\newcommand\RlVajtV{{R}_{V_{\alpha_j}}(\lambda)V_{\alpha_j}}
\newcommand\RlVaktV{{R}_{V_{\alpha_k}}(\lambda)V_{\alpha_k}}
\newcommand\RlVakktV{{R}_{V_{\alpha_{k+1}}}(\lambda)V_{\alpha_{k+1}}}
\newcommand\Tl{T(\lambda)}
\newcommand\Tlt{\tilde\Tl}
\newcommand\Tltp{\tilde T'(\lambda)}
\newcommand\Tltpp{\tilde T''(\lambda)}
\newcommand\Tli{T_i(\lambda)}
\newcommand\Tlit{\tilde\Tli}
\newcommand\Tlip{T_i'(\lambda)}
\newcommand\Tlipp{T_i''(\lambda)}
\newcommand\Tlj{T_j(\lambda)}
\newcommand\Tla{T_{\alpha}(\lambda)}
\newcommand\Tlaa{T_{\alpha_1}(\lambda)}
\newcommand\Tlab{T_{\alpha_2}(\lambda)}
\newcommand\Tlak{T_{\alpha_k}(\lambda)}
\newcommand\Tlakt{\tilde\Tlak}
\newcommand\Tlaj{T_{\alpha_j}(\lambda)}
\newcommand\Tlajj{T_{\alpha_{j+1}}(\lambda)}
\newcommand\Tlajp{T_{\alpha_j}'(\lambda)}
\newcommand\Tlajpt{\tilde\Tlajp}
\newcommand\Tlajt{\tilde\Tlaj}
\newcommand\Tlakk{T_{\alpha_{k+1}}(\lambda)}
\newcommand\Tlakkp{T_{\alpha_{k+1}}'(\lambda)}
\newcommand\Tlap{T_{\alpha'}(\lambda)}
\newcommand\Tlapt{\tilde\Tlap}
\newcommand\Tlapp{T_{\alpha''}(\lambda)}
\newcommand\Tkl{T^{(k)}(\lambda)}
\newcommand\Tcl{T^{\flat}(\lambda)}
\newcommand\Fl{F(\lambda)}
\newcommand\BlVt{\tilde B_V(\lambda)}
\newcommand\KBlVt{K_{\BlVt}}
\newcommand\BlVaat{B_{V_{\alpha_1}}(\lambda)}
\newcommand\BV{B_V}
\newcommand\Bone{B_1}
\newcommand\Btwo{B_2}
\newcommand\Bthree{B_3}
\newcommand\Banyj{B_j}
\newcommand\PlV{P_V(\lambda)}
\newcommand\PlVc{P_V^{\flat}(\lambda)}
\newcommand\Pl{P_0(\lambda)}
\newcommand\SVl{S_V(\lambda)}
\newcommand\Sjr{S_j^{\reduced}}
\newcommand\Rkp{{\mathcal R}^k_+}
\newcommand\Rkm{{\mathcal R}^k_-}
\newcommand\Rkpm{{\mathcal R}^k_{\pm}}
\newcommand\Phys{{\mathcal P}}
\newcommand\Pc{\overline{\mathcal P}}
\newcommand\pip{\pi^{\perp}}
\newcommand\pipa{\pi_\partial}
\newcommand\gammapa{\gamma_\partial}
\newcommand\pipah{\hat\pi_\partial}
\newcommand\pit{\tilde\pi}
\newcommand\xit{\tilde\xi}
\newcommand\zetat{\tilde\zeta}
\newcommand\etat{\tilde\eta}
\newcommand\sigmat{\tilde\sigma}
\newcommand\sigmahat{\hat\sigma}
\newcommand\thetat{\tilde\theta}
\newcommand\psit{\tilde\psi}
\newcommand\phit{\tilde\phi}
\newcommand\chit{\tilde\chi}
\newcommand\rhot{\tilde\rho}
\newcommand\xib{\bar\xi}
\newcommand\zetab{\bar\zeta}
\newcommand\thetab{\bar\theta}
\newcommand\etab{\bar\eta}
\newcommand\iotal{\iota_{\lambda}}
\newcommand\rhoat{\rhot_{\alpha_1}}
\newcommand\Lambdat{\tilde\Lambda}
\newcommand\poles{\Lambda'}
\newcommand\rpoles{\Lambda_p}
\newcommand\thresholds{\Lambda}
\newcommand\Vt{\tilde V}
\newcommand\It{\tilde I}
\newcommand\half{{\frac{1}{2}}}
\newcommand\sigmah{\sigma^{1/2}}
\newcommand\bX{\partial X}
\newcommand\bXb{\partial \Xb}
\newcommand\Deltabt{\tilde\Delta_0}
\newcommand\strip{\Omega_T}
\newcommand\Kf{K^{\flat}}
\newcommand\Gs{G^{\sharp}}
\newcommand\Gt{\tilde G}
\newcommand\Osc{\sci\Omega}
\newcommand\OSc{{}^\Scl\Omega}
\newcommand\Osch{\sci\Omega^{\half}}
\newcommand\Oscmh{\sci\Omega^{-\half}}
\newcommand\Isc{I_{sc}}
\newcommand\Qzl{Q^0_{-\lambda}}
\newcommand\Lie{{\mathcal L}}
\newcommand\bl{{\text b}}
\newcommand\scl{{\text{sc}}}
\newcommand\sccl{{\text{scc}}}
\newcommand\Scl{{\text{Sc}}}
\newcommand\ScLl{{\text{Sc,L}}}
\newcommand\ScRl{{\text{Sc,R}}}
\newcommand\Sccl{{\text{Scc}}}
\newcommand\sus{{\text{sus}}}
\newcommand\XXb{X^2_\bl}
\newcommand\XXbt{\Xt^2_\bl}
\newcommand\XXsc{X^2_\scl}
\newcommand\XXsct{\Xt^2_\scl}
\newcommand\XXSc{X^2_\Scl}
\newcommand\XXSct{\Xt^2_\Scl}
\newcommand\XXScL{X^2_\ScLl}
\newcommand\XXScR{X^2_\ScRl}
\newcommand\MMsc{M^2_\scl}
\newcommand\Deltab{\Delta_\bl}
\newcommand\Deltasc{\Delta_\scl}
\newcommand\DeltaSc{\Delta_\Scl}
\newcommand\DeltaScL{\Delta_\ScLl}
\newcommand\DeltaScR{\Delta_\ScRl}
\newcommand\prs{\sigma}
\newcommand\Nsc{N_\scl}
\newcommand\Nscp{N_{\scl,p}}
\newcommand\Nff{N_{\ff}}
\newcommand\Nffz{N_{\ff,0}}
\newcommand\Nffzp{N_{\ff,0,p}}
\newcommand\Nffl{N_{\ff,l}}
\newcommand\Nffml{N_{\ff,-l}}
\newcommand\Nmf{N_{\mf}}
\newcommand\Nmfz{N_{\mf,0}}
\newcommand\Nmfl{N_{\mf,l}}
\newcommand\Nmfml{N_{\mf,-l}}
\newcommand\ffb{\operatorname{bf}}
\newcommand\Ffb{\operatorname{bf'}}
\newcommand\ffsc{\operatorname{sf}}
\newcommand\ffSc{\operatorname{sf_C}}
\newcommand\Ffsc{\operatorname{sf'}}
\newcommand\rff{\rho_{\ff}}
\newcommand\rmf{\rho_{\mf}}
\newcommand\rffb{\rho_{\ffb}}
\newcommand\rffsc{\rho_{\ffsc}}
\newcommand\rFfsc{\rho_{\Ffsc}}
\newcommand\rffSc{\rho_{\ffSc}}
\newcommand\rinf{\rho_{\infty}}
\newcommand\CL{C_L}
\newcommand\CR{C_R}
\newcommand\betab{\beta_\bl}
\newcommand\betasc{\beta_\scl}
\newcommand\betaSc{\beta_\Scl}
\newcommand\BetaSc{\bar\beta_\Scl}
\newcommand\betaScL{\beta_\ScLl}
\newcommand\betaScR{\beta_\ScRl}
\newcommand\ScT{{}^\Scl T^*}
\newcommand\SccT{{}^\Scl \bar T^*}
\newcommand\ScS{{}^\Scl S^*}
\newcommand\Tb{{}^\bl T}
\newcommand\Tsc{{}^\scl T}
\newcommand\TSc{{}^\Scl T}
\newcommand\CSc{C_\Scl}
\newcommand\Lambdasc{{}^\scl\Lambda}
\newcommand\XXXb{X^3_\bl}
\newcommand\XXXsc{X^3_\scl}
\newcommand\XXXSc{X^3_\Scl}
\newcommand\XXXScO{X^3_{\Scl,O}}
\newcommand\XXXScF{X^3_{\Scl,F}}
\newcommand\XXXScS{X^3_{\Scl,S}}
\newcommand\XXXScC{X^3_{\Scl,C}}
\newcommand\KDsc{\operatorname{KD^{\half}_\scl}}
\newcommand\KDSc{\operatorname{KD^{\half}_\Scl}}
\newcommand\KDScEF{\operatorname{KD^{E,F}_\Scl}}
\newcommand\Oh{\operatorname{\Omega^{\half}}}
\newcommand\WFSc{\WF_\Scl}
\newcommand\WFtSc{\WF_{\text 3sc}}
\newcommand\WFScmf{\WF_{\Scl,\mf}}
\newcommand\WFScff{\WF_{\Scl,\ff}}
\newcommand\WFScs{\WF_{\Scl,\prs}}
\newcommand\WFScp{\WF'_\Scl}
\newcommand\WFScmfp{\WF'_{\Scl,\mf}}
\newcommand\WFScffp{\WF'_{\Scl,\ff}}
\newcommand\WFScsp{\WF'_{\Scl,\prs}}
\newcommand\Diffscc{\Diff_\sccl}
\newcommand\DiffSc{\Diff_\Scl}
\newcommand\DiffScc{\Diff_\Sccl}
\newcommand\DiffscI{\Diff_{\scl,\text{I}}}
\newcommand\VscI{\Vf_{\scl,\text{I}}}
\newcommand\DiffsV{\operatorname{Diff}_{\sus(V)}}
\newcommand\DiffsVsc{\operatorname{Diff}_{\sus(V),\scl}}
\newcommand\DiffsVCsc{\operatorname{Diff}_{\sus(V)-C,\scl}}   
\newcommand\Psisc{\Psop_\scl}
\newcommand\Psiscc{\Psop_\sccl}
\newcommand\Psisch{\Psop_{\scl,h}}
\newcommand\Psiscch{\Psop_{\sccl,h}}
\newcommand\PsiSc{\Psop_\Scl}
\newcommand\PsiScph{\Psop_{\Scl,\phi}}
\newcommand\PsiScra{\Psop_{\Scl,\rho^\sharp_a}}
\newcommand\PsiScc{\Psop_\Sccl}
\newcommand\PsiSccml{\Psop^{m,l}_\Sccl}
\newcommand\PsiScxx{\Psop^{*,*}_\Scl}
\newcommand\PsiScml{\Psop^{m,l}_\Scl}
\newcommand\PsiScmz{\Psop^{m,0}_\Scl}
\newcommand\PsiScmmz{\Psop^{-m,0}_\Scl}
\newcommand\PsiSckz{\Psop^{k,0}_\Scl}
\newcommand\PsiScmmml{\Psop^{-m,-l}_\Scl}
\newcommand\Psiscmkk{\Psop^{-k,k}_\scl}
\newcommand\Psiscmmmkk{\Psop^{-m-k,k}_\scl}
\newcommand\Psiscmoo{\Psop^{-1,1}_\scl}
\newcommand\Psiscmz{\Psop^{m,0}_\scl}
\newcommand\Psiscmmz{\Psop^{-m,0}_\scl}
\newcommand\PsiSckmkl{\Psop^{km,kl}_\Scl}
\newcommand\PsiScmplp{\Psop^{m',l'}_\Scl}
\newcommand\PsiScmmpllp{\Psop^{m+m',l+l'}_\Scl}
\newcommand\Psiscml{\Psop^{m,l}_\scl}
\newcommand\PsiScid{\Psop^{0,0}_\Scl}
\newcommand\PsiSczo{\Psop^{0,1}_\Scl}
\newcommand\PsiScmii{\Psop^{-\infty,\infty}_\Scl}
\newcommand\PsiScmiz{\Psop^{-\infty,0}_\Scl}
\newcommand\PsiScmoo{\Psop^{-1,1}_\Scl}
\newcommand\PsisCid{\Psop^{0,0}_{\scl-C}}
\newcommand\PsisC{\Psop_{\scl-C}}
\newcommand\Psiinf{\Psop_{\infty}}
\newcommand\Psiinfid{\Psop_{\infty}^0}
\newcommand\PsiFinf{\Psop_{\infty-\Fr}}
\newcommand\PsisVscml{\Psop^{m,l}_{\sus(V),\scl}}
\newcommand\PsisVsc{\Psop_{\sus(V),\scl}}
\newcommand\PsisVpsc{\Psop_{\sus(V_p),\scl}}
\newcommand\PsisVCSc{\Psop_{\sus(V)-C,\scl}}
\newcommand\SFinf{S_{\infty-\Fr}}
\newcommand\YsVC{Y^2_{\sus(V)-C,\scl}}
\newcommand\ffYsc{\ffsc_{\sus(V)}}
\newcommand\SXC{S(X;C)}
\newcommand\Ios{I_{\text{os}}}
\newcommand\pbL{\pi^2_{\bl,{\text L}}}
\newcommand\pbR{\pi^2_{\bl,{\text R}}}
\newcommand\pscL{\pi^2_{\scl,{\text L}}}
\newcommand\pscR{\pi^2_{\scl,{\text R}}}
\newcommand\PbO{\pi^3_{\bl,{\text O}}}
\newcommand\PscO{\pi^3_{\scl,{\text O}}}
\newcommand\PScO{\pi^3_{\Scl,{\text O}}}
\newcommand\PScF{\pi^3_{\Scl,{\text F}}}
\newcommand\PScC{\pi^3_{\Scl,{\text C}}}
\newcommand\PScS{\pi^3_{\Scl,{\text S}}}
\newcommand\pScL{\pi^2_{\Scl,{\text L}}}
\newcommand\pScR{\pi^2_{\Scl,{\text R}}}
\newcommand\CLF{\CL^F}
\newcommand\CLO{\CL^O}
\newcommand\CLS{\CL^S}
\newcommand\CLC{\CL^C}
\newcommand\DeltaYb{\Delta_{\bl,Y}}
\newcommand\DeltaYsc{\Delta_{\sus-\scl}}
\newcommand\diag{\operatorname{diag}}
\newcommand\Vf{{\mathcal V}}
\newcommand\Vb{{\mathcal V}_{\bl}}
\newcommand\Vsc{{\mathcal V}_{\scl}}
\newcommand\VSc{{\mathcal V}_{\Scl}}
\newcommand\VfI{\Vf_{\text{I}}}
\newcommand\VfIq{\Vf_{\text{I},q}}
\newcommand\scH{{}^\scl H}
\newcommand\scHg{\scH_g}
\newcommand\xh{\hat x}
\newcommand\Yh{\hat Y}
\newcommand\Zh{\hat Z}
\newcommand\Yb{\bar Y}
\newcommand\hb{\bar h}
\newcommand\xih{\hat\xi}
\newcommand\etah{\hat\eta}
\newcommand\muh{\hat\mu}
\newcommand\mub{\bar\mu}
\newcommand\nub{\bar\nu}
\newcommand\mubh{\widehat{\bar\mu}}
\newcommand\yb{\bar y}
\newcommand\zb{\bar z}
\newcommand\ub{\bar u}
\newcommand\Qb{\bar Q}
\newcommand\Wbp{{\bar W}^\perp}
\newcommand\Wp{W^\perp}
\newcommand\Kt{\tilde K}
\newcommand\Wt{\tilde W}
\newcommand\Ut{\tilde U}
\newcommand\yt{\tilde y}
\newcommand\ut{\tilde u}
\newcommand\vt{\tilde v}
\newcommand\ft{\tilde f}
\newcommand\htil{\tilde h}
\newcommand\St{\tilde S}
\newcommand\Pt{\tilde P}
\newcommand\Rt{\tilde R}
\newcommand\qt{\tilde q}
\newcommand\Qt{\tilde Q}
\newcommand\Xb{\bar X}
\newcommand\lambdat{\tilde\lambda}
\newcommand\betat{\tilde\beta}
\newcommand\epst{\tilde\epsilon}
\newcommand\ep{\epsilon}
\newcommand\bt{\tilde b}
\newcommand\Xt{\tilde X}
\newcommand\At{\tilde A}
\newcommand\at{\tilde a}
\newcommand\Ct{\tilde C}
\newcommand\pih{\hat\pi}
\newcommand\Rh{\hat R}
\newcommand\Ah{\hat A}
\newcommand\Bh{\hat B}
\newcommand\Ch{\hat C}
\newcommand\Gh{\hat G}
\newcommand\Hh{\hat H}
\newcommand\Qh{\hat Q}
\newcommand\Ph{\hat P}
\newcommand\Nh{\hat N}
\newcommand\Sh{\hat S}
\newcommand\Gcal{{\mathcal G}}
\newcommand\GcalC{{\mathcal G}_C}
\newcommand\Jcal{{\mathcal J}}
\newcommand\JcalC{{\mathcal J}_C}
\setcounter{secnumdepth}{3}
\newtheorem{lemma}{Lemma}[section]
\newtheorem{prop}[lemma]{Proposition}
\newtheorem{thm}[lemma]{Theorem}
\newtheorem{cor}[lemma]{Corollary}
\newtheorem{result}[lemma]{Result}
\newtheorem*{thm*}{Theorem}
\newtheorem*{prop*}{Proposition}
\newtheorem*{conj*}{Conjecture}
\numberwithin{equation}{section}
\theoremstyle{remark}
\newtheorem{rem}[lemma]{Remark}
\theoremstyle{definition}
\newtheorem{Def}[lemma]{Definition}
\newtheorem*{Def*}{Definition}
\def\signature#1#2{\par\noindent#1\dotfill\null\\*
{\raggedleft #2\par}}

\renewcommand{\theenumi}{\roman{enumi}}
\renewcommand{\labelenumi}{(\theenumi)}

\title[Many-body scattering]
{Propagation of singularities 
in many-body scattering in the presence of bound states}
\author[Andras Vasy]{Andr\'as Vasy}
\date{September 16, 1999}
\thanks{{\em Address:} Department of Mathematics, University of California, 
Berkeley, CA 94720-3840. {\em E-mail:} \texttt{andras@math.berkeley.edu}.}
\thanks{Partially supported by NSF grant \#DMS-99-70607.}

\begin{abstract}
In this paper we describe the propagation of singularities of
tempered distributional solutions $u\in\temp$ of $(H-\lambda)u=0$,
where $H$ is
a many-body Hamiltonian $H=\Delta+V$, $\Delta\geq 0$, $V=\sum_a V_a$,
and $\lambda$ is not a threshold of $H$,
under the assumption
that the inter-particle (e.g.\ two-body) interactions $V_a$
are real-valued
polyhomogeneous symbols of order $-1$ (e.g.\ Coulomb-type with the singularity
at the origin removed). Here the term `singularity' refers to a
microlocal description of the lack of decay at infinity.
Thus, we prove that the set of singularities of $u$ is a union of maximally
extended broken bicharacteristics of $H$. These are curves in the
characteristic variety of $H-\lambda$,
which can be quite complicated due to the
existence of bound states.
We use this result to describe
the wave front relation of the S-matrices. We also analyze Lagrangian
properties of this relation, which shows that the relation
is not `too large' in terms of its dimension.
\end{abstract}

\maketitle

\section{Introduction}
In this paper we describe the propagation of singularities of generalized
eigenfunctions of a many-body Hamiltonian $H=\Delta+V$, $V=\sum_a V_a$,
on $\Rn$
under the assumption that the inter-particle interactions $V_a$ are real-valued
polyhomogeneous symbols of order $-1$ (e.g.\ Coulomb-type with the singularity
at the origin removed). More precisely,
we use the `many-body scattering wave front set' $\WFSc(u)$ at infinity for
tempered distributions $u\in\temp(\Rn)$, and
prove that for $u\in\temp(\Rn)$ satisfying $(H-\lambda)u=0$,
$\WFSc(u)$ is a union of maximally
extended generalized
broken bicharacteristics of $H$, broken at the collision planes. Here
$\WFSc(u)$ provides a microlocal description of the lack of decay of
$u$ modulo $\Sch(\Rn)$, similarly to how the usual wave front set describes
distributions modulo $\Cinf$ functions.

The definition of generalized broken bicharacteristics is quite
technical due to the presence of bound states in the subsystems.
However, if these bound
states are absent, our definition becomes completely analogous to Lebeau's
definition \cite{Lebeau:Propagation}
for the wave equation in domains with corners. Indeed, in this
case the propagation result itself, which was proved
in \cite{Vasy:Propagation-Many}, is a direct ($\Cinf$-type)
analogue of Lebeau's result
for the propagation of analytic singularities for solutions of the wave
equation in domains with corners.

If there are bound states in the subsystems, but
either the set of thresholds
is discrete, or $H$ is a four-body Hamiltonian, the geometry of generalized
broken bicharacteristics is not much more complicated than in Lebeau's
setting.
The general
definition reflects that when particles collide, the total energy as
well as the external momentum is preserved. The complication in the
presence of bound states is that kinetic energy is not preserved, even
asymptotically. In summary, our results provide a connection between
quantum and classical objects, just as Lebeau's results connect the
wave equation and geometric optics. Note that these results only provide
the answer as to where in phase space the
singularities of generalized eigenfunctions may be located; they leave
open the question of what these singularities are like, i.e.\ we do not
have FIO-type results. Such results exist at least in certain
3-body settings \cite{Vasy:Structure, Hassell:Plane}, see the remarks
in the next section. In addition, we show that in our setting,
if the set of thresholds
is discrete, then the broken bicharacteristic relation has the
correct Lagrangian geometry to give rise to such FIO results.

We also prove the corresponding result in the `limiting absorption
principle' setting, namely that under certain assumptions on $\WFSc(f)$,
$R(\lambda\pm i0)f$ are defined, and $\WFSc(R(\lambda+i0)f)$
is a subset of the image of $\WFSc(f)\cup R_-(\lambda)$
under the forward broken bicharacteristic
relation. Here $R_-(\lambda)$ is the outgoing `radial set'.
Such a result makes the `radial-variable' propagation estimates
that have been used in many-body scattering, especially as derived in
the works of G\'erard, Isozaki and Skibsted \cite{GerComm, GIS:N-body},
more precise.

We use this result to analyze that
the wave front relation of the scattering matrices (S-matrices).
These connect the incoming and outgoing data of generalized eigenfunctions of
$H$, so one expects that their singularities are described by
considering the limit points
of generalized broken bicharacteristics $\gamma=\gamma(t)$ as $t\to
\pm\infty$. In fact,
in addition to the propagation of singularities result, the only ingredient
that
is required for this analysis is a good approximation for the incoming
Poisson operators
with incoming state $\alpha$ near the incoming region, and similar results
for the outgoing Poisson operators with, say, outgoing state $\beta$.
In general, one expects a WKB-type construction,
essentially as in Hadamard's parametrix construction. Indeed, this is what
Melrose and Zworski do in the geometric two-body type
setting, \cite{RBMZw}. In the
Euclidean many-body setting this construction has been done by Skibsted
\cite{Skibsted:Smoothness} in the short-range and by
Bommier~\cite{Bommier:Proprietes} in the long-range setting, in the
latter case by adopting
the construction of Isozaki and Kitada \cite{Isozaki-Kitada:Scattering},
at least under
the assumption that the energies of the states $\alpha$, $\beta$ are below
the continuous spectrum of the corresponding subsystem Hamiltonians.
Such a construction is unnecessary if $V_c$ are Schwartz, for then the
product decomposition is sufficiently accurate to give a good approximation
for the Poisson operator.
We thus obtain the following result.

\begin{thm*}
Suppose that $H$ is a many-body Hamiltonian, and
$\lambda$ is not a threshold or ($L^2$-)eigenvalue of $H$.
Suppose also that either $\alpha$ and $\beta$ are channels such that the
corresponding
eigenvalues $\ep_\alpha$, $\ep_\beta$, of the subsystem Hamiltonians $H^a$,
$H^b$, are in the discrete spectrum of $H^a$ and $H^b$ respectively, or
that $V_c$ is Schwartz for all $c$.
Then the wave front relation of the S-matrix $S_{\beta\alpha}(\lambda)$,
is given by the generalized broken bicharacteristic relation of $H$
as stated precisely below in Theorem~\ref{thm:sc-matrix}.
\end{thm*}

Special cases, which have already been analyzed, include the
free-to-free S-matrix in three-body scattering
\cite{Vasy:Structure, Hassell:Plane, Vasy:Propagation-2},
or indeed in many-body
scattering under the additional assumption that there are no
bound states in any subsystem  \cite{Vasy:Propagation-Many}.
In these cases the wave front relation
is given by the broken geodesic
relation, broken at the collision planes,
on $\Sn$ at distance $\pi$. In both cases, one can naturally extend the
results to geometric many-body type problems on asymptotically
Euclidean manifolds.

Also, Bommier \cite{Bommier:Proprietes} and Skibsted \cite{Skibsted:Smoothness}
have shown that the kernels of the 2-cluster to free cluster and 2-cluster
to 2-cluster S-matrices are smooth (except for the diagonal singularity
if the 2-clusters are the same), and previously Isozaki had showed this
in the three-body setting \cite{Isozaki:Structures}. We remark that
(under our polyhomogeneous assumption) the proofs of Bommier and
Skibsted in fact show that the 2-cluster to same 2-cluster S-matrices
are (non-classical, if the potentials are long-range)
pseudo-differential operators, which differ from the identity operator
by operators of order $-1$ if the potentials are short range.
In our geometric normalization this
means that they are Fourier integral operators associated to the
geodesic flow on the sphere at infinity to distance $\pi$ (along the
cluster).

It may seem that the (generalized) broken bicharacteristic relation
is rather large since a single bicharacteristic can be continued in
many ways after it hits a collision plane. In fact, this is not the
case, since many
bicharacteristics do not hit these collision planes at all.
A more precise statement would be the that the wave front relation of
the S-matrices is given by a union of Lagrangian submanifolds,
so at least it has the same structure as if the S-matrices were (sums of)
Fourier integral operators.
Indeed, we prove that if the set of thresholds of $H$ is discrete,
e.g.\ if there are no bound states in any subsystems, then
the wave front relations of the S-matrices are given by
{\em finite} unions of smooth Lagrangian relations.

We remark that the results of this paper would remain valid if we assumed
only that $V_a\in S^{-\rho}(X^a)$, $\rho>0$, as customary. In fact, the
proof of the propagation of singularities for generalized eigenfunctions
remains essentially unchanged, and the only difference in the above Theorem
is that the parametrix for the Poisson operators is not as explicit,
cf.\ \cite{Vasy:Propagation-Many};
instead, one needs to use
the constructions of Isozaki-Kitada \cite{Isozaki-Kitada:Scattering}
(as presented by Skibsted and Bommier) directly. The reason for the
polyhomogeneous assumption is that the proofs are somewhat nicer,
especially in notation, and it is a particularly natural assumption
to make in the compactification approach we adopt.

Our main tool in proving the propagation of singularities results
consists of microlocally positive commutator estimates, i.e.\ on
the construction of operators which have a positive commutator with $H$
in the part of phase space, say $U$,
where we wish to conclude that a generalized
eigenfunction $u$ has no scattering wave front set. These commutators are
usually negative in another region of phase space, namely backwards
(or forwards, depending on the construction) along
generalized broken bicharacteristics through $U$.
We thus assume the absence
of this region from $\WFSc(u)$, and conclude that the positive
commutator region, $U$, is also missing from $\WFSc(u)$.
Such techniques have been used by H\"ormander, Melrose and
Sj\"ostrand \cite{Hor, Melrose-Sjostrand:I}
to show the propagation of singularities
for hyperbolic equations (real principal type propagation) such as the
wave equation, including
in regions with smooth boundaries. Indeed, the best way to interpret
our results is to say that $H-\lambda$ is hyperbolic at infinity. In two-body
scattering the analogy with the wave equation in domains without boundary
is rather complete; this was the basis of Melrose's proof of propagation
estimates for scattering theory for `scattering metrics' in
\cite{RBMSpec}. In many-body scattering, the lack of commutativity of
the appropriate pseudo-differential algebra, even to top order, makes
the estimates (and their proofs) more delicate. We remark that,
as can be seen directly from the approach we take, the wave front set
estimates can be easily turned into microlocal estimates on the resolvent
considered as an operator between weighted Sobolev spaces; wave front set
statements are a particularly convenient way of describing propagation.

Indeed, there is some freedom in the precise definition of the wave front
set; see the remarks preceeding the statement of Theorem~\ref{thm:prop-sing}.
The alternative
definitions differ slightly, but agree for generalized eigenfunctions,
and the corresponding propagation results are based on the same
positive commutator estimates. Thus, the reader may find
the explicit microlocal estimates of Section~\ref{sec:propagation}
particularly clear. However, piecing together these estimates to describe
propagation is more cumbersome than making simple geometric statements
based on wave front sets; for this reason we emphasize the latter.

Positive commutator estimates have also played a major role in
many-body scattering
starting with the work of Mourre \cite{Mourre-Absence}, Perry, Sigal
and Simon \cite{Perry-Sigal-Simon:Spectral}, Froese and Herbst
\cite{FroMourre}, Jensen \cite{Jensen:Propagation},
G\'erard, Isozaki and Skibsted \cite{GerComm, GIS:N-body}
and Wang \cite{XPWang}.
In particular,
the Mourre estimate is one of them; it estimates $i[H,w\cdot D_w+D_w\cdot w]$.
This and some other {\em global} positive commutator results
have been used to prove the global results mentioned in the first
paragraph about some of the S-matrices
with initial state
in a two-cluster. They also give the basis for the existence, uniqueness
and equivalence statements in our definition of the S-matrix by
asymptotic expansions; these statements are discussed in
\cite{IsoRad, IsoUniq, Vasy:Scattering} in more detail.

More delicate (and often time-dependent)
commutator estimates have been used in the proof
of asymptotic completeness. This completeness property of many-body
Hamiltonians was proved by Sigal and Soffer, Graf, Derezi\'nski and Yafaev
under different assumptions on the potentials and by different techniques
\cite{Sigal-Soffer:N, Sigal-Soffer:Long-range, Sigal-Soffer:Asymptotic,
Sigal-Soffer:Asymptotic-4,
Graf:Asymptotic, Derezinski:Asymptotic, Yaf}. In
particular, Yafaev's paper \cite{Yaf}
shows quite explicitly the importance of the
special structure of the Euclidean Hamiltonian. This
structure enables him to obtain
a (time-independent)
positive commutator estimate, which would not follow from the
indicial operator arguments of
\cite{Vasy:Propagation-2, Vasy:Propagation-Many} and the present paper,
and which is then
used to prove asymptotic completeness.

We briefly outline of the positive commutator proofs. They consist of two
parts: first, the
construction of a symbol (and an associated pseudo-differential
operator) that we claim has a microlocally positive commutator with $H$,
when localized in energy, and second, the proof that the commutator
is indeed positive in the appropriate part of the phase space. The first
part has much in common with the analysis of generalized broken
bicharacteristics, since both are intimately connected to various Hamilton
vector fields associated to the subsystems. The second part is essentially
a microlocal version of the proof of the Mourre estimate that was obtained
by Froese and Herbst \cite{FroMourre}. Indeed, we could follow the
full indicial operator version of this proof, as was done in
\cite{Vasy:Propagation-2} and \cite{Vasy:Propagation-Many}.
However, in the presence of bound states
in the proper
subsystems we would have to rely more heavily on the
approximate product structure of the Hamiltonian at each cluster.
Indeed, in the full geometric problem considered in
\cite{Vasy:Propagation-Many}, the $L^2$ eigenvalues and eigenfunctions
of the indicial operator of $H$ can vary as one moves along the collision
planes, which makes even the description of the bicharacteristics
more complicated.
However, it turns out that in the Euclidean many-body setting,
after an explicit calculation of the indicial operators,
we can use the Mourre estimate explicitly in
the normal (non-commutative)
variables at each cluster, and use the standard Poisson bracket
formula for the commutator in the tangential (commutative)
variables. Since it eliminates
the need to present arguments that are essentially simple (microlocal)
modifications of the Froese-Herbst proof of
the Mourre estimate, we adopt the second approach.

The structure of the paper is the following. In the next section we set up
the framework of many-body scattering, mostly following the book
of Derezi\'nski and G\'erard, \cite{Derezinski-Gerard:Scattering},
and we state the precise results on the propagation of singularities.
In Section~\ref{sec:sc-geometry} we relate this to the Melrose's
approach to scattering via compactification \cite{RBMSpec}. We also include
in this section the outline of the resolved space construction and
the definition of the algebra of
many-body scattering differential operators from
\cite{Vasy:Propagation-Many}. As we indicate in this section, from an
algebraic point of view, it is the lack of commutativity of this algebra
to `top weight' (at infinity) that gives rise to the breaks in the
generalized broken bicharacteristics along which singularities propagate.
In Section~\ref{sec:indicial} we recall
the definition and basic properties of
many-body scattering pseudo-differential operators as well as the
definition of the corresponding wave front set from
\cite{Vasy:Propagation-Many}, and we analyze the characteristic variety
of many-body Hamiltonians. Microlocal elliptic regularity is stated
here in Corollary~\ref{cor:micro-ell-reg}.
In Section~\ref{sec:br-bichar} we describe
generalized broken bicharacteristics, and in Section~\ref{sec:commutators}
we explain the positive commutator argument that is the key to our
propagation results. Propagation of singularities itself, stated in
Theorem~\ref{thm:prop-sing}, is proved
in Section~\ref{sec:propagation}.
Sections~\ref{sec:resolvent}-\ref{sec:S-matrix} turn this result into
theorems on the resolvent and the scattering matrices, stated in
Theorems~\ref{thm:res-WF}-\ref{thm:sc-matrix-disc}.
Finally, in Appendix~\ref{app:geometric} we show that if the set of
thresholds of $H$ is discrete, or if $H$ is a four-body Hamiltonian,
generalized broken bicharacteristics are piecewise integral curves
of the Hamilton vector fields, with only a finite number of breaks,
and then in Appendix~\ref{app:Lagrangian} we analyze the Lagrangian
structure of the broken bicharacteristic relation.

Most results of this paper were announced in \cite{Vasy:SJM-1999}; this
paper contains the detailed proofs.

I am very grateful to Andrew Hassell,
Richard Melrose and Maciej Zworski for numerous
very fruitful discussions; in particular, I would like to thank Richard
Melrose for his comments on this paper.
I am grateful to Maciej Zworski
for introducing me to the work of Gilles Lebeau \cite{Lebeau:Propagation}.
If there are no bound states in any subsystems, many-body scattering is
`philosophically' and, to a certain extent, technically (e.g.\ the structure
of generalized broken bicharacteristics) is very similar to the wave
equation in domain with corners. Thus, Lebeau's paper played an important
direct
role in my paper \cite{Vasy:Propagation-Many}, and remained philosophically
important while working on the present manuscript.
My joint projects with Andrew Hassell, as well as our discussions
in general, provided very valuable insights
into the broken bicharacteristic geometry, especially towards understanding
their Lagrangian structure.
I would also like
to thank Rafe Mazzeo, Erik Skibsted
and Jared Wunsch for helpful
discussions, their encouragement and for their interest in this research.

\section{Notation and detailed statement of results}\label{sec:results}
Before we can state the
precise definitions, we need to introduce some basic (and mostly
standard) notation. We refer to \cite{Derezinski-Gerard:Scattering} for
a very detailed discussion of the setup and the basic results.
We consider the Euclidean space $\Rn$, and let
$g$ be the standard Euclidean metric
on it. We assume also that we
are given a (finite) family $\calX$ of linear subspaces
$X_a$, $a\in I$,
of $\Rn$ which is closed under intersections and includes the subspace
$X_1=\{0\}$ consisting of the origin, and the whole space $X_0=\Rn$.
Let $X^a$ be the
orthocomplement of $X_a$.
We write $g_a$ and $g^a$ for the induced metrics on $X_a$ and $X^a$
respectively.
We let $\pi^a$ be the orthogonal
projection to $X^a$, $\pi_a$ to $X_a$. A many-body Hamiltonian is an operator
of the form
\begin{equation}
H=\Delta+\sum_{a\in I} (\pi^a)^*V_a;
\end{equation}
here $\Delta$ is the positive Laplacian, $V_0=0$,
and the $V_a$ are real-valued
functions in an appropriate class which we take here to be
polyhomogeneous symbols of
order $-1$ on the vector space $X_a$
to simplify the problem:
\begin{equation}
V_a\in S^{-1}_\phg(X^a).
\end{equation}
In particular, smooth potentials $V_a$
which behave at infinity like the Coulomb
potential are allowed.
Since $(\pi^a)^*V_a$ is bounded and self-adjoint and $\Delta$ is self-adjoint
with domain $H^2(\Rn)$ on $L^2=L^2(\Rn)$, $H$ is also a self-adjoint operator
on $L^2$ with domain $H^2(\Rn)$.
We let $R(\lambda)=(H-\lambda)^{-1}$ for $\lambda
\in\Cx\setminus\Real$ be the resolvent of $H$.

There is a natural partial
ordering on $I$ induced by the ordering of $X^a$
by inclusion. (Though the ordering based on inclusion of the $X_a$
would be sometimes more natural, here
we use the conventional ordering, we simply write $X_a\subset X_b$ if the
opposite ordering is required.) 
Let $I_1=\{1\}$ (recall that
$X_1=\{0\}$); $1$ is the maximal element of $I$. A maximal
element of $I\setminus I_1$ is called a 2-cluster; $I_2$ denotes the set of
2-clusters. In general, once $I_k$ has been defined for
$k=1,\ldots,m-1$, we let $I_m$ (the set of $m$-clusters)
be the set of maximal elements of
$I'_{m}=I\setminus\cup_{k=1}^{m-1} I_k$, if $I'_{m}$ is not empty.
If $I'_m=\{0\}$ (so $I'_{m+1}$ is empty), we call $H$ an $m$-body
Hamiltonian. For example, if $I\neq\{0,1\}$, and for all $a,b\nin\{0,1\}$
with $a\neq b$
we have
$X_a\cap X_b=\{0\}$, then $H$ is a 3-body Hamiltonian.
The $N$-cluster of an $N$-body Hamiltonian
is also called the free cluster, since it corresponds to the particles
which are asymptotically free.

Our goal is to study generalized eigenfunctions of $H$, i.e.\ solutions
$u\in\temp(\Rn)$ of $(H-\lambda)u=0$.
Since $H-\lambda$ is an elliptic partial differential operator with smooth
coefficients, $(H-\lambda)u\in\Cinf(\Rn)$ implies that $u\in
\Cinf(\Rn)$. Thus, the place where such $u$ can have interesting behavior
is at infinity. Analysis at infinity can be viewed either as analysis
of uniform properties, or as that of
properties in the appropriate compactification
of $\Rn$. We adopt the second point of view
by compactifying $\Rn$ as in \cite{RBMSpec}. Thus, we let
\begin{equation}
\Xb=\Xb_0=\Snp
\end{equation}
to be the radial compactification of $\Rn$ (also called the
geodesic compactification) to a closed
hemisphere, i.e.\ a ball,
and $\Sn=\partial\Snp$.
Recall from
\cite{RBMSpec} that $\SP:\Rn\to\Snp$ is given by
\begin{equation}\label{eq:SP-def}
\SP(w)=(1/(1+|w|^2)^{1/2},w/(1+|w|^2)^{1/2})\in\Snp\subset\Real^{n+1},
\quad w\in\Rn.
\end{equation}
Here we use the notation $\SP$ instead of $\operatorname{SP}$, used
in \cite{RBMSpec}, to avoid confusion with the standard stereographic
projection giving a one-point compactification of $\Rn$.
We write the coordinates on $\Rn=X_a\oplus X^a$ as $(w_a,w^a)$.
We let
\begin{equation}
\Xb_a=\cl(\SP(X_a)),\quad C_a=\Xb_a\cap\partial\Snp.
\end{equation}
Hence, $C_a$ is a sphere of dimension $n_a-1$ where
$n_a=\dim X_a$. We also let
\begin{equation}
\calC=\{C_a:\ a\in I\}.
\end{equation}
Thus, $C_0=\partial\Snp=\Sn$,
and $a\leq b$ if and only if $C_b\subset C_a$.
Since throughout this paper we work in the Euclidean setting,
where the notation $X$, $X_a$, etc., has been used for the (non-compact)
vector spaces, we always use a bar, as in $\Xb$, $\Xb_a$, etc., to
denote the corresponding compact spaces. In particular, even when talking
about general compact manifolds with boundary in the following sections,
recalling the results of \cite{Vasy:Propagation-Many}, we will write them
as $\Xb$.

We note that if $a$ is a 2-cluster then
$C_a\cap C_b=\emptyset$ unless $C_a\subset C_b$, i.e.\ $b\leq a$.
We also define the `singular part' of $C_a$ as the set
\begin{equation}
C_{a,\sing}=\cup_{b\not\leq a}(C_b\cap C_a),
\end{equation}
and its `regular part' as the set
\begin{equation}
C'_a=C_a\setminus\cup_{b\not\leq a} C_b=C_a\setminus C_{a,\sing}.
\end{equation}
For example, if $a$ is a 2-cluster then $C_{a,\sing}=\emptyset$ and
$C'_a=C_a$.

We usually identify (the interior of)
$\Snp$ with $\Rn$. A particularly useful boundary defining function of
$\Snp$ is given by $x\in\Cinf(\Snp)$ defined as
$x=r^{-1}=|w|^{-1}$ (for $r\geq 1$, say, smoothed out near
the origin); so $\Sn=\partial\Snp$ is given by $x=0$, $x>0$ elsewhere,
and $dx\neq 0$ at $\Sn$.
We write $S^{m}_\phg(\Snp)$ and $S^{m}
_\phg(\Rn)$
interchangeably.
We also remark that
\begin{equation}
S^m_{\phg}(\Snp)=x^{-m}\Cinf(\Snp).
\end{equation}
We recall that under $\SP$, $\dCinf(\Snp)$, the space of smooth functions
on $\Snp$ vanishing to infinite order at the boundary, corresponds
to the space of Schwartz functions $\Sch(\Rn)$, and its dual, $\dist(\Snp)$,
to tempered distributions $\temp(\Rn)$.
We also have the following correspondence of weighted
Sobolev spaces
\begin{equation}\label{eq:Sob-def}
\Hsc^{k,l}(\Snp)=H^{k,l}=H^{k,l}(\Rn)=\langle w\rangle^{-l}H^k(\Rn)
\end{equation}
where $\langle w\rangle=(1+|w|^2)^{1/2}$.

Corresponding to each cluster $a$ we introduce the cluster
Hamiltonian $H^a$ as an operator on $L^2(X^a)$ given by
\begin{equation}
H^a=\Delta+\sum_{b\leq a} V_b,
\end{equation}
$\Delta$ being the Laplacian of the induced metric on $X^a$.
Thus, if $H$ is a $N$-body Hamiltonian and $a$ is a $k$-cluster,
then $H^a$ is a $(N+1-k)$-body Hamiltonian. The $L^2$
eigenfunctions
of $H^a$ (also called bound states)
play an important role in many-body scattering; we
remark that by a result of Froese and Herbst, \cite{FroExp},
$\pspec(H^a)\subset(-\infty,0]$ (there are no positive
eigenvalues). Moreover, $\pspec(H^a)$ is bounded below since
$H^a$ differs from $\Delta$ by a bounded operator. Note that
$X^0=\{0\}$, $H^0=0$, so the unique eigenvalue of $H^0$ is $0$.

The eigenvalues of $H^a$ can be used to define the set of
thresholds of $H^b$. Namely, we let
\begin{equation}
\Lambda_a=\cup_{b<a}\pspec(H^b)
\end{equation}
be the set of thresholds of $H^a$, and we also let
\begin{equation}
\Lambda'_a=\Lambda_a\cup\pspec(H^a)
=\cup_{b\leq a}\pspec(H^b).
\end{equation}
Thus, $0\in\Lambda_a$ for $a\neq 0$ and $\Lambda_a\subset(-\infty,0]$.
It follows from the Mourre
theory (see e.g.\ \cite{FroMourre, Perry-Sigal-Simon:Spectral})
that $\Lambda_a$ is closed,
countable, and $\pspec(H^a)$ can only accumulate at $\Lambda
_a$. Moreover, $R(\lambda)$, considered as an operator on weighted
Sobolev spaces, has a limit
\begin{equation}
R(\lambda\pm i0):\Hsc^{k,l}(\Snp)\to \Hsc^{k+2,l'}(\Snp)
\end{equation}
for $l>1/2$, $l'<-1/2$,
from either half of the complex plane away from
\begin{equation}
\Lambda=\Lambda_1\cup\pspec(H).
\end{equation}
In addition, $L^2$ eigenfunctions
of $H^a$ with eigenvalues which are not thresholds are necessarily
Schwartz functions on $X^a$ (in fact, they decay exponentially,
see \cite{FroExp}). We also label the
eigenvalues of $H^a$, counted with multiplicities, by integers $m$,
and we call the pairs $\alpha=(a,m)$ channels. We denote the eigenvalue
of the channel $\alpha$ by $\epsilon_\alpha$, write $\psi_\alpha$ for
a corresponding normalized eigenfunction, and let $e_\alpha$ be the
orthogonal projection to $\psi_\alpha$ in $L^2(X^a)$.

The phase space in scattering theory is the cotangent bundle
$T^*\Rn$. Again, it is convenient to consider its appropriate partial
compactification, i.e.\ to consider it as a vector bundle over $\Snp$.
Thus, consider the set of all one-forms on $\Rn$ of the form
\begin{equation}
\sum_{j=1}^n a_j\, dw_j
\end{equation}
where $a_j\in\Cinf(\Snp)$ (we drop $\SP$ from the notation as usual).
This is then the set of all smooth sections of a trivial vector bundle
over $\Snp$, with basis $dw_1,\ldots,dw_n$.
Following Melrose's geometric approach to scattering theory,
see \cite{RBMSpec}, we consider this as the (dual) structure bundle,
and call it the scattering cotangent bundle of $\Snp$, denoted by $\sct\Snp$.
Note that $T^*\Rn$ can be identified with $\Rn\times\Rn$ via the metric $g$;
correspondingly $\sct\Snp$ is identified with $\Snp\times\Rn$, i.e.\ we
simply compactified the base of the standard cotangent bundle.
We remark that the construction of $\sct\Snp$
is completely natural and geometric, just like the following ones,
see \cite{RBMSpec}, or Section~\ref{sec:sc-geometry} for a summary.

However, in many-body scattering $\sct\Snp$ is
{\em not} the natural place for microlocal analysis for the very same
reason that introduces the compressed cotangent bundle
in the study of the wave equation on bounded domains.
We can see what causes trouble from both the dynamical and
the quantum point of view. Regarding dynamics,
the issue is that only the external part of the momentum
is preserved in a collision, the internal part is not; while from
the quantum point of view the problem is that there is only partial
commutativity in the algebra of
the associated pseudo-differential operators, even to top
order. To rectify this, we replace the full bundle
$\sct_{C'_a}\Snp=C'_a\times\Rn$
over $C'_a\subset\Sn$ 
by $\sct_{C'_a}\Xb_a=C'_a\times X_a$, i.e.\ we consider
\begin{equation}
\scdt \Snp=\cup_a \sct_{C'_a}\Xb_a.
\end{equation}
Over $C'_a$, there is a natural projection $\pi_a:\sct_{C'_a}\Snp
\to\sct_{C'_a}\Xb_a$ corresponding to the pull-back of one-forms;
in the trivialization given by the metric it is
induced by the orthogonal projection to $X_a$ in
the fibers. By putting the $\pi_a$ together, we obtain a projection
$\pi:\sct_{\Sn}\Snp\to\scdt\Snp$. We put the topology induced by
$\pi$ on $\scdt\Snp$. This definition is analogous to that of the compressed
cotangent bundle in the works of Melrose, Sj\"ostrand
\cite{Melrose-Sjostrand:I} and Lebeau \cite{Lebeau:Propagation} on
the wave equation in domains with smooth boundaries or corners,
respectively.

We also recall from \cite{RBMSpec} that the characteristic variety
$\Sigma_0(\lambda)$
of $\Delta-\lambda$ is simply the subset of $\sct_{\Sn}\Snp$ where
$g-\lambda$ vanishes; $g$ being the metric function. If $\Lambda_1=\{0\}$,
the compressed characteristic set of $H-\lambda$ will be simply
$\pi(\Sigma_0(\lambda))\subset\scdt\Snp$. In general, all the bound
states contribute to the characteristic variety. Thus, we let
\begin{equation}
\Sigma_b(\lambda)=\{\xi_b\in\sct_{C_b}\Xb_b:
\ \lambda-|\xi_b|^2\in\pspec{H^b}\}\subset\sct_{C_b}\Xb_b;
\end{equation}
note that $|\xi_b|^2$ is the kinetic energy of a particle in a bound
state of $H^b$.
If $C_a\subset C_b$, there is also a natural projection
$\pi_{ba}:\sct_{C'_a}\Xb_b\to\sct_{C'_a}\Xb_a$
(in the metric trivialization we can use the orthogonal projection
$X_b\to X_a$ as above), and then we define
the characteristic set of $H-\lambda$ to be
\begin{equation}
\dot\Sigma(\lambda)=\cup_a\dot\Sigma_a(\lambda),\qquad\dot\Sigma_a(\lambda)
=\cup_{C_b\supset C_a}
\pi_{ba}(\Sigma_b(\lambda))\cap\sct_{C'_a}\Xb_a,
\end{equation}
so $\dot\Sigma(\lambda)\subset\scdt \Snp$. We let $\pih_b$ be the restriction
of $\pi_b:\sct_{C'_a}\Xb_b\to\dot\Sigma(\lambda)$ to $\Sigma_b(\lambda)$.

We next recall from \cite{Vasy:Propagation-Many} the definition of generalized
broken bicharacteristics in case there are
no bound states in any of the subsystems. In fact, in this case the
word `generalized' can be dropped; for the generalized broken
bicharacteristics have a simple geometry as stated below.
First, note that
the rescaled Hamilton vector field of the metric function $g$,
i.e.
\begin{equation}
2\langle w\rangle\xi\cdot\partial_w\in\Vf(T^*\Rn)=\Vf(\Rn\times\Rn)
\end{equation}
extends to a smooth vector field, $\scHg\in\Vf(\sct\Snp)=\Vf(\Snp\times\Rn)$,
with $\Snp$ considered as the radial compactification of $\Rn$; in
fact, $\scHg$ is tangent to the boundary $\sct_{\Sn}\Snp=\Sn\times\Rn$.

\begin{Def*}
Suppose $\Lambda_1=\{0\}$, and $I=[\alpha,\beta]$ is an interval.
We say that a continuous map $\gamma:I\to\dot\Sigma(\lambda)$
is a broken bicharacteristic of $H-\lambda$ if
there exists a finite
set of points $t_j\in I$, $\alpha=t_0<t_1<\ldots<t_{k-1}<t_k=\beta$
such that for each $j$,
$\gamma|_{(t_j,t_{j+1})}$ is the image of an integral curve of
$\scHg$ in $\Sigma_0(\lambda)$ under $\pi$.
If $I$ is an interval (possibly $\Real$), we say that $\gamma:I\to
\dot\Sigma(\lambda)$ is a broken bicharacteristic of $H-\lambda$, if the
restriction of $\gamma$ to every compact subinterval of $I$ is a
broken bicharacteristic in the above sense.
\end{Def*}

Here $\gamma(I)\subset\dot\Sigma(\lambda)=\pi(\Sigma_0(\lambda))$
corresponds to the conservation of kinetic energy in
collisions (since there are no bound states), and
the use of the compressed space $\scdt\Snp$
shows that external momentum is conserved
in the collisions. It turns out, see \cite{Vasy:Propagation-Many}, that
$\gamma$ is essentially the lift of a broken geodesic on $\Snp$ of length
$\pi$ (if $I=\Real$, otherwise shorter),
broken at the collision planes, i.e.\ at $\calC$. In particular, even if
$I=\Real$, it has only a finite number of breaks, and in fact, there is a
uniform bound on the number of such breaks (depending only on the geometry,
i.e.\ on $\calC$, not on $\gamma$).

The definitions are less explicit if $\Lambda_1\neq\{0\}$,
but they essentially still state that the total energy and the external
momentum are preserved in collisions. Thus, generalized
broken bicharacteristics
will be continuous maps $\gamma$ defined
on intervals $I$, $\gamma:I\to\dot\Sigma(\lambda)$ with certain appropriate
generalization of the integral curve condition described above.
In order to take the bound states into consideration, we also need
to consider the rescaled Hamilton vector fields $\scHg^b$ of the metric $g_b$
in the subsystem $b$. Thus, under the inclusion map
\begin{equation}\label{eq:imath_b-def}
\imath_b:\sct_{C_b}\Xb_b\hookrightarrow\sct_{C_b} \Snp
\end{equation}
induced by the inclusion $X_b\hookrightarrow\Rn$
in the fibers, $(\imath_b)_*\scHg^b=\scHg$ (i.e.\ the restriction
of the vector field $\scHg$ to $\sct_{C_b} \Xb_b$, considered as a subset
of $\sct_{C_b} \Snp$).
Thus,
we require that lower bounds on the Hamilton vector fields $\scHg^b$
applied to $\pi$-invariant
functions, i.e.\ to functions $f\in\Cinf(\sct_{\Sn}\Snp)$ such that
$f(\xi)=f(\xi')$ if $\pi(\xi)=\pi(\xi')$, imply lower bounds on the
derivatives of $f_\pi$ along $\gamma$. Here $f_\pi$ is the function
induced by $f$ on $\scdt\Snp$, so $f=f_\pi\circ\pi$.

\begin{Def}\label{Def:gen-br-bichar}
A generalized broken bicharacteristic of $H-\lambda$
is a continuous map
$\gamma:I\to\dot\Sigma(\lambda)$,
where $I\subset\Real$ is an interval, such that for all $t_0\in I$
and for each sign $+$ and $-$ the following holds. Let
$\xi_0=\gamma(t_0)$, suppose that $\xi_0\in\sct_{C'_a}\Xb_a$. Then for
all $\pi$-invariant functions
$f\in\Cinf(\sct_{\Sn}\Snp)$,
\begin{equation}
D_\pm(f_\pi\circ \gamma)(t_0)
\geq\inf\{\scHg^b f(\xit_0):\ \xit_0\in\pih_{b}^{-1}(\xi_0),\ C_a\subset C_b\}.
\end{equation}
Here $D_\pm$ are the one-sided lower derivatives: if $g$ is defined on
an interval $I$, $(D_\pm g)(t_0)=\liminf_{t\to t_0\pm} (g(t)-g(t_0))/(t-t_0)$.
\end{Def}

Although it is not apparent, this definition is equivalent to the previous
one if $\Lambda_1=\{0\}$. Moreover, in four-body scattering, even if
$\Lambda_1\neq\{0\}$, one can describe
the generalized
broken bicharacteristics piecewise as projections of integral curves
of $\scHg$. In general many-body scattering, the lack of conservation of
kinetic energy makes such a description harder, but if $\Lambda_1$ is
discrete, we obtain a description that parallels the one above.
More precisely, in Theorem~\ref{thm:geom-br-bichar} we prove the
following. Suppose that $\Lambda_1$ is discrete and
$\gamma:\Real\to\dot\Sigma(\lambda)$
is a continuous curve. Then $\gamma$
is a generalized
broken bicharacteristic of $H-\lambda$ if and only if there exist
$t_0<t_1<t_2<\ldots<t_k$ such that $\gamma|_{[t_j,t_{j+1}]}$,
as well as $\gamma|_{(-\infty,t_0]}$ and $\gamma|_{[t_k,+\infty)}$, are
the projections of integral curves
of the Hamilton vector field $\scHg^a$ for some $a$. In addition,
there is a uniform bound on $k$ (independent of $\gamma$),
depending only on $\calC$ and $\Lambda_1$. Similar results
hold if the interval of definition, $\Real$, is replaced by any interval.

As mentioned in the introduction, `singularities' (i.e.\ lack of
decay at infinity) of $u\in\temp$ are described
by the many-body scattering wave front set, $\WFSc(u)$, which was
introduced in \cite{Vasy:Propagation-Many}, and which describes $u$ modulo
Schwartz functions, similarly to how the usual wave front set
describes distributions modulo smooth functions. Just as for the image
of the bicharacteristics, $\scdt\Snp$ provides the natural setting
in which $\WFSc$ is defined: $\WFSc(u)$ is a closed subset of $\scdt\Snp$.
The definition of $\WFSc(u)$ relies on the algebra of many-body scattering
pseudo-differential operators, also introduced in
\cite{Vasy:Propagation-Many}. There are several possible definitions of
$\WFSc$,
all of which agree for generalized eigenfunctions of $H$, but the
one given in \cite{Vasy:Propagation-Many} that is modelled on the
fibred-cusp wave front set of Mazzeo and Melrose \cite{Mazzeo-Melrose:Fibred}
enjoys many properties of the usual wave front set. A slightly
different definition, for which the crucial property in the last line
of \eqref{eq:WF-props}
still holds, was discussed in \cite{Vasy:Propagation-2}; this is the only
essential property for the positive commutator estimates in this paper.
Also, the discussion after \eqref{eq:Pt-error} indicates why
`finite order' versions of the wave front set,
i.e.\ versions in which we only require a
fixed (though possibly high) number of the symbol estimates for ps.d.o.'s,
somewhat akin to the constructions of \cite{GerComm},
would be helpful; nonetheless, the Mazzeo-Melrose definition appears to be
the most natural one from the point of view of general microlocal
analysis.

We recall the precise definitions in Section~\ref{sec:indicial};
here we also
translate our results into statements on the S-matrices where the usual
wave front set can be used.
We remark that in
the two-body setting, when $\scdt\Snp=\sct\Snp$, $\WFSc$ is just the
scattering wave front set $\WFsc$ introduced by Melrose, \cite{RBMSpec},
which in turn is closely related to the usual wave front set via the Fourier
transform. Thus, for $(\omega,\xi)\in\sct_{\Sn}\Snp$, considered as
$\Sn\times\Rn=\partial\Snp\times\Rn$,
$(\omega,\xi)\nin\WFsc(u)$ means that there exists
$\phi\in\Cinf(\Snp)$ such that $\phi(\omega)\neq 0$ and $\Fr(\phi u)$ is
$\Cinf$ near $\xi$. If we employed the usual conic terminology instead
of the compactified one, we would think of $\phi$ as a conic cut-off
function in the direction $\omega$. Thus, $\WFsc$ at infinity is
analogous to $\WF$ with the role of position and momentum reversed.
The definition of $\WFSc(u)$ is more complicated, but if
$u=\psi(H)v$ for some $\psi\in\Cinf_c(\Real)$ (any other operator
in $\PsiSc^{-\infty,0}(\Xb,\calC)$ would do instead of $\psi(H)$), then
the following is
a sufficient condition for $(\omega,\xi_a)\in\sct_{C'_a}\Xb_a$,
considered as $C'_a\times X_a$, not
to be in $\WFSc(u)$. Suppose that there exists
$\phi\in\Cinf(\Snp)$, $\phi(\omega)\neq 0$, and $\rho\in\Cinf_c(X_a)$,
$\rho(\xi_a)\neq 0$, and $((\pi^a)^*\rho)\Fr(\phi u)\in\Sch(\Rn)=\Sch(X_0)$.
Then $(\omega,\xi_a)\nin\WFSc(u)$.
We also remark that
we state all of the following results for the absolute wave front sets
(i.e.\ we work modulo Schwartz functions), but they have complete analogues
for the relative wave front sets (working modulo weighted Sobolev spaces);
indeed, it is the latter that is used to prove the results on the former.

Our main result is then the following theorem, in which we allow
arbitrary thresholds, and which describes the relationship between
$\WFSc(u)$ and generalized
broken bicharacteristics, if, for example, $(H-\lambda)u=0$. Note that
if $(H-\lambda)u=0$, then $u=\psi(H)u$ for $\psi\in\Cinf_c(\Real)$,
$\psi\equiv 1$ near $\lambda$, so the above description of $\WFSc(u)$ is
applicable.

\begin{thm}\label{thm:prop-sing}
Let $u\in\temp(\Rn)$,
$\lambda\nin\Lambda_1$. Then
\begin{equation}
\WFSc(u)\setminus\WFSc((H-\lambda)u)
\end{equation}
is a union of maximally extended generalized
broken bicharacteristics of $H-\lambda$ in
$\dot\Sigma(\lambda)\setminus\WFSc((H-\lambda)u)$.
\end{thm}

We remark that the statement of the theorem is empty at points 
$\xi_0\in\sct_{C'_a}\Xb_a$ at which $\scHg^b (\xit_0)=0$ for
some $\xit_0\in\pih_b^{-1}(\xi_0)$ and some $b$ with $C_a\subset C_b$.
Indeed, at such points the
constant curve ($\gamma(t)=\xit_0$ for all $t$ in some interval)
is a generalized broken bicharacteristic. A simple calculation shows
that the set of these points $\xi_0$ is $R_+(\lambda)\cup
R_-(\lambda)$, where
\begin{equation}
R_\pm(\lambda)=\{\xi\in\sct_{C'_a}\Xb_a:\ \exists b,\ C_a\subset C_b,
\lambda-\tau(\xi)^2\in\pspec(H^b),\ \pm\tau(\xi)\geq 0\}
\end{equation}
are the incoming ($+$) and outgoing ($-$) radial sets respectively, and
$\tau$ is the $\scl$-dual variable of the boundary defining function $x$,
so in terms of the Euclidean variables
\begin{equation}
\tau=-\frac{w\cdot\xi}{|w|};
\end{equation}
see the next section for further details.
Hence, the theorem permits
singularities to emerge `out of nowhere' at the radial sets. Although we do
not prove that this indeed does happen, based on general principles, this
appears fairly likely. Moreover, the optimality of Theorem~\ref{thm:prop-sing}
if $\Lambda_1=\{0\}$ follows from \cite{Hassell:Plane, Vasy:Structure},
see the remarks
about this in \cite{Vasy:Propagation-Many}; the amplitude of the reflected
`wave' is given (to top order) by the appropriate subsystem S-matrix.

There is a similar result for $\WFSc(u)$, $u=R(\lambda+i0)f$; namely that
$\WFSc(u)\setminus\WFSc(f)$ is the image of $\WFSc(f)\cup R_-(\lambda)$
under forward propagation,
if e.g.\ $f\in H^{r,s}$, $s>1/2$. The set $R_-(\lambda)$ appears here since
there can be maximally extended generalized broken bicharacteristics
which are either not disjoint from $R_-(\lambda)$, or simply whose closure
is not disjoint from $R_-(\lambda)$.
In particular, even if $f$ is Schwartz, $\WFSc(u)$ is not necessarily
a subset of $R_-(\lambda)$, rather a subset of its image under forward
propagation. Indeed, by duality, this is exactly what gives rise to
the conditions on $\WFSc(f)$ under which $u=R(\lambda+i0)f$ can be defined.
To make it easier to state these results, we make the following definition.

\begin{Def}
Suppose $K\subset\dot\Sigma(\lambda)$. The image $\Phi_+(K)$ of $K$ under the
forward broken bicharacteristic relation is defined as
\begin{equation}\begin{split}
\Phi_+(K)=&\{\xi_0\in\dot\Sigma(\lambda):\ \exists
\ \text{a generalized broken bicharacteristic}\ \gamma:(-\infty,t_0]\to
\dot\Sigma(\lambda)\\
&\ \text{s.t.}
\ \gamma(t_0)=\xi_0,\ \overline{\gamma((-\infty,t_0])}\cap K\neq\emptyset\}.
\end{split}\end{equation}
The image $\Phi_-(K)$ of $K$ under the backward broken bicharacteristic
relation
is defined similarly, with $[t_0,+\infty)$ in place of $(-\infty,t_0]$.
\end{Def}

Note that $\Phi_+(K)=\cup_{\xi\in K}\Phi_+(\{\xi\})$ directly from the
definition.
The result on the boundary values of the resolvent is then:

\begin{thm}\label{thm:res-WF}
Suppose that $\lambda\nin\Lambda$, $f\in\Sch(\Rn)$, and let
$u=R(\lambda+i0)f$. Then $\WFSc(u)\subset
\Phi_+(R_-(\lambda))$. Moreover, $R(\lambda+i0)$ extends by continuity
to $v\in\temp(\Rn)$ with $\WFSc(v)\cap\Phi_-(R_+(\lambda))=\emptyset$,
and for such $v$,
\begin{equation}
\WFSc(R(\lambda+i0)v)\subset\Phi_+(\WFSc(v))\cup\Phi_+(R_-(\lambda)).
\end{equation}
\end{thm}

The scattering matrices $S_{\beta\alpha}(\lambda)$
of $H$ with incoming channel $\alpha$, outgoing channel $\beta$ can be defined
either via the wave operators, or via the asymptotic behavior of generalized
eigenfunctions. It was shown in \cite{Vasy:Scattering} that the two are
the same, up to normalization (free motion is factored out in the wave
operator definition); here we briefly recall the second definition.
We first state it for short-range $V_c$ ($V_c$ polyhomogeneous of order
$-2$ for all $c$).
Thus, for $\lambda\in(\epsilon_\alpha,\infty)\setminus\Lambda$
and $g\in\Cinf_c(C_a')$,
there is a unique $u\in\temp(\Rn)$
such that $(H-\lambda)u=0$, and $u$ has the form
\begin{equation}\label{eq:intro-21}
u=e^{-i\sqrt{\lambda-\ep_\alpha} r}r^{-\dim C_a/2}((\pi^a)^*\psi_\alpha)v_-
+R(\lambda+i0)f,
\end{equation}
where $v_-\in\Cinf(\Snp)$, $v_-|_{C_a}=g$, and
$f\in\Hsc^{\infty,1/2+\epsilon'}(\Snp)$, $\epsilon'>0$.
The Poisson operator $P_{\alpha,+}(\lambda)$ is the map
\begin{equation}\label{eq:Poisson-def}
P_{\alpha,+}(\lambda):\Cinf_c(C_a')\to\temp(\Rn)\ \text{defined by}
\ P_{\alpha,+}(\lambda)g=u.
\end{equation}
The term $R(\lambda+i0)f$ has distributional asymptotics of a similar form
`at the channel $\beta$', i.e.\ of the form
$e^{i\sqrt{\lambda-\ep_\beta} r}r^{-\dim C_b/2}
((\pi^b)^*\psi_\beta)v_{\beta,+}$,
see \cite{Vasy:Scattering} for the precise definitions.

Only minor modifications are necessary for $V_c\in S^{-1}_\phg(X^c)$.
Namely, write
\begin{equation}
I_a=\sum_{b\not\leq a}V_b,\ \It_a=(r_a I_a)|_{C_a}\in\Cinf
(C_a),
\end{equation}
$I_a$ is $\Cinf$ near $C'_a$ with simple vanishing at $C'_a$ (since
$b\not\leq a$ means
$C_a\not\subset C_b$, hence $C_a\cap C_b\subset C_{a,\sing}$), so
$r_a I_a$ is $\Cinf$ there. Then the asymptotics in \eqref{eq:intro-21}
must be replaced by
\begin{equation}
e^{-i\sqrt{\lambda-\ep_\beta} r}r^{-\dim C_a/2}
r^{i\It_a/2\sqrt{\lambda-\ep_\alpha}}((\pi^a)^*\psi_\alpha)v_-
\end{equation}
plus lower order terms; see Section~\ref{sec:Poisson} for details.

The scattering matrix $S_{\beta\alpha}(\lambda)$
maps $g=v_-|_{C_a}$ to $v_{\beta,+}|_{C'_b}$.
It is also given by the formula
\begin{equation}\label{eq:equiv-7}
S_{\beta\alpha}(\lambda)=\frac{1}{2i\sqrt{\lambda-\ep_\beta}}
((H-\lambda)\Pt_{\beta,-}(\lambda))^*P_{\alpha,+}(\lambda),
\end{equation}
$\lambda>\max(\epsilon_\alpha,\epsilon_\beta)$, $\lambda\nin\Lambda$.
Here $\Pt_{\beta,-}(\lambda)$ is a microlocalized version of the
outgoing Poisson operator, microlocalized near the outgoing region
for $\beta$, i.e.\ where $\tau$ is near $-\sqrt{\lambda-\ep_\beta}$, see
\cite{Vasy:Scattering}.
In fact, we can simply take $\Pt_{\beta,-}(\lambda)$
to be a microlocal (cut-off) parametrix for $P_{\beta,-}(\lambda)$.
This formula is closely related to that of Isozaki and Kitada
\cite{Isozaki-Kitada:Scattering}.

A very good parametrix, $\Pt_{\alpha,+}(\lambda)$,
for $P_{\alpha,+}(\lambda)$ in the region
of phase space where $\tau$ is close to $\sqrt{\lambda-\ep_\alpha}$ has
been constructed by Skibsted \cite{Skibsted:Smoothness} in the short-range
and by Bommier \cite{Bommier:Proprietes} in the long-range setting,
under the assumption that $\ep_\alpha\in\spec_d(H^a)$. If we instead
assume that the $V_c$ are all Schwartz, then the trivial (product
type) construction gives the desired parametrix.
Their constructions enable us to deduce the structure of the S-matrices
immediately from our propagation theorem,
Theorem~\ref{thm:res-WF} via \eqref{eq:equiv-7} and
\begin{equation}\label{eq:Poisson-intro-par}
P_{\alpha,+}(\lambda)=\Pt_{\alpha,+}(\lambda)-R(\lambda+i0)(H-\lambda)
\Pt_{\alpha,+}(\lambda).
\end{equation}
Since the parametrix (near the incoming or outgoing sets)
is important for turning the results on the propagation
of singularities to wave front set results, in all our results on the
Poisson operators and S-matrices $S_{\beta\alpha}(\lambda)$ in this paper
we make the following assumption:
\begin{equation}\label{eq:param-Poisson}
\text{either}\ \ep_\alpha\in\spec_d(H^a)\Mand\ep_\beta\in\spec_d(H^b),
\Mor V_c\in\Sch(X^c)\ \text{for all}\ c.
\end{equation}

It is easy to describe the wave front set of $\Pt_{\alpha,+}(\lambda)g$,
$g\in\dist_c(C'_a)$, near its `beginning point', i.e.\ near the
$(\alpha,+)$-incoming set. Namely, it is the union
of integral curves of $\scHg^a$ (in $\sct_{C'_a}\Xb_a$) (which are in
particular bicharacteristics of $H-\lambda$, hence broken bicharacteristics),
one integral curve for each $\zeta\in\WF(g)\subset S^*C'_a$; we denote
these integral curves by $\gamma_{\alpha,-}(\zeta)$. It is actually
convenient to replace the parameter $t$ of the integral curve by $s$, the
arclength parameter of its projection to $C_a$. The relationship between
these two is that if we write $s=S(t)$, then $S$ solves the ODE
$dS/dt=2(\lambda-\ep_\alpha-\tau(\gamma(t))^2)^{1/2}$. The reparameterized
integral curves are then given by
\begin{equation}
\tau_a=\sqrt{\lambda-\ep_\alpha}\cos(s-s_0),
\ (y_a,\mu_a)=
\sqrt{\lambda-\ep_\alpha}\sin(s-s_0)\exp((s-s_0)\scHg^a)(\zeta)
\end{equation}
where $s\in(s_0,s_0+\pi)$.
This defines $\gamma_{\alpha,-}(\zeta)$ up to replacing $t$ by $t-t_1$ for any
fixed $t_1\in\Real$, so we are abusing the notation slightly.
Due to \eqref{eq:Poisson-intro-par},
Theorem~\ref{thm:res-WF} describes $\WFSc(P_{\alpha,+}(\lambda)g)$ elsewhere.
A similar result also applies for $\Pt_{\beta,-}(\lambda)$; in this case
one simply has to replace the range of the arclength parameter by
$s\in(s_0-\pi,s_0)$. We denote the corresponding integral curves by
$\gamma_{\beta,+}(\zeta)$.

\begin{Def}
The forward broken bicharacteristic relation with initial channel $\alpha$
is defined to be the relation $\calR_{\alpha,+}\subset 
\dot\Sigma(\lambda)\times S^*C'_a$ given by
\begin{equation}
\calR_{\alpha,+}=\{(\zeta,\xi)\in \dot\Sigma(\lambda)\times S^*C'_a:
\ \Phi_-(\{\xi\})\cap\gamma_{\alpha,-}(\zeta)\neq\emptyset\}.
\end{equation}
The backward broken bicharacteristic relation with initial channel $\beta$,
denoted by $\calR_{\beta,-}$ is defined similarly, with $\Phi_+$ in place
of $\Phi_-$ and $\gamma_{\beta,+}$ in place of $\gamma_{\alpha,-}$.
Finally, the forward broken bicharacteristic relation with initial channel
$\alpha$, final channel $\beta$, $\calR_{\beta\alpha}\subset S^*C'_b
\times S^*C'_a$ is defined as the composite relation of $\calR_{\alpha,+}$
and $\calR_{\beta,-}^{-1}$:
\begin{equation}
\calR_{\beta\alpha}
=\{(\zeta,\zeta')\in S^*C'_b\times S^*C'_a:\ \exists \xi\in\dot\Sigma(\lambda),
(\zeta,\xi)\in\calR_{\alpha,+},\ (\zeta',\xi)\in\calR_{\beta,-}\}.
\end{equation}
\end{Def}

Note that $(\zeta,\xi)\in\calR_{\alpha,+}$ thus means that there
exists a generalized broken bicharacteristic
$\gamma:\Real\to\dot\Sigma(\lambda)$ and $t_0\in\Real$ such that
$\gamma|_{(-\infty,t_0]}$ is given by $\gamma_{\alpha,-}(\zeta)$,
and $\xi\in\overline{\gamma([t_0,+\infty))}$. Thus, $\calR_{\alpha,+}$
should be thought of as the relation induced by $\Phi_-$ `at channel
$\alpha$' as time goes to $-\infty$.

If $\calR\subset B\times A$ is a relation, $K\subset A$, by $\calR(K)$
we mean $\{\xi\in B:\ \exists\zeta\in K,\ (\xi,\zeta)\in\calR\}$.
Similarly, if $U\subset B$, by $\calR^{-1}(U)$ we mean
$\{\zeta\in A:\ \exists\xi\in U,\ (\xi,\zeta)\in\calR\}$. We call $\calR(K)$
the image of $K$ under $\calR$. Thus,
if $K\subset S^*C'_a$,
\begin{equation*}
\calR_{\alpha,+}(K)=\{\xi\in\dot\Sigma(\lambda):\ \exists\zeta\in K,
\ \Phi_-(\{\xi\})\cap\gamma_{\alpha,-}(\zeta)\neq\emptyset\},
\end{equation*}
and if $U\subset\dot\Sigma(\lambda)$, then
\begin{equation*}
\calR^{-1}_{\alpha,+}(U)=\{\zeta\in S^*C'_a:\ \exists\xi\in U,
\ \Phi_-(\{\xi\})\cap\gamma_{\alpha,-}(\zeta)\neq\emptyset\}.
\end{equation*}

This definition, \eqref{eq:Poisson-intro-par} and
Theorem~\ref{thm:res-WF} immediately prove the following proposition.

\begin{prop}\label{prop:Poisson-extend}
Suppose that $H$ is a many-body Hamiltonian, $\lambda\nin\Lambda$,
and \eqref{eq:param-Poisson} holds.
Suppose also that $g\in\Cinf_c(C'_a)$. Then
\begin{equation}\label{eq:Poisson-WF-smooth}
\WFSc(P_{\alpha,+}(\lambda)g)
\setminus R_+(\lambda)\subset\Phi_+(R_-(\lambda)).
\end{equation}

In addition, $P_{\alpha,+}(\lambda)$ extends by continuity from
$\Cinf_c(C'_a)$ to distributions
$g\in\dist_c(C'_a)$ with $(R_+(\lambda)\times\WF(g))\cap\calR_{\alpha,+}
=\emptyset$ (i.e.\ $\calR_{\alpha,+}(\WF(g))\cap R_+(\lambda)=\emptyset$).
If $g$ is such a distribution, then
\begin{equation}
\WFSc(P_{\alpha,+}(\lambda)g)\setminus R_+(\lambda)
\subset\Phi_+(R_-(\lambda))\cup
\calR_{\alpha,+}(\WF(g)).
\end{equation}
\end{prop}

One of the main features of \eqref{eq:Poisson-WF-smooth} is that
in general $\WFSc(P_{\alpha,+}(\lambda)g)$ cannot be expected to be
contained in the radial sets; one also
has to include the image of the outgoing
radial set under forward propagation in the statement.
As a corollary, \eqref{eq:equiv-7} shows that in general
$S_{\beta\alpha}(\lambda)$ does not
map smooth incoming data to smooth outgoing data. However, if $\beta$
is a two-cluster channel, every generalized broken bicharacteristic
$\gamma$ such that for some $t_0\in\Real$,
$\gamma|_{(-\infty,t_0]}$ is given by $\gamma_{\beta,+}(\zeta)$,
$\zeta\in S^* C'_b$, is actually equal to $\gamma_{\beta,+}(\zeta)$
for all times, and, $b$ being a 2-cluster, $\gamma_{\beta,+}(\zeta)$
never intersects the radial sets, and as $t\to\pm\infty$,
$\gamma_{\beta,+}(\zeta)(t)$ goes to $R_\mp(\lambda)$. Thus, if
$\beta$ is a 2-cluster, $\alpha$ is any cluster, $S_{\beta\alpha}(\lambda)$
maps smooth functions to smooth functions. On the other hand,
if $\beta$ is the free channel $0$, then the absence of positive
thresholds gives a similar conclusion.

\begin{cor}
Suppose that $H$ is a many-body Hamiltonian, $\lambda\nin\Lambda$,
and \eqref{eq:param-Poisson} holds.
Suppose $g\in\Cinf_c(C'_a)$.
Then $\WF(S_{\beta\alpha}(\lambda)g)\subset\calR^{-1}_{\beta,-}(R_-(\lambda))$.
Thus, if $\beta$ is the free channel or it is
a two-cluster channel, then $S_{\beta\alpha}(\lambda)g$ is $\Cinf$.
\end{cor}

Our theorem on the wave front relation of the S-matrix is then the following.

\begin{thm}\label{thm:sc-matrix}
Suppose that $H$ is a many-body Hamiltonian, $\lambda\nin\Lambda$,
and \eqref{eq:param-Poisson} holds.
Then $S_{\beta\alpha}(\lambda)$ extends by continuity from
$\Cinf_c(C'_a)$ to distributions
$g\in\dist_c(C'_a)$ with $\calR_{\alpha,+}(\WF(g))\cap R_+(\lambda)=\emptyset$.
If $g$ is such a distribution, then
\begin{equation}
\WF(S_{\beta\alpha}(\lambda)g)\subset\calR^{-1}_{\beta,-}(R_-(\lambda))\cup
\calR_{\beta\alpha}(\WF(g)).
\end{equation}
\end{thm}

If $\Lambda_1=\{0\}$, then maximally extended
generalized broken bicharacteristics are
essentially the lift of generalized broken geodesics on $\Sn$ of length
$\pi$, so in this case we recover the following
result of \cite{Vasy:Propagation-Many}.

\begin{cor}
If no subsystem of $H$ has bound states and $\lambda>0$, then
the wave front relation of $S_{00}(\lambda)$ is given by the broken geodesic
relation on $\Sn$, broken at $\calC$, at distance $\pi$.
\end{cor}

It may seem that our results are too weak in the sense that
some broken bicharacteristics can be continued in many ways when they
hit a collision plane. However, it should be noted that not every
broken bicharacteristic hits a collision plane; indeed, if the dimension
of a collision plane increases, more broken bicharacteristics will hit it,
but each will generate a lower dimensional family of continuations.
We can make this precise in terms of a Lagrangian characterization
of the wave front relation of the S-matrices.
To do so, recall first that the wave front relation of operators mapping
$\Cinf_c(C'_a)$ to $\dist(C'_b)$ is
often understood as a subset of $T^*C'_b\times T^* C'_a\setminus 0$
corresponding to
the wave front set of its Schwartz kernel, with the sign of the second
component, the one in $T^*C'_a$, switched. For example, in this sense,
the wave front relation of an FIO, say $P$, is given by a homogeneous
canonical relation, i.e.\ by a Lagrangian submanifold $\Lambdat$ of
$T^*C'_b\times T^* C'_a$, Lagrangian with respect to the usual twisted
symplectic form $\omega_b-\omega_a$, that is conic with respect to
the diagonal $\Real^+$ action in $T^*C'_b\times T^* C'_a$. Here
$\omega_a$ denotes the canonical symplectic form on $T^*C_a$, etc.
The wave front set mapping properties of such an operator are
that $Pu$ is defined if $(0\times\WF(u))\cap\Lambdat=\emptyset$,
and if $Pu$ is thus defined,
\begin{equation}\begin{split}\label{eq:WF-map-3}
\WF(Pu)\subset &\Lambdat(\WF(u))\cup\{\xi\in T^*C'_b\setminus 0:
\ (\xi,0)\in\Lambdat\}\\
&\quad=\{\xi\in T^*C'_b\setminus 0:\ (\xi,0)\in\Lambdat\Mor\exists
\xi'\in\WF(u),\ (\xi,\xi')\in\Lambdat\}.
\end{split}\end{equation}
Note that $\WF(u)$ is a conic subset of $T^*C'_a$, so it is better to regard
it as a subset of $S^*C'_a$ (quotient out by the $\Real^+$ action on
$T^*C'_a\setminus 0$);
indeed, this is what we have done. But every element $\zeta$ of
$S^*C'_a$ (thought of as the quotient bundle) has several representatives
$\xi$ in $T^*C'_a$. Correspondingly, the wave front relation, understood
as a relation connecting $S^*C'_a$ and $S^*C'_b$,
i.e.\ as a subset of $S^*C'_b\times S^*C'_a$, only determines
the corresponding conic relation on $(T^*C'_b\setminus 0)\times
(T^*C'_a\setminus 0)$ up to the
rescaling of one factor with respect to the other. Thus, our results
by themselves cannot pinpoint the wave front set of the
Schwartz kernel of $S_{\beta\alpha}(\lambda)$
itself, only its image under the quotient map with respect to the
$\Real^+\times\Real^+$ action.
Note, however, that $\Lambdat\cap (T^*C'_b\times 0)$ and $\Lambdat
\cap (0\times T^* C'_a)$ show up in the mapping properties of $P$, namely
whether it can be applied to all distributions, and whether it maps
smooth functions to smooth functions. Thus, by Theorem~\ref{thm:sc-matrix}
we would expect, if $S_{\beta\alpha}(\lambda)$ were an FIO, that
the (twisted) wave front set of its kernel lies in
\begin{equation}\label{eq:WF-SK-16}
(\calR^{-1}_{\beta,-}(R_-(\lambda))\times 0)\cup(0\times \calR_{\alpha,+}^{-1}
(R_+(\lambda))\cup \Lambdat_0,
\end{equation}
where $\Lambdat_0\subset (T^*C'_b\setminus 0)\times (T^*C'_a\setminus 0)$
is conic Lagrangian and projects to $\calR_{\beta\alpha}$ under the quotient
map in both factors.
While, as explained, our results as stated do not
prove that the wave front set of the kernel of $S_{\beta\alpha}(\lambda)$
lies in \eqref{eq:WF-SK-16}, we can prove its `quotiented' version.

\begin{thm}\label{thm:sc-matrix-disc}
Suppose that $H$ is a many-body Hamiltonian,
$\lambda\nin\Lambda$, and \eqref{eq:param-Poisson} holds.
Suppose also that $\Lambda_1$ is discrete.
The wave front relation of $S_{\beta\alpha}(\lambda)$ is a subset of
the projection of a finite union of
conic Lagrangian submanifolds of $T^*(C'_a\times C'_b)\setminus 0$,
given by the
broken bicharacteristic relation of $H$. Some of the Lagrangians may lie
in $T^*C'_a\times 0$ or $0\times T^*C'_b$ where $0$ denotes the zero
section.
\end{thm}

To see how this theorem is proved, recall that under our
discreteness assumption the generalized
broken bicharacteristics are actually broken bicharacteristics
with a finite number of breaks (and a uniform bound on
this finite number). Correspondingly, we can associate to
each generalized broken bicharacteristic one of a finite
number of collision patterns which describe from which clusters
in which sequence the bicharacteristics reflect, and along which
collision plane at which kinetic energy they travel
meanwhile. We prove in Appendix~\ref{app:Lagrangian}
that each of these patterns corresponds to a
(typically not closed) conic Lagrangian submanifold of
$T^*(C'_b\times C'_a)\setminus 0$,
the projection of the union of which gives the
wave front relation (which is, however, closed), thereby proving this theorem.

We remark that if there are no bound
states in the proper subsystems, the proof of the result is
essentially the same as in the three-body setting
\cite{Vasy:Structure, Hassell:Plane}; Hassell's proof employing Jacobi
vector fields on the sphere is particularly easy (and transparent!) to
adapt.

The statement of this theorem shows that, under the discreteness assumption,
the wave front relation of
the S-matrices has the correct structure to be the wave front
relation of a Fourier integral operator (FIO). However, the construction
of such an FIO would require a much better understanding of the symplectic
geometry of, and related analysis on,
manifolds with corners that are equipped with certain boundary fibrations
(corresponding to the resolved space
$[\Xb;\calC]$ and the fibration given by the blow-down map which
are discussed in the next section). In relatively simple settings
this has been discussed in \cite{Mazzeo-Melrose:Fibred,
Hassell-Vasy:Spectral}, but much progress is needed
for dealing with the general case.

\section{Scattering geometry and analysis}\label{sec:sc-geometry}
Next, we recall from \cite{RBMSpec} Melrose's definition of
the Lie algebra of `scattering vector fields' $\Vsc(\Xb)$, defined for
every manifold with boundary $\Xb$. Thus,
\begin{equation}
\Vsc(\Xb)=x\Vb(\Xb)
\end{equation}
where $\Vb(\Xb)$ is the set of smooth vector fields
on $\Xb$ which are tangent to $\bXb$. If $(x,y_1,\ldots, y_{n-1})$
are coordinates on $\Xb$ where
$x$ is a boundary defining function, then locally a basis of $\Vsc(\Xb)$
is given by
\begin{equation}
x^2\partial_x,\ x\partial_{y_j},\ j=1,\ldots,n-1.
\end{equation}
Correspondingly, there is a vector bundle $\Tsc \Xb$ over $\Xb$, called
the scattering tangent bundle of $\Xb$, such that $\Vsc(\Xb)$ is the set of
all smooth sections of $\Tsc \Xb$:
\begin{equation}
\Vsc(\Xb)=\{v\in\Cinf(\Xb;\Tsc \Xb):\ \forall p\in\Xb,\ v_p\in\Tsc_p\Xb\}.
\end{equation}
The
dual bundle of $\Tsc \Xb$ (called the scattering cotangent bundle) is
denoted by $\sct \Xb$.
Thus, covectors $v\in\sct_p \Xb$, $p$ near
$\bXb$, can be written as 
\begin{equation}\label{eq:sct-coords}
v=\tau\,\frac{dx}{x^2}+\mu\cdot\frac{dy}{x}.
\end{equation}
Hence, we have local coordinates $(x,y,\tau,\mu)$ on $\sct \Xb$ near $\bXb$.
The scattering density bundle $\Osc \Xb$
is the density bundle associated to
$\sct \Xb$, so locally near $\bXb$ it is spanned by $x^{-n-1}\,dx\,dy$ over
$\Cinf(\Xb)$.
Finally, $\Diffsc(\Xb)$ is the algebra of differential operators generated
by the vector fields in $\Vsc(\Xb)$; $\Diffsc^m(\Xb)$ stands for scattering
differential operators of order (at most) $m$.

An example is provided by the radial compactification of
Euclidean space, $\Xb=\Snp$. 
We can use `inverse' polar coordinates on $\Rn$ to induce local
coordinates on $\Snp$ near $\partial\Snp$ as above. Thus, we let
$x=r^{-1}=|w|^{-1}$ (for $r\geq 1$, say, smoothed out near
the origin), as in the introduction,
write
$w=x^{-1}\omega$, $\omega\in\Sn$, $|w|>1$, and let $y_j$,
$j=1,\ldots,n-1$, be local coordinates on $\Sn$. For
example, one can take the $y_j$ to be $n-1$ of the $w_k/|w|$, $k=1,\ldots,
n$. Then $x\in\Cinf(\Snp)$ is a boundary defining function of
$\Snp$, and $x$ and the $y_j$ give local coordinates near $\partial\Snp=\Sn$.

To establish the
relationship between the scattering structure and
the Euclidean scattering theory, we
introduce local coordinates on $\Xb$ near $p\in\bXb$ as above, and
use these to identify the coordinate neighborhood $U$ of $p$ with
a coordinate patch $U'$ on the closed upper hemisphere $\Snp$ (which is
just a closed ball) near its boundary.
Such an identification preserves the scattering structure since this
structure is
completely natural. We
further identify $\Snp$ with $\Rn$ via the
radial compactification $\SP$ as in \eqref{eq:SP-def}.
The constant coefficent vector fields $\partial_{w_j}$ on $\Rn$ lift
under $\SP$
to give a basis of $\Tsc\Snp$. Thus, $V\in\Vsc(\Snp)$ can
be expressed as (ignoring the lifting in the notation)
\begin{equation}\label{eq:sc-vfields}
V=\sum_{j=1}^n a_j\partial_{w_j},\quad a_j\in\Cinf(\Snp).
\end{equation}
As mentioned above, $a_j\in\Cinf(\Snp)$ is equivalent to
requiring that $\SP^* a_j$ is a classical (i.e.\ one-step polyhomogeneous)
symbol of order $0$ on $\Rn$. This description also shows that
the positive Euclidean Laplacian, $\Delta$, is an element of
$\Diffsc^2(\Snp)$, and that $\Osc\Snp$ is spanned by the pull-back of
the standard Euclidean density $|dw|$.
In terms of the `inverted' polar coordinates on $\Rn$,
covectors $\xi\cdot dw=\sum_j\xi_j\,dw_j$
take the form \eqref{eq:sct-coords} with
\begin{equation}
\tau=-\frac{w\cdot\xi}{|w|}=-y\cdot \xi,\ \tau^2+|\mu|^2=|\xi|^2.
\end{equation}
Here $\mu$ is the orthogonal
projection of $\xi$ to the tangent space of the unit sphere
$\Sn$ at $y\in\Sn$, and
$|\mu|$ denotes the length of a covector on $\Sn$ with respect to the
standard metric $h$ on the unit sphere.

We next show that polyhomogeneous symbols on $X^a$, pulled back to $\Rn$
by $\pi^a$, are smooth on the blown-up space $[\Snp;C_a]$. Recall that
the blow-up process is simply an invariant way of introducing
polar coordinates about a submanifold. A full description appears
in \cite{RBMDiff} and a more concise one in \cite[Appendix~A]{RBMSpec},
but we give a brief summary here. Near $C_a$, $\Snp$ can be identified
with the inward-pointing normal bundle $N^+ C_a$ of $C_a$. Here
$N^+ C_a$ is the image of the inward-pointing tangent bundle,
$T^+\Snp$, in the quotient bundle $N C_a=T_{C_a}\Snp/T C_a$; $W\in T_q\Snp$
is inward pointing if $Wx(q)\geq 0$. The blow
up $[\Snp;C_a]$ of $C_a$ in $\Snp$ is then locally given by the blow up
of $N^+ C_a$ at the zero section, which amounts to introducing
polar coordinates in its fibers. The front face of the blow up (the lift
of $C_a$) is then identified with the inward pointing sphere bundle,
$S^+ NC_a$, which is the quotient of $N^+C_a\setminus 0$ by the natural
$\Real^+$ action in its fibers ($0$ denotes the zero section).
We denote the blow-down map (which is $\Cinf$)
by $\beta[\Snp;C_a]:[\Snp;C_a]\to \Snp$.
Now $S^+ NC_a$ is a hemisphere bundle
over $C_a$, which can be identified with the
radial compactification of the normal bundle of $C_a$ in $\Sn$
whose fibers can in turn be identified with  $X^a$.

To see this
in more concrete terms,
we proceed by finding local coordinates on $[\Xb;C_a]$ explicitly.
Near $C_a$ in $\Snp$ we have $|w_a|>c|w^a|$ for some $c>0$.
Hence, near any point $p\in C_a$ one of the coordinate functions $(w_a)_j$
which we take to be $(w_a)_m$ for the sake of definiteness,
satisfies $|(w_a)_m|>c'|(w_a)_j|$,
$|(w_a)_m|>c'|w^a|$ for some $c'>0$. Taking into account the
coordinate form of $\SP$ we see that
\begin{equation}\label{eq:Snp-C_a-coords}
x=|w_a|^{-1},\ y_j=\frac{(w_a)_j}{|w_a|}\ (j=1,\ldots,m-1),
\ z_j=\frac{(w^a)_j}{|w_a|}\ (j=1,\ldots,n-m)
\end{equation}
give coordinates on $\Snp$ near $p$ where we think of the
$y_j$ as coordinates on the unit sphere $C_a$.
In these coordinates $C_a$ is
defined by $x=0$, $z=0$. Sometimes we write these coordinates as
$x_a$, $y_a$, $z_a$ to emphasize the collision plane. {\em Note that
here the notation regarding $y$ and $z$ is reversed as compared to
\cite{Vasy:Propagation-Many}.}
We call these the standard coordinates around $C_a$; we will always use these
in this paper.
Correspondingly, we have coordinates
\begin{equation}\label{eq:ff-coords}
x,\ y_j\ (j=1,\ldots,m-1),\ Z_j=z_j/x\ (j=1,\ldots,n-m),
\end{equation}
i.e.
\begin{equation}\label{eq:ff-coords-2}
x=|w_a|^{-1},\ y_j=\frac{(w_a)_j}{|w_a|}\ (j=1,\ldots,m-1),
\ Z_j=(w^a)_j\ (j=1,\ldots,n-m)
\end{equation}
near the interior of the front face $\ff$
of the blow-up $[\Xb;C_a]$, i.e.\ near the
interior of $\ff=\beta[\Xb;C_a]^*C_a$.
Similarly, one can easily write down local coordinates near
the corner $\ff\cap\beta[\Xb;C_a]^*\bXb$, see \cite[Section~2]
{Vasy:Propagation-Many}. As a result of such calculations, we deduce
the following lemma and its corollary.

\begin{lemma}\label{lemma:blowup-1}(\cite[Lemma~2.5]{Vasy:Propagation-Many})
Suppose that $\Xb=\Snp$ and let $\beta=\beta[\Xb;C_a]$ be the
blow-down map. Then the pull-back $\beta^*(\SP^{-1})^*\pi^a$ of
$\pi^a:\Rn\to X^a$ extends to a $\Cinf$ map, which we also
denote by $\pi^a$,
\begin{equation}\label{eq:geom-11}
\pi^a:[\Xb;C_a]\to \Xb^a.
\end{equation}
Moreover, if $x^a$ is a boundary defining function on $\Xb^a$
(e.g.\ $x^a=|w^a|^{-1}$ for $|w^a|>1$), then $\rho_{\bXb}=
(\pi^a)^* x^a$ is a defining function for
the lift of $\partial \Xb$ to $[\Xb;C_a]$, i.e.\ for $\beta^*\partial \Xb$.
\end{lemma}

\begin{cor}(\cite[Corollary~2.6]{Vasy:Propagation-Many})
Suppose that $\Xb=\Snp$, $f\in S^r_\phg(X^a)$. Then
\begin{equation}
(\pi^a)^*f\in\rho_{\bXb}^{-r}\Cinf([\Xb;C_a]).
\end{equation}
Here, following the previous lemma, we regard $\pi^a$ as the map
in \eqref{eq:geom-11},
and $\rho_{\bXb}$ is
the defining function of $\beta[\Xb;C_a]^*\bXb$, i.e.\ of the lift
of $\bXb$, and the
subscript $\phg$ refers to classical (one-step polyhomogeneous) symbols.
\end{cor}

This corollary shows that for a Euclidean many-body Hamiltonian, $H=\Delta+
\sum_a V_a$, $V_a$ becomes a smooth function on the compact resolved
space $[\Snp;C_a]$. Thus, to understand $H$, we need to blow up
{\em all} the $C_a$.
The iterative construction was carried out in detail in
\cite[Section~2]{Vasy:Propagation-Many}; we refer to the discussion given
there for details.
However, we remind the reader that the $C_a$ are blown up in the order
of inclusion ({\em opposite} to the usual order on the clusters). That is, one
starts with the blow-up of 2-clusters (which are disjoint);
3-clusters become disjoint upon this blow-up. One proceeds to blow-up
the 3-clusters; 4-clusters become disjoint now. One proceeds this way,
finally blowing up the $N-1$-clusters. (The blow-up of the $N$-cluster is
a diffeomorphism, hence can be neglected.) We thus obtain a manifold
with corners which is denoted by $[\Snp;\calC]$, and the blow-down map
(which is $\Cinf$) is written as
\begin{equation}
\betaSc=\beta[\Snp;\calC]:[\Snp;\calC]\to\Snp.
\end{equation}
The algebra of many-body differential operators is then defined as
\begin{equation}
\DiffSc(\Snp,\calC)=\Cinf([\Snp;\calC])\otimes_{\Cinf(\Snp)}
\Diffsc(\Snp).
\end{equation}
That is, similarly to \eqref{eq:sc-vfields},
$P\in\DiffSc^m(\Snp,\calC)$ means that
\begin{equation}
P=\sum_{|\alpha|\leq m} a_\alpha D^\alpha_w,\quad
a_\alpha\in\Cinf([\Snp;\calC])
\end{equation}
where we again ignored the pull-back by $\SP$ in the notation.

The coordinates on $\sct \Xb$ induced by the standard coordinates around $C_a$
are
\begin{equation}
\tau_a\frac{dx_a}{x_a^2}+\mu_a\cdot\frac{dy_a}{x_a}+\nu_a\cdot\frac{dz_a}{x_a};
\end{equation}
we included the subscript $a$ to emphasize
the element $C=C_a$ of $\calC$ around which the local
coordinates are centered. {\em Again, the roles of $\mu$ and $\nu$ have
been switched from \cite{Vasy:Propagation-Many}.}
In terms of the splitting $\Rn=X_0=X_a\oplus X^a$,
and the corresponding splitting $\xi=(\xi_a,\xi^a)$ of the dual
coordinates,
this gives
\begin{equation}
\xi^a=\nu_a,\ \text{and, at}\ C_a,\ \tau_a=\tau,\ \mu_a\ \text{the
orthogonal projection of}\ \xi_a\ \text{to}\ T_y C_a.
\end{equation}

One of the main differences between $\DiffSc(\Snp,\calC)$ and
$\Diffsc(\Snp)$ is
that the former is not commutative to `top weight'. That is, while
for $P\in\Diffsc^m(\Snp)$, $Q\in\Diffsc^{m'}(\Snp)$, we have $[P,Q]\in x\Diffsc
^{m+m'-1}(\Snp)$, this is replaced by $[P,Q]\in\rho_{C_0}
\DiffSc^{m+m'-1}(\Snp,\calC)$
for $P\in\DiffSc^m(\Snp,\calC)$, $Q\in\DiffSc^{m'}(\Snp,\calC)$
with $\rho_{C_0}$ a defining function for the lift of
$C_0$. Thus, a vanishing factor
(such as $x$ above) is only present
at the lift of the free face $C_0$, i.e.\ there is no gain
of a weight factor at the front faces $\ff$. From an algebraic point of
view, it is this non-commutativity that causes the breaks in the propagation
of singularities by necessitating the use of the compressed bundle
for wave front sets; one can only microlocalize in the commutative
(tangential) variables, i.e.\ at $C_a$ (where $z_a=0$) only in $x_a$, $y_a$
(i.e.\ $w_a$), and $\tau_a$, $\mu_a$ (i.e.\ $\xi_a$).

For a closed embedded submanifold $C$ of $\bXb$,
the relative scattering tangent bundle $\Tsc(C;\Xb)$ of $C$ in $\Xb$ is
the subbundle
of $\Tsc_C \Xb$ consisting of
$v\in\Tsc_p \Xb$, $p\in C$, for which there exists
\begin{equation}
V\in\Vsc(\Xb;C)\subset\Vsc(\Xb)
\end{equation}
with $V_p=v$ (equality understood in $\Tsc_p \Xb$). Here
\begin{equation}
\Vsc(\Xb;C)=x\Vb(\Xb;C)=x\{V\in\Vb(\Xb):\ V\ \text{is tangent to}\ C\}
\end{equation}
and tangency is defined using the (non-injective) inclusion map
$\Tb X\to T\Xb$.

Thus, in local coordinates $(x,y,z)$ near $p\in C$ such that
$C$ is defined by $x=0$, $z=0$, $\Tsc(C;\Xb)$ is spanned by
$x^2\partial_x$ and $x\partial_{ y_j}$, $j=1,\ldots, m-1$ where
$n-m$ is the codimension on $C$ in $\bXb$.
Thus, in our
Euclidean setting, $\Xb=\Snp$, $C=C_a=\partial\Xb_a$, $g$ the Euclidean
metric, $\Tsc(C;\Xb)$ is naturally isomorphic to $\Tsc_C \Xb_a$, i.e.\ it
should be regarded as the bundle of scattering
tangent vectors of the collision
plane at infinity, spanned by $\partial_{(w_a)_j}$, $j=1,\ldots,m$, $m
=\dim X_a$.

For $C=C_a\in\calC$,
the metric $g$ defines the orthocomplement $(\Tsc(C;\Xb))^\perp$
of $\Tsc (C;\Xb)$ in $\Tsc_C \Xb$.
The subbundle of $\sct_C \Xb$
consisting of covectors that
annihilate $(\Tsc(C;\Xb))^\perp$
is denoted by $\sct(C;\Xb)$; it is called
the relative scattering cotangent bundle of
$C$ in $\Xb$.
In our
Euclidean setting, $\sct(C_a;\Xb)$ is naturally isomorphic to $\sct_{C_a} \Xb_a$
and is spanned by $d(w_a)_j$, $j=1,\ldots,m$; so we simply write
\begin{equation}
\sct(C_a;\Xb)=\sct_{C_a} \Xb_a.
\end{equation}

In our standard local coordinates $(x,y,z)$ near $p\in C$, $C=C_a$,
$C$ defined by $x=0$, $z=0$, and $x\partial_{z_j}$ give
an orthonormal basis of $(\Tsc (C;\Xb))^\perp$. Note that a basis
of $\Tsc (C;\Xb)$ is given by $x^2\partial_x$ and $x\partial_{y_j}$,
while a basis of $\sct (C;\Xb)$ is given by $x^{-2}\,dx$, $x^{-1}\,dy_j$.
A covector in $\sct \Xb$ can be written in these local coordinates as
\begin{equation}
\tau_a\frac{dx_a}{x_a^2}+\mu_a\cdot\frac{dy_a}{x_a}+\nu_a\cdot\frac{dz_a}{x_a}.
\end{equation}
Thus, local coordinates on $\sct_{\bXb} \Xb$ are given by $(y,z,\tau,\mu,\nu)$,
while on $\sct (C_a;\Xb)$ by $(y,\tau,\nu)=(y_a,\tau_a,\mu_a)$.

Now if $C=C_a$, $C_b\in\calC$ with $C_a\subset C_b$, we can further adjust
our coordinates so that $C_b$ is defined by $x=0$, $z''=0$, for some
splitting $z=(z',z'')$. In fact, we can use the decomposition
\begin{equation}\label{eq:X_a-X_b}
\Rn=X_a\oplus X^a_b\oplus X^{ab},\quad X^a_b=X^a\cap X_b,
\ X^{ab}=X^a\cap X^b;
\end{equation}
the sum is of course still orthogonal with respect to $g$.
Then the splitting $z=(z',z'')$ corresponds to the splitting $X^a=X^a_b\oplus
X^{ab}$ in terms of the coordinates \eqref{eq:Snp-C_a-coords}.
Following \eqref{eq:X_a-X_b}, we write these coordinates as
$z_a=(z_{ab},z_a^b)$

With the corresponding splitting of the dual
variable, $\nu=(\nu',\nu'')=(\nu_{ab},\nu_a^b)$,
we obtain a well-defined projection
\begin{equation}
\pi_{ba}:\sct_{C_a} (C_b;\Xb)\to\sct (C_a;\Xb),
\end{equation}
\begin{equation}
\pi_{ba}(0,y_a,\tau_a,\mu_a,\nu_{ab})=(y_a,\tau_a,\mu_a).
\end{equation}
In our Euclidean setting this is just the obvious projection
\begin{equation}
\pi_{ba}:\sct_{\partial \Xb_a}\Xb_b\to\sct_{\partial\Xb_a} \Xb_a
\end{equation}
under the inclusion $\Xb_a\subset \Xb_b$.
We write $\pi$ for the collection of these maps.

\section{Indicial operators and the characteristic set}\label{sec:indicial}
We start by recalling from \cite{Vasy:Propagation-Many}
the definition of the many-body pseudo-differential
calculus via the quantization of symbols.
Thus, the non-polyhomogeneous space
$\PsiScc^{m,l}(\Snp,\calC)$ is the following.
Suppose that, with the notation of \cite{RBMCalcCon},
\begin{equation}\label{eq:PsiScc-5}
a\in\bcon^{-m,l}([\Snp;\calC]\times\Snp).
\end{equation}
That is,
identifying $\interior(\Snp)$ and $\interior([\Snp;\calC])$ with $\Rn$
as usual (via $\SP^{-1}$), suppose that
$a\in\Cinf(\Rn_w\times\Rn_\xi)$ has the following property.
For every $P\in\Diffb^k(\Snp)$, acting on the second factor of
$\Snp$ (i.e.\ in the $\xi$ variable), and
$Q\in\Diffb^{k'}([\Snp;\calC])$, acting
on the first factor of $\Snp$ (i.e.\ in the $w$ variable),
$k,k'\in\Nat$,
\begin{equation}\label{eq:symb-est-3}
PQa\in\rho_{\infty}^{-m}\rho_{\partial}^l L^{\infty}(\Snp\times\Snp)
\end{equation}
where $\rho_\infty$ and $\rho_\partial$ are defining functions of the first
and second factors of $\Snp$ respectively. Let $A=q_L(a)$ denote the
left quantization of $a$:
\begin{equation}\label{eq:q_L-def}
Au(w)=(2\pi)^{-n}\int e^{i(w-w')\cdot\xi}a(w,\xi)u(w')\,dw'\,d\xi,
\end{equation}
understood as an oscillatory integral. Then $A\in\PsiScc^{m,l}(\Snp,\calC)$.
We could have equally well used other (right, Weyl, etc.) quantizations as
well.

In addition, the polyhomogeneous class $\PsiSc^{m,l}(\Snp,\calC)$
is given by the quantization of symbols
\begin{equation}\label{eq:cl-symb-def}
a\in\rho_{\infty}^{-m}\rho_{\partial}^l\Cinf([\Snp;\calC]\times\Snp).
\end{equation}
Since
differential operators $\sum a_\alpha(w) D^\alpha$ are just the
left quantization of the symbols $a(w,\xi)=\sum a_\alpha(w)\xi^\alpha$, it
follows immediately that
\begin{equation}
\DiffSc^m(\Snp,\calC)\subset\PsiSc^m(\Snp,\calC).
\end{equation}

It was shown in \cite{Vasy:Propagation-Many} that $\PsiScc(\Snp,\calC)$
and $\Psisc(\Snp,\calC)$ are $*$-algebras with respect to composition
and taking adjoints. Moreover, $\PsiScc^{m,l}(\Snp,\calC)$ is
bounded from $\Hsc^{r,s}(\Snp)$ to $\Hsc^{r-m,s+l}(\Snp)$; this
follows from the inclusion $\PsiScc^{m,0}(\Snp,\calC)\subset
\Psi_\infty^m(\Rn)$ where the latter is H\"ormander's `uniform' calculus
in $\Rn$, see \cite[Section~18.1]{Hor}.

In the standard pseudo-differential calculus on compacts manifolds without
boundary Fredholm properties are captured by the principal symbol. In our
setting,
H\"ormander's principal symbol map on $\Psi_\infty^m(\Rn)$ restricts to
a principal symbol map
\begin{equation}\label{eq:calc-23}
\prs_{\Scl,m}:\PsiSc^{m,0}(\Snp,\calC)\to\Shom^m(\ScT[\Snp;\calC]),
\end{equation}
$\Shom(\ScT[\Snp;\calC])$
denoting the space of smooth symbols which are homogeneous
of degree $m$ on $\ScT[\Snp;\calC]$, the pull-back of $\sct\Snp$ by the
blow-down map. However, the invertibility of the principal symbol
is not sufficient for an operator to be Fredholm, since the behavior
near infinity must be taken into account as well. The additional
piece of information that is necessary is the invertibility of the
indicial operators whose definition we briefly recall below.
First, however, we remind the reader what the
structure of the lift of $C_a$ to $[\Xb;\calC]$ is.

For $C_a\in\calC$ let
\begin{equation}
\calC_a=\{C_b\in\calC:\ C_b\subsetneq C_a\},
\end{equation}
\begin{equation}
\calC^a=\{C_b\in\calC:\ C_a\subsetneq C_b\}.
\end{equation}
We carry out the blow-up $[\Snp;\calC]$ by first blowing up $\calC_a$. Since
all elements of $\calC_a$ are p-submanifolds (i.e.\ product submanifolds)
of $C_a$, the lift
$\beta[\Snp;\calC_a]^*C_a$ of $C_a$ to $[\Snp;\calC_a]$ is naturally
diffeomorphic to
\begin{equation}
\Ct_a=[C_a;\calC_a].
\end{equation}
Thus, over $C'_a$, the regular part of $C_a$, $\Ct_a$
can be identified with $C_a$.
The front face of the new blow-up, i.e.\ of the blow up of
$\beta[\Snp;\calC_a]^*C_a$ in $[\Snp;\calC_a]$ is thus a hemisphere
(i.e.\ ball) bundle
over $\Ct_a$, namely $S^+N\Ct_a$. We write the bundle projection, which
is just the restriction of the new blow-down map to the front face,
$S^+N\Ct_a$ as
\begin{equation}
\rho_a:S^+N\Ct_a\to\Ct_a.
\end{equation}
In our Euclidean
setting, these fibers can be naturally identified with $\Xb^a$ via
the projection $\pi^a$ (extended as in Lemma~\ref{lemma:blowup-1}).

Every
remaining blow up in $[\Snp;\calC]$ concerns submanifolds that are either
disjoint from this new front face or are the lift of elements of $\calC^a$.
The former do not affect the structure near the new front face,
$S^+N\Ct_a=\beta[\Snp;\calC_a;C_a]^*C_a$, while
the latter, which are given by the lifts of elements of $\calC^a$,
correspond to blow ups that can be performed in the fibers of $S^+N\Ct_a$.
Note that
the lift of $C_b\in\calC^a$, meets the new front face only at
its boundary since all $C_b$ are subsets of $\Sn=\partial\Snp$.
In particular, the lift $\betaSc^*C_a$
of $C_a$ to $[\Snp;\calC]$ fibers over $\Ct_a$ and the
fibers are diffeomorphic to a hemisphere (i.e.\ ball) with certain
boundary submanifolds blown up. More specifically, the intersection
of $\beta[\Snp;\calC_a;C_a]^*C_b$, $C_b\in\calC^a$, with the
front face $S^+N\Ct_a$ is the image of $T\beta[\Snp;\calC_a]^*C_b$
under the
quotients defining the spherical normal bundle;
$\betaSc^*C_a$ is obtained by blowing these up
in $S^+N\Ct_a$. Hence, the fiber of $\betaSc^*C_a$ over $p\in\Ct_a$
is given by $[S^+N_q C_a;T_q \calC^a]$ where $q=\beta[\Snp;\calC_a](p)\in
C_a$.
In particular, in our Euclidean
setting, the fibers of $\betaSc^*C_a$ over $\Ct_a$
can be naturally identified with
$[\Xb^a;\calC^a]$ via $\pi^a$. 

We now define $\sct (\Ct_a;\Snp)$ denote the pull-back of
$\sct(C_a;\Snp)$ by the blow-down
map $\beta[C_a;\calC_a]$:
\begin{equation}
\sct (\Ct_a;\Snp)=\beta[C_a;\calC_a]^*\sct(C_a;\Snp).
\end{equation}
Thus, $\sct (\Ct_a;\Snp)$ is, via the metric $g$,
naturally diffeomorphic to $\Ct_a\times X_a$, where the metric is used
to identify $X_a$ with its dual.
If $C_a\subset C_b$ then $\pi_{ba}$ lifts to a map
\begin{equation}
\pit_{ba}:\sct_{\beta[C_b;\calC_b]^*C_a} (\Ct_b;\Xb)\to\sct(\Ct_a;\Xb).
\end{equation}

For each operator $A\in\PsiSc^{m,l}(\Xb,\calC)$,
the $C_a$-indicial operator of $A$,
denoted by $\Ah_{a,l}$, will be a collection of operators, one for each $\zeta
\in\sct_p(\Ct_a;\Xb)$, acting on functions on the fiber $\rho_a^{-1}(p)$
of $\rho_a$. Recall that the interior of
these fibers can be naturally identified with
$X^a$. So suppose that $u\in\dCinf(\rho_a^{-1}(p))$, i.e.\ $u$ is
a Schwartz function on the the fiber above $p$ which is a compactification
of $X^a$; we
need to define $\Ah_a(\zeta)u$. For this purpose
choose $\ft\in\Cinf(\Snp;\Real)$ such that
$d(\ft/x)$, evaluated at $\beta[C_a;\calC_a](p)$,
is equal to $\zeta$. Then let
$\At=e^{-i\ft/x}x^{-l}Ae^{i\ft/x}\in\PsiSc^{m,0}(\Xb,\calC)$,
and choose $u'\in\Cinf([\Xb;\calC])$ such that
$u'|_{\rho_a^{-1}(p)}=u$. Then
\begin{equation}
\Ah_{a,l}(\zeta)u=(\At u')|_{\rho_a^{-1}(p)},
\end{equation}
which is independent of all the choices we made.

In our Euclidean setting we can simply take $f$ to be the pull-back
of a function in $\Cinf(\Xb_a)$ (which we write, abusing the notation,
as $f\in\Cinf(\Xb_a)$). Thus, at least in the regular part $C'_a$ of
$\Ct_a$, the $a$-indicial operator of $A$ arises
by thinking of $A$ as a scattering pseudo-differential operator
on $\Xb_a$ with values in operators on $X^a$, and finding its scattering
symbol at $C'_a\subset C_a=\partial\Xb_a$. In fact, $\Ah_{a,l}(\xi)$ is
a smooth function of $\xi\in\sct(\Ct_a;\Xb)$, so the above description also
determines $\Ah_{a,l}$ at $\partial\Ct_a$.
We often simplify (and thereby abuse) the
notation and drop the index $l$, i.e.\ we write $\Ah_a=\Ah_{a,l}$, when
the value of $l$ is understood.

In the case of Euclidean many-body scattering,
$\Ah_{a,l}$ is
a function on $\beta_a^*\sct_{C_a}\Xb_a$
with values in operators on $\Sch(X^a)$; here
\begin{equation}
\beta_a=\beta[C_a;\calC_a]:\Ct_a=[C_a;\calC_a]\to C_a
\end{equation}
is the blow-down map. Thus, $\beta_a$ is simply the
restriction of $\beta[\Xb_a;\calC_a]$ to the lift $\Ct_a
=\beta[\Xb_a;\calC_a]^* C_a$. Then
\begin{equation}\label{eq:ind-35}
\Ah_{a,l}\in
\Cinf(\beta_a^*\sct_{\partial\Xb_a}\Xb_a,\PsiSc^{m,0}(\Xb^a,\calC^a)).
\end{equation}
We write the space on the right as
\begin{equation}
\Ah_{a,l}\in\PsiScra^{r,0}(\beta_a^*\sct(\Ct_a;\Xb),\calCt_a).
\end{equation}
An explicit formula for the indicial operator of $A$, given by
\eqref{eq:q_L-def}, at $\Ct_a$ is
\begin{equation}\label{eq:Psop-46}
\Ah_a(y_a,\xi_a)u(Z_a)=(2\pi)^{-\dim X^a}
\int e^{i(Z_a-w^a)\cdot\xi^a}a(0,y_a,Z_a,\xi_a,\xi^a) u(w^a)\,d\xi^a\,dw^a.
\end{equation}
Note that the oscillatory testing definition that we adopted automatically
shows that the indicial operator map is multiplicative.

The principal symbol and the indicial operators together
indeed describe compactness properties of
ps.d.o's. In particular, if $A\in\PsiSc^{m,l}(\Snp,\calC)$,
$\prs_{\Scl,m}(A)$ never vanishes, and for  every $a$, $\Ah_{a,l}$ is
invertible in $\PsiScra^{r,0}(\beta_a^*\sct(\Ct_a;\Xb),\calCt_a)$, then
there is a two-sided parametrix $P$ for $A$ such that $PA-\Id,
AP-\Id\in\PsiSc^{-\infty,\infty}(\Snp,\calC)$.
For $A\in\PsiSc^{m,0}(\Snp,\calC)$ self-adjoint, $m>0$, such that
$\prs_{\Scl,m}(A)$
never vanishes (but no assumptions made on the indicial operators of $A$),
$\lambda\in\Cx\setminus\Real$,
we conclude that $(A-\lambda)^{-1}\in\PsiSc^{-m,0}(\Snp,\calC)$. The
invertibility of the indicial operators of $A-\lambda$
is automatic since $\widehat{(A-\lambda)}_{a,0}=\Ah_{a,0}-\lambda$, and
then an iterative argument, starting at $a=0$, using the self-adjointness
of $A$, shows the required invertibility. In addition, an argument using
almost analytic extensions and the Cauchy representation also shows
that certain  functions of self-adjoint ps.d.o's are ps.d.o's themselves.
In particular, for
$A\in\PsiSc^{m,0}(\Snp,\calC)$ self-adjoint, $m>0$, such that
$\prs_{\Scl,m}(A)$ never vanishes,
$\phi\in\Cinf_c(\Real)$, we conclude
that $\phi(A)\in\PsiSc^{-\infty,0}(\Snp,\calC)$
(see \cite[Propositions~4.7-4.8]{Vasy:Propagation-Many}, and
see \cite{Derezinski-Gerard:Scattering} and \cite{Hassell-Vasy:Symbolic}
for a general discussion of functions of ps.d.o.'s).

Now it is particularly easy
to identify the indicial operators: for each $\xi=(y,\tau,\mu)
\in\sct_p(\Ct_a;\Xb)$,
they are operators on $X^a$, identified as the fiber of the front face
$\rho_a^{-1}(p)$ over $p$.
A simple calculation, see \cite[Sections~4 and 11]{Vasy:Propagation-2}
for more details, shows that
the indicial operators of $H$ are given by
\begin{equation}\label{eq:ham-4}
\Hh_{a,0}(\xi)=\Hh_{a,0}((p,0))+
|\xi|^2,\ \xi=(y,\tau,\mu)\in\sct(\Ct_a;\Xb),
\end{equation}
\begin{equation}\label{eq:ham-5}
\Hh_{a,0}(p,0)=\Delta_Z+V(p,Z)=H^a
\end{equation}
where $Z$ are Euclidean
coordinates on the interior of $\rho_a^{-1}(p)$, i.e. on
$X^a$, and $\Delta_Z$ is the
Euclidean Laplacian.
Note that under these natural identifications (that $\rho_a^{-1}(p)$ is
regarded as $X^a$ for all $p$) the $a$-indicial operators of $H$ only depend
on $|\xi|^2$, i.e.\ on the metric function of $g_a$.

Equations \eqref{eq:ham-4}-\eqref{eq:ham-5} show
that $\Hh_{a,0}(p,0)$ is uniformly bounded
below, so for any $\psi\in\Cinf_c(\Real)$ the set
\begin{equation}
\cup_a\cl(\{\xi\in\sct(\Ct_a;\Xb):\ \psi(\Hh_{a}(\xi))\neq 0\})
\end{equation}
is compact.

Note that $\Hh_{a,0}(\xi)-\lambda$ is invertible in $\PsiSc(\Xb^a,\calC^a)$
if and only if $\lambda-|\xi_a|^2$ is not in the spectrum of $H^a$.
Indeed, by the HVZ theorem, this
means that $\lambda-|\xi_a|^2$ is not an $L^2$ eigenvalue of
$H^a$, and it is smaller than all of its thresholds, $\lambda-|\xi_a|^2<\inf
\Lambda_a$,
i.e.\ smaller
than the eigenvalues of subsystems of $H^a$. By an iterative argument
corresponding to $C_a\subsetneq C_b$,
$\lambda-|\xi_a|^2<\inf\Lambda_a$
guarantees the invertibility of the indicial operators of
$\Hh_{a,0}(\xi)-\lambda$, so $\Hh_{a,0}(\xi)-\lambda$ is Fredholm, and
hence invertible if and only if it has no $L^2$ eigenfunctions. Thus, we define
the characteristic variety $\dot\Sigma(\lambda)\subset\scdt \Xb$
of $H-\lambda$ to be the union image of these
sets under $\pi$:
\begin{equation}
\dot\Sigma(\lambda)=\cup_a \dot\Sigma_a(\lambda)\subset\scdt \Xb,
\end{equation}
\begin{equation}\begin{split}
\dot\Sigma_a(\lambda)=&\{\xi_a\in\sct_{C'_a}
(C_a;\Xb):\ \lambda-|\xi_a|^2\geq
\sigma\ \text{for some}\ \sigma\in\pspec{H^b},\ C_a\subsetneq C_b\}\\
&\cup\{\xi_a\in\sct_{C'_a}
(C_a;\Xb):\ \lambda-|\xi_a|^2\in\pspec{H^a}\}.
\end{split}\end{equation}

If there are no bound states in the subsystems, a compactness argument
as in \cite{Vasy:Propagation-Many}
allows one to show that the free Laplacian governs the propagation of
singularities except that breaks occur at the collision planes $C_a$ where
the usual law of reflection is satisfied. This indicates that one should
think of the characteristic variety $\dot\Sigma(\lambda)$ as the union
of projection of the characteristic sets corresponding to the bound states
of the various subsystems. Thus, let
\begin{equation}
\Sigma_b(\lambda)=\{\xi_b\in\sct_{C_b}\Xb_b:
\ \lambda-|\xi_b|^2\in\pspec{H^b}\}.
\end{equation}
For $a$ with $C_a\subset C_b$,
the projection of $\Sigma_b(\lambda)$ under $\pi_{ba}$ is
\begin{equation}\begin{split}
\pi_{ba}(\Sigma_b(\lambda))&=\{\xi_a\in\sct_{C_a}\Xb_a:
\ \exists \xi_b\in\sct_{C_b}\Xb_b\Mwith \pi_{ba}(\xi_b)=\xi_a\Mand
\lambda-|\xi_b|^2\in\pspec{H^b}\}\\
&=\{\xi_a\in\sct_{C_a}\Xb_a:\ \lambda-|\xi_a|^2\geq
\sigma\ \text{for some}\ \sigma\in\pspec{H^b}\}.
\end{split}\end{equation}
Thus,
\begin{equation}\label{eq:dSa-as-proj-Sb}
\dot\Sigma_a(\lambda)
=\cup_{C_b\supset C_a}
\pi_{ba}(\Sigma_b(\lambda))\cap\sct_{C'_a}(C_a;\Xb).
\end{equation}
In fact, letting $\pi_b:\sct_{C_b}\Xb_b\to\scdt \Xb$ be the projection, we
thus conclude that
\begin{equation}\label{eq:dot-Sigma-Sigma_b}
\dot\Sigma(\lambda)=\cup_b \pi_b(\Sigma_b(\lambda)).
\end{equation}
We write $\pih_b$ for the restriction of $\pi_b$ to $\Sigma_b(\lambda)$.

We proceed to recall the basic topological properties of $\scdt\Xb$ from
\cite[Section~5]{Vasy:Propagation-Many}.
We put the topology induced by $\pi:\sct_{\bXb}\Xb\to\scdt \Xb$
on $\scdt \Xb$.
Thus, $C\subset\scdt\Xb$ is closed if and only if
$\pi^{-1}(C)$ is closed, so if $f$ is continuous on $\sct_{\bXb}\Xb$
and $\pi$-invariant, then $f_\pi$ is continuous on $\scdt\Xb$.
In particular, there are always continuous functions separating points
on $\scdt \Xb$:
if $p(\xib)\neq p(\xib')$ (here $p:\scdt\Xb\to\bXb$ stands
for projection to the base), one can use the pull-back of an appropriate
function on $\bXb$, and if $p(\xib)=p(\xib')\in C'_a$, then $\xib,\xib'\in
\sct_{C'_a}\Xb_a$ are of the form $(\yb_a,\xib_a)$ and $(\yb'_a,\xib'_a)$,
$\xib_a\neq\xib'_a$,
and the function $\xi\mapsto\xi_a$ is well-defined and $\pi$-invariant
on a neighborhood
of $p(\xib)$, so multiplying it by the pull-back of a cutoff on $\bXb$
gives a globally well-defined separating continuous function. Thus,
$\scdt\Xb$ is Hausdorff.

By the continuity of $\pi$, if $K\subset\sct_{\bXb}\Xb$ is compact, the same
holds for $\pi(K)$. But by \eqref{eq:dot-Sigma-Sigma_b},
$\dot\Sigma(\lambda)=\pi(\cup_b\imath_b(\overline{\Sigma_b(\lambda)}))$
(see \eqref{eq:imath_b-def}),
and $\imath_b(\overline{\Sigma_b(\lambda)})$ is compact for each $b$, so
$\dot\Sigma(\lambda)$ is compact, hence closed in $\scdt\Xb$.
In particular $C\subset\dot\Sigma(\lambda)$
is closed in $\dot\Sigma(\lambda)$ if and only if it is closed
in $\scdt \Xb$, i.e.\ if and only if $\pi^{-1}(C)$ is closed in $\sct_{\bXb}
\Xb$. Also,
if $C\subset\dot\Sigma(\lambda)$ is closed, then $\pih_0^{-1}(C)$ is compact.

Note that for all $C\subset\dot\Sigma(\lambda)$ and for
any $R>0$ with $R>\lambda-\inf\Lambda_1$,
\begin{equation}\label{eq:C-cl-comp}
\pi(\pi^{-1}(C))=\pi(\pi^{-1}(C)\cap\{\xi\in\sct_{\bXb}\Xb:\ |\xi|^2\leq R\}).
\end{equation}
Indeed, if $\xi\in\pi^{-1}(C)$, $\pi(\xi)\in\sct_{C'_a}\Xb_a$, then
by $C\subset\dot\Sigma(\lambda)$, $|\imath_a(\pi(\xi))|^2<R$,
$\pi(\imath_a(\pi(\xi)))=\pi(\xi)$, 
so $\pi(\imath_a(\pi(\xi)))$ is in the
right hand side, showing one inclusion, and the other is clear.
For any $b$ then
\begin{equation}
\pih_b^{-1}(C)=\pi_{0b}(\pi^{-1}(C)\cap\{\xi:\ |\xi|^2\leq R\}\cap
\sct_{C_b}\Xb)\cap\Sigma_b(\lambda),
\end{equation}
and $\pi_{0b}$ is continuous, so $\pih_b^{-1}(C)$ is also compact
(and similarly $\pi_b^{-1}(C)$ is closed). Thus, all the maps
$\pi_b$, $\pih_b$, are continuous.

Also, fix $\xib\in\sct_{C'_a}\Xb_a$, write $\xib=(\yb_a,\xib_a)$, and
choose a neighborhood $U_0$ of $\yb_a=p(\xib)$
in $\bXb$ such that $\overline{U_0}
\cap C_b=\emptyset$ unless $C_a\subset C_b$.
Let $\omega_{\xib}=\omega:\scdt_{\overline{U_0}} \Xb\to\Real$
be given by the following
$\pi$-invariant function on $\sct_{\overline{U_0}}\Xb$
(also denoted by $\omega$):
\begin{equation}
\omega(\xi)=(y_a-\yb_a)^2+|z_a|^2+(\xi_a-\xib_a)^2\geq 0
\end{equation}
in the coordinates $(y_a,z_a,\xi_a,\xi^a)$.
Suppose
that $U$ is open in $\dot\Sigma(\lambda)$, $\xib\in U$. Thus,
$K=(\dot\Sigma(\lambda)\cap\scdt_{\overline{U_0}}\Xb)\setminus U$
is compact, so, unless $K$ is empty, $\omega$ assumes a minimum on
it which thus has to be non-negative. But $\omega(\xi)=0$ implies
$z_a=0$, so $p(\xi)\in C_a$, and then $y_a=\yb_a$, $\xi_a=\xib_a$
show that $\xi=\xib$. Since $\xib\in U$, this shows that there exists
$\delta>0$ such that $\omega\geq \delta$ on $K$. Replacing $\delta>0$
by possibly a smaller number, we can also assume that $\omega(\xi)<\delta$
implies $p(\xi)\in U_0$. We thus conclude that
if $U$ is a neighborhood of $\xib$, then there exists $\delta>0$ such that
\begin{equation}
\{\xi\in\dot\Sigma(\lambda):\ \omega(\xi)<\delta\}\subset U.
\end{equation}
Of course, we could have equally well used $(\tau_a,\mu_a)$ in place of
$\xi_a$. These sets are open since $\omega$ is continuous, hence
they form a basis for the topology as $\xib$ and $\delta$ vary; it is
easy to see that if one restricts both of these to suitable countable
sets, one still has a basis. Note that, separating the complement of
$U$ from $\xib$ by a level set of $\omega$ shows explicitly
that $\dot\Sigma(\lambda)$ is regular, and a simple compactness argument
using these $\omega_{\xib}$ (composed with cut-off functions on the reals
as in the next paragraph) shows that $\dot\Sigma(\lambda)$ is normal,
hence it is a compact metrizable space.

Composing $\omega$ with a $\Cinf$ function on $\Real$ supported
near $0$ also shows that given any $\xib\in\dot\Sigma(\lambda)$ and
any neighborhood $U$ of $\xib$ in $\dot\Sigma(\lambda)$,
one can construct a $\pi$-invariant
$\Cinf$ function $f$
on $\sct_{\bXb} \Xb$ for which $f_\pi(\xib)\neq 0$ and $\supp f_\pi\subset U$.
This also shows the existence of smooth partitions of unity on
$\dot\Sigma(\lambda)$, smoothness understood as smoothness for the
pull-back to $\sct\Xb$.

If we define
\begin{equation}
\dot\Sigma(I)=\cup_{\lambda\in I}\dot\Sigma(\lambda)
\end{equation}
for $I\subset\Real$ compact, all the previous statements remain valid
with trivial modifications with $\dot\Sigma(I)$ in place of
$\dot\Sigma(\lambda)$. We write
\begin{equation}
\Sigma_b(I)=\cup_{\lambda\in I}\Sigma_b(\lambda),
\ \pih_{b,I}=\pi_b|_{\Sigma_b(I)}.
\end{equation}
Thus, $\dot\Sigma(I)=\cup_b\pi_b(\Sigma_b(I))$,
\begin{equation}\label{eq:pi_b_I-char}
\xit\in\pih_{b,I}^{-1}(\{\xi\})\Miff \pi_b(\xit)=\xi\Mand
\exists\sigma\in\pspec(H^b)\Mst
|\xit|^2+\sigma\in I.
\end{equation}

We end this section by recalling from \cite[Section~5]{Vasy:Propagation-Many}
the definition of the operator wave front set and the wave front set
of distributions in the many-body setting. Both of these are closed
subsets of $\scdt \Snp$ which were modelled on the fibred cusp wave front set
introduced by Mazzeo and Melrose in \cite{Mazzeo-Melrose:Fibred}.
Rather than using the invariant
definition given there in terms of oscillatory testing, we recall the one
given in terms of representation of ps.d.o.'s via quantization.

So suppose that $A\in\PsiSc^{-\infty,l}(\Snp,\calC)$ is the left
quantization of a symbol $a\in\rho_\partial^l\Cinf
([\Snp;\calC]\times\Snp)$ which vanishes to infinite order
at $[\Snp;\calC]\times\partial\Snp$. Then $\xi\nin\WFScp(A)$,
$\xi\in\sct_p(C_a;\Xb)$, $p\in C'_a$, if
and only if
there exists a neighborhood $U$ of $\xi$ in $\scdt\Snp$ such that
$a$ vanishes at $U'\subset(\partial[\Snp;\calC])\times\Rn$ to infinite
order where $U'$ is the
inverse image of $U$
under the composite map
\begin{equation}\begin{CD}\label{eq:WFScp-20}
\pit:(\partial[\Snp;\calC])\times\Rn@>\betaSc\times\identity>>
(\partial\Snp)\times\Rn=\sct_{\Sn}\Snp
@>\pi>>\scdt\Snp.
\end{CD}\end{equation}
It follows immediately from the usual formulae relating quantizations
that this definition is independent
of such choices. For example, we could have equally well written
$A$ as the right quantization of a symbol with similar properties.

The general definition for $A\in\PsiSc^{m,l}(\Xb,\calC)$,
again following the paper \cite{Mazzeo-Melrose:Fibred},
in the explicit quantization form as in the previous paragraph, would
also require the rapid decay of $a$ in an open cone (conic in the cotangent
variable, $\xit$, i.e.\ in the second factor, $\Rn$, in \eqref{eq:WFScp-20})
that includes $U'$. For $A\in\PsiSc^{-\infty,l}(\Xb,\calC)$,
$a$ is rapidly decreasing in every direction as $|\xi|\to\infty$,
so this statement is vacuous, and we recover the above definition.
The main point is that if $A\in\PsiSc^{0,0}(\Xb,\calC)$,
$\Ah_a(\xi)$ is invertible, then there exists a microlocal parametrix
for $A$, i.e.\ there exists $G\in\PsiSc^{0,0}(\Xb,\calC)$ such that
$\Id=AG+R_R$, $\Id=GA+R_L$, with $R_R,R_L\in\PsiSc^{0,0}(\Xb,\calC)$,
$\xi\nin\WFScp(R_R)$, $\xi\nin\WFScp(R_L)$; see
\cite[Lemmas~14-15]{Mazzeo-Melrose:Fibred}.

We recall that the many-body wave front set $\WFSc(u)$ of
$u\in\dist(\Xb)$, defined in \cite[Section~5]{Vasy:Propagation-Many},
is a subset of $\scdt \Xb$, defined as follows. The definition is
somewhat delicate due to the relation among the principal symbols
$\prs_{\Scl,0}(\Ah_a(\xi))$ as $\xi$ varies, see \cite{Vasy:Propagation-Many}
for a detailed discussion, but it becomes much simpler for approximate
generalized eigenfunctions of $H$ as the proposition following the definition
shows. Although the definition is complicated, it ensures that all the desired
properties hold: for all $u\in\dist(X)$, $\WFSc(u)$ is closed in $\scdt \Xb$,
\begin{equation}\begin{split}\label{eq:WF-props}
&u_1,u_2\in\dist(\Xb)\Rightarrow
\WFSc(u_1+u_2)\subset\WFSc(u_1)\cup\WFSc(u_2),\\
&u\in\dist(\Xb),\ A\in\PsiSc^{m,l}(\Xb,\calC)\Rightarrow\WFSc(Au)
\subset\WFSc(u),\\
&A\in\PsiSc^{-\infty,l}(\Xb,\calC),\ u\in\dist(\Xb),
\ \WFScp(A)\ \text{compact},
\ \WFScp(A)\cap\WFSc(u)=\emptyset\\
&\qquad\qquad\Rightarrow Au\in\dCinf(\Xb).
\end{split}\end{equation}
The central one of these means that `$\PsiSc^{*,*}(\Xb,\calC)$
is microlocal'.

\begin{Def}\label{Def:WFSc}(\cite[Definition~5.2]{Vasy:Propagation-Many})
We say that
\begin{equation}\begin{split}
\xi\nin\WFSc(u)\cap\sct_{C'_a} (C_a;\Xb)\Miff
&\exists A\in\PsiSc^{0,0}(\Xb,\calC),
\ \Ah_{a,0}(\xi)\ \text{invertible in}\ \PsiSc^{0,0}(\Xb^a,\calC^a),\\
&\exists B_j\in\PsiSc^{-\infty,0}(\Xb,\calC),\ \xi\nin\WFScp(B_j),\\
&\exists u_j\in\dist(\Xb),\ j=1,\ldots,s,\ f\in\dCinf(\Xb),\\
&Au=\sum_{j=1}^s B_j u_j+f.
\end{split}\end{equation}
Similarly, the filtered version of the $\Scl$-wave front set is given by
\begin{equation}\begin{split}
\xi\nin\WFSc^{m,l}(u)\cap\sct_{C'_a} (C_a;\Xb)\Miff
&\exists A\in\PsiSc^{0,0}(\Xb,\calC),
\ \Ah_{a,0}(\xi)\ \text{invertible in}\ \PsiSc^{0,0}(\Xb^a,\calC^a),\\
&\exists B_j\in\PsiSc^{-\infty,0}(\Xb,\calC),\ \xi\nin\WFScp(B_j),\\
&\exists u_j\in\dist(\Xb),\ j=1,\ldots,s,\ f\in\Hsc^{m,l}(\Xb),\\
&Au=\sum_{j=1}^s B_j u_j+f.
\end{split}\end{equation}
\end{Def}

The wave front set has the following property
which is useful
for approximate generalized eigenfunctions of $H$.

\begin{prop}\label{prop:WF-3}(\cite[Proposition~5.5]{Vasy:Propagation-Many})
Suppose that $u\in\dist(\Xb)$, $\lambda\in\Real$, and define
$W\subset\scdt \Xb$ by
\begin{equation}\begin{split}
\xi\nin W\cap\sct_{C'_a}(C_a;\Xb)
\Miff&\exists\psi\in\Cinf_c(\Real),\ \psi(\lambda)=1,\\
&\exists A\in\PsiSc^{-\infty,0}(\Xb,\calC),
\ \Ah_a(\xi)=\widehat{\psi(H)}_a(\xi),
Au\in\dCinf(\Xb).
\end{split}\end{equation}
Then
\begin{equation}
\WFSc(u)\subset\WFSc((H-\lambda)u)\cup W.
\end{equation}
The same conclusion holds with $\WFSc$ replaced by $\WFSc^{m,l}$ and
$Au\in\dCinf(\Xb)$ by $Au\in\Hsc^{m,l}(\Xb)$.
\end{prop}

\begin{cor}\label{cor:micro-ell-reg} (Microlocal elliptic regularity.)
Suppose that $u\in\dist(\Xb)$, $\lambda\in\Real$. Then $\WFSc(u)\subset
\WFSc((H-\lambda)u)\cup\dot\Sigma(\lambda)$.
\end{cor}

\begin{proof}
If $\xib_a\nin\dot\Sigma(\lambda)$,
$\xib_a\in\sct_{C'_a}\Xb_a$, one can choose $\psi\in\Cinf_c(\Real)$
such that $\psi(\lambda)=1$ and $\widehat{\psi(H)}_a(\xib_a)
=\psi(H^a+|\xib_a|^2)=0$ by making the support of $\psi$ sufficiently
small. Then we can take $A$ to be the zero operator, showing that
$\xib_a\nin W$. Thus, Proposition~\ref{prop:WF-3} proves the
corollary.
\end{proof}

\section{Broken bicharacteristics}\label{sec:br-bichar}
Since at $C'_a$ a particle can be in a bound state of any of the
clusters  $C_b$ with
$C_a\subset C_b$, it can be expected that the Hamilton vector fields
associated to the $b$-external kinetic energy, i.e.\ to $\Delta_{X_b}$,
govern the propagation of singularities. If the particle is in a
$b$-bound state of energy $\ep_\beta$, its $b$-external kinetic energy is
$\lambda-\ep_\beta$.
The
symbol of $\Delta_{X_b}-(\lambda-\ep_\beta)$ at $C_b=\partial\Xb_b$ (i.e.\ its
$\scl$-indicial operator) is $|\xi_b|^2+\ep_\beta-\lambda$. The Hamilton
vector field of $g_b=|\xi_b|^2$ (or $|\xi_b|^2+\ep_\beta-\lambda$)
on $T^*X_b$ is
\begin{equation}
H_g^b=2\xi_b\cdot\partial_{w_b}.
\end{equation}
Following the general principle that it is more convenient to do analysis
on compact spaces than analyzing uniform properties in a non-compact
setting, we rescale $H_g^b$ and regard the result as a smooth vector field
on $\sct\Xb_b$, i.e.\ as an element of $\Vf(\sct \Xb_b)$.
Thus, the rescaled Hamilton vector field
\begin{equation}
\scHg^b=x_b^{-1}H_g^b
\end{equation}
of $g_b=|\xi_b|^2$, introduced in \cite{RBMSpec}, is
\begin{equation}\begin{split}
\scHg^b&=2|w_b| \xi_b\cdot\partial_{w_b}\\
&=2\tau_b(x_b\partial_{x_b}+\mu_b\cdot\partial_{\mu_b})-2|\mu_b|^2
\partial_{\tau_b}+H_{|\mu_b|^2}+x_b W',\quad W'\in\Vb(\sct \Xb_b),
\end{split}\end{equation}
so its restriction to $C_b$, also denoted by $\scHg^b$, is
\begin{equation}
\scHg^b=2\tau_b \mu_b\cdot\partial_{\mu_b}-2|\mu_b|^2
\partial_{\tau_b}+H_{|\mu_b|^2}.
\end{equation}
Thus, for $C_a\subset C_b$,
\begin{equation}
\scHg^b=2x_b^{-1} (\xi_a\cdot\partial_{w_a}+\xi_{ab}\cdot\partial_{w_{ab}})
=(x_a/x_b)(\scHg^a+2\xi_{ab}\cdot\partial_{w_{ab}}).
\end{equation}
If $b=0$, then in agreement with \cite{RBMSpec} we write
\begin{equation}
\scHg=\scHg^0.
\end{equation}

Before defining generalized
broken bicharacteristics, we discuss how the rescaled
Hamilton vector fields are related for various $b$.
First, note that
$\sct_{C_b}\Xb$ is the bundle direct sum of $\sct_{C_b}\Xb_b$ and
the annihilator of $\Tsc_{C_b}(C_b;\Xb)$. Let
\begin{equation}
\imath_b:\sct_{C_b}\Xb_b\hookrightarrow\sct_{C_b} \Xb
\end{equation}
be the inclusion map as the zero section in the second summand, so
\begin{equation}
\imath_b(y_b,\tau_b,\mu_b)=(y_b,0,\tau_b,\mu_b,0),
\end{equation}
i.e.\ $\nu_b(\imath_b(\xi))=0$ for all $\xi\in\sct_{C_b}\Xb_b$. Then
\begin{equation}\label{eq:scHg-scHg_b}
(\imath_b)_*\scHg^b=\scHg,
\end{equation}
since the corresponding statement for the non-rescaled Hamilton vector
fields is easy to see and $x_b/x$ is $\Cinf$ near $C_b$, equal to $1$ at
$C_b$.

Next, note that
\begin{equation}\label{eq:T^*X_a-X^a}
T^* \Rn=T^*X_a\times T^*X^a.
\end{equation}
Thus, if $f$ is a smooth function on $\sct_{C_a}\Xb_a$, we can extend it to
a smooth function on $\sct\Xb_a$ using polar coordinates (i.e.\ so that
it is independent of $x_a$, at least near $C_a$),
and pull-back the result to $T^*\Rn$
using the projection to the first factor. The result is then smooth
on $\sct \Snp$ in a neighborhood of $C_a$ (in fact, everywhere away from
the closure of $X^a$ in $\Snp$). Thus, we have a pull-back, which
we denote as $F=(\pi_a^e)^*f$. Note that $F=(\pi_a^e)^*f$ is $\pi$-invariant
near $C'_a$. Moreover,
\begin{equation}\label{eq:pi_a^e-scHg}
\scHg F=\frac{x_a}{x}\,(\pi_a^e)^*\scHg^a f.
\end{equation}
One often drops the pull-back notation, and simply writes $f$ in place of $F$.
A similar formula also holds for $\scHg^b f$ by \eqref{eq:scHg-scHg_b}.

It is worth calculating the derivatives of two functions in particular
along the rescaled Hamilton vector fields. First,
\begin{equation}
\scHg\tau=-2|\mu|^2,
\end{equation}
so by \eqref{eq:scHg-scHg_b},
\begin{equation}
\scHg^b\tau=-2|\mu_b|^2.
\end{equation}
Next, consider
\begin{equation}
\eta_a=z_a\cdot\nu_a=\frac{\xi^a\cdot w^a}{|w_a|}.
\end{equation}
Then
\begin{equation}
\scHg\eta_a=\frac{x_a}{x}(2\tau_a\eta_a+2|\xi^a|^2),
\end{equation}
so
\begin{equation}\label{eq:scHg^b-eta_a}
\scHg^b\eta_a=\frac{x_a}{x}(2\tau_a\eta_a+2|\xi^a_b|^2)
=\frac{x_a}{x}(2\tau_a\eta_a+2|\nu_{ab}|^2).
\end{equation}

We define generalized broken
bicharacteristics of $H-\lambda$ partly following Lebeau
\cite{Lebeau:Propagation}, but in such a way
that all the analytic properties (compactness, applicability of
positive commutator estimates) will be clear. However, the geometric
properties will be less apparent, and we will devote some time to
clarifying these.

First, we say that a function
$f\in\Cinf(\sct_{\bXb} \Xb)$ is $\pi$-invariant if for
$\xit,\xit'\in\sct_{\bXb} \Xb$, $\pi(\xit)=\pi(\xit')$ implies $f(\xit)=f(\xit')$.
A $\pi$-invariant function $f$ naturally defines a function $f_\pi$ on
$\scdt \Xb$ by $f_\pi(\xi)=f(\xit)$ where $\xit\in\sct_{\bXb}\Xb$ is chosen so
that $\pi(\xit)=\xi$.

We also need to introduce notation for the
upper and lower one-sided derivatives of functions
$f$ defined on an interval $I$:
\begin{equation}\begin{split}\label{eq:one-sided-derivs}
(D_+ f)(t_0)&=\liminf_{t\to t_0+}\frac{f(t)-f(t_0)}{t-t_0},\\
(D_- f)(t_0)&=\liminf_{t\to t_0-}\frac{f(t)-f(t_0)}{t-t_0},\\
(D^+ f)(t_0)&=\limsup_{t\to t_0+}\frac{f(t)-f(t_0)}{t-t_0},\\
(D^- f)(t_0)&=\limsup_{t\to t_0-}\frac{f(t)-f(t_0)}{t-t_0}.
\end{split}\end{equation}
A function $f$ is thus differentiable at $t_0$ if and only if these four
assume the same (finite) value; then we write these as $Df(t_0)$.
Similarly, $f$ is differentiable from the left at $t_0$ if and only if
$(D_- f)(t_0)=(D^- f)(t_0)$ is finite; we write the common value as
$Df(t_0-)$.

\begin{Def*}(Definition~\ref{Def:gen-br-bichar})
A generalized broken bicharacteristic of $H-\lambda$
is a continuous map
$\gamma:I\to\dot\Sigma(\lambda)$,
where $I\subset\Real$ is an interval, such that for all $t_0\in I$
and for each sign $+$ and $-$ the following holds. Let
$\xi_0=\gamma(t_0)$, suppose that $\xi_0\in\sct_{C'_a}\Xb_a$. Then for
all $\pi$-invariant functions
$f\in\Cinf(\sct_{\bXb}\Xb)$,
\begin{equation}
D_\pm(f_\pi\circ \gamma)(t_0)
\geq\inf\{\scHg^b f(\xit_0):\ \xit_0\in\pih_{b}^{-1}(\xi_0),\ C_a\subset C_b\}.
\end{equation}
In this section we usually drop the word `generalized' for the sake of
brevity.
\end{Def*}

\begin{rem}
Considering $-f$ in place of $f$, the definition immediately
gives a similar estimate for $D^\pm(f_\pi\circ \gamma)(t_0)$. In addition,
the function $f$ only has to be defined on $\sct_U \Xb$ where $U$ is
a neighborhood of $p$ in $\bXb$, $\gamma(t_0)\in\sct_p (C_a;\Xb)$. Indeed,
otherwise one can consider $\phi f$, where $\phi$ is (the pull-back of)
a cut-off function supported in $U$ that is identically $1$ near $p$,
without affecting any of the statements above.
\end{rem}

\begin{rem}\label{rem:tgt-smoothness}
The local coordinates
$y_a$ and $\xi_a$ are actually differentiable along $\gamma$
at $t_0$ since the Hamilton vector fields $\scHg^b$ applied to them give
$\pi$-invariant results which agree with each other by
\eqref{eq:scHg-scHg_b} and \eqref{eq:pi_a^e-scHg}.
In addition, over $C'_0$, $\pih$ is one-to-one,
so if $J$ is an interval and $\gamma|_J$ lies
in $\sct_{C'_0}\Xb$ then $\gamma|_J$ is a bicharacteristic
of $\scHg$.
\end{rem}

Moreover, since $\cup_b\pih^{-1}_b(\dot\Sigma(\lambda))$ is
compact, it follows that for any $\pi$-invariant function $f$, the
derivatives
$D_\pm(f_\pi\circ \gamma)$,
$D^\pm(f_\pi\circ \gamma)$ are bounded (independently of $t_0$, and
even of $\gamma$), hence $f_\pi\circ \gamma$ is
Lipschitz; in fact, uniformly so (i.e.\ the Lipschitz constant depends on $f$,
but not on $\gamma$).
This at once implies the equicontinuity of the set of broken
bicharacteristics as in \cite{Lebeau:Propagation}.

\begin{lemma}(cf.\ Lebeau's proof, \cite[Proposition~6]{Lebeau:Propagation})
For any compact interval $I=[T_1,T_2]$, the set
$\calR$ of broken bicharacteristics $\gamma:I\to\dot\Sigma(\lambda)$ is
equicontinuous.
\end{lemma}

\begin{proof}
Let $d$ be a metric defining the topology on $\dot\Sigma(\lambda)$.
Suppose $\calR$ is not equicontinuous. That is, suppose that
there exists $\ep_0>0$,
sequences $s_k$, $s'_k$ in $I$, $\gamma_k\in\calR$, such that
$|s_k-s'_k|\leq 1/k$ but $d(\gamma_k(s_k),\gamma_k(s'_k))\geq\ep_0$.
Since $\dot\Sigma(\lambda)$ is compact, one can pass to subsequences (which
we do not show in the notation) such that
$\gamma_k(s_k)$ and $\gamma_k(s'_k)$ converge to points $\xi$ and
$\xi'$ respectively. Note that $\xi\neq\xi'$ since $d(\xi,\xi')\geq\ep_0$.
For any $\pi$-invariant function $f$,
$f_\pi\circ\gamma_k$ is uniformly Lipschitz (i.e.\ the Lipschitz constant
is independent of $k$), so $|s_k-s'_k|\leq 1/k$ shows that
$|f_\pi(\gamma_k(s_k))-f_\pi(\gamma_k(s'_k))|\leq M/k$;
$M$ independent of $k$. But $f_\pi$ is continuous, so
$\lim_{k\to\infty}f_\pi(\gamma_k(s_k))=f_\pi(\xi)$,
$\lim_{k\to\infty}f_\pi(\gamma_k(s'_k))=f_\pi(\xi')$, so we conclude
$f_\pi(\xi)=f_\pi(\xi')$. But let $a$ be such that
$\xi\in\sct_p\Xb_a$, $p\in C'_a$. All functions on $\Xb$, pulled back by
the bundle projection to $\sct \Xb$, are $\pi$-invariant, so we see that
$\xi'\in\sct_p\Xb_a$ as well. But $(\xi_a)_j$ is
$\pi$-invariant near $\xi$ for all $j$,
so this shows $\xi=\xi'$, a contradiction.
\end{proof}

Next, we note that the uniform limit of broken bicharacteristics is a broken
bicharacteristic.

\begin{prop}
Suppose that $I$ is a compact interval $I=[T_1,T_2]$,
$\gamma_k:I\to\dot\Sigma(\lambda)$
are broken bicharacteristics, $\gamma:I\to\dot\Sigma(\lambda)$,
and $\gamma_k\to\gamma$ uniformly on $I$.
Then $\gamma$ is a broken bicharacteristic.
\end{prop}

\begin{proof}
Let $f$ be a $\pi$-invariant function, $t_0\in I$,
$\xi_0=\gamma(t_0)$, and let
\begin{equation}
c_0=\inf\{
\scHg^b f(\xit_0):\ \xit_0\in\pih_b^{-1}(\xi_0):\ C_a\subset C_b\}.
\end{equation}
We need to show that $D_\pm (f_\pi\circ\gamma)(t_0)\geq c_0$.
We only consider $D_+$ for the sake of definiteness.
So let $\ep>0$; we need to prove that there exists $\delta>0$ such
that for all $t\in (t_0,t_0+\delta)$,
$f_\pi\circ\gamma(t)-f_\pi\circ\gamma(t_0)\geq (c_0-\ep)(t-t_0)$.

But by the continuity of $\scHg^b f$ on $\sct_{C_b}\Xb_b$, there exists
a neighborhood $U$ of $\xi_0$ in $\dot\Sigma(\lambda)$ such that
$\scHg^b f(\xit)\geq c_0-\ep/3$ for all $b$, $\xit\in\pih_b^{-1}(\xi)$,
$\xi\in U$. Next, by the uniform convergence of the $\gamma_k$,
there exist $\delta>0$ and $M\in\Nat$
such that for $t\in(t_0-\delta,t_0+\delta)$, $k\geq M$,
$\gamma_k(t)\in U$. Let
\begin{equation}
F_k:I\to\Real,\quad F_k(t)=f_\pi\circ\gamma_k(t)-
(c_0-\ep/3)t.
\end{equation}
Thus, for $t\in(t_0,t_0+\delta)$, $k\geq M$, $F_k$ satisfies 
$D_+F_k\geq 0$, hence
$F_k$ is non-decreasing, so
for $t\in(t_0,t_0+\delta)$, $k\geq M$,
\begin{equation}\label{eq:unif-bich-3}
f_\pi\circ\gamma_k(t)-f_\pi\circ\gamma_k(t_0)\geq (c_0-\ep/3)(t-t_0).
\end{equation}
Now given $t\in(t_0,t_0+\delta)$, simply choose $k\geq M$ such that
$|f_\pi\circ\gamma_k(t)-f_\pi\circ\gamma(t)|<\ep (t-t_0)/3$,
$|f_\pi\circ\gamma_k(t_0)-f_\pi\circ\gamma(t_0)|<\ep (t-t_0)/3$,
which is possible since $f_\pi$ is continuous, so $f_\pi\circ\gamma_k(t')
\to f_\pi\circ\gamma(t')$ for all $t'$. Thus, using this particular
value of $k$ and \eqref{eq:unif-bich-3}, we conclude that
\begin{equation}
f_\pi\circ\gamma(t)-f_\pi\circ\gamma(t_0)\geq (c_0-\ep)(t-t_0).
\end{equation}
This holds for every $t\in (t_0,t_0+\delta)$, completing the proof.
\end{proof}

We thus deduce the compactness of the set of broken bicharacteristics.

\begin{prop}\label{prop:Lebeau-compactness}
For any compact interval $I=[T_1,T_2]$, the set
$\calR$ of broken bicharacteristics $\gamma:I\to\dot\Sigma(\lambda)$ is
compact in the topology of uniform convergence.
\end{prop}

\begin{proof}
Since $\calR$ is equicontinuous and $\dot\Sigma(\lambda)$ is
compact, any sequence of broken bicharacteristics
in $\calR$ has a convergent subsequence, converging uniformly over $I$,
by the theorem of Ascoli-Arzel\'a. But the uniform limit of broken
bicharacteristics is a broken bicharacteristic, proving the
proposition.
\end{proof}

\begin{cor}(cf.\ Lebeau, \cite[Corollaire~7]{Lebeau:Propagation})
\label{cor:Lebeau-extension}
If $\gamma:(T_1,T_2)\to\dot\Sigma(\lambda)$
is a generalized broken bicharacteristic then
$\gamma$ extends to a generalized broken bicharacteristic on $[T_1,T_2]$.
\end{cor}

\begin{proof}
Let $\gamma_n=[T_1,T_1+\delta]\to\dot\Sigma(\lambda)$ be given by
$\gamma_n(t)=\gamma(t+1/n)$, where $\delta>0$ is chosen sufficiently small.
By the previous proposition, there exists a subsequence of $\{\gamma_n\}$
that converges to a broken bicharacteristic $\gammat:[T_1,T_1+\delta]\to
\dot\Sigma(\lambda)$. But then $\gammat|_{(T_1,T_1+\delta]}=\gamma$, so
$\gammat$ gives the desired extension of $\gamma$ to $[T_1,T_2)$. The other
endpoint can be dealt with similarly.
\end{proof}

\begin{rem}\label{rem:ext-to-Real}
This corollary shows that every generalized broken bicharacteristic
$\gamma_0:J\to\dot\Sigma(\lambda)$
can be extended to another one, $\gamma$, defined over $\Real$
(meaning that $\gamma:\Real\to\dot\Sigma(\lambda)$ and
$\gamma|_J=\gamma_0$). 
Indeed, by the corollary we may assume that $J$ is closed. It suffices
to show that, say, $\gamma$ can be extended to an interval that is
unbounded from below; the extension in the other direction is similar.
So suppose the lower endpoint of $J$ is $T_1\in\Real$ (if it is $-\infty$,
we are done), and suppose that $\gamma_0(T_1)\in\sct_{C'_a}\Xb_a\cap
\pi_{ba}(\Sigma_b(\lambda))$; such a $b$ exists by \eqref{eq:dSa-as-proj-Sb}.
Then choose $\xit\in\Sigma_b(\lambda)$ such that $\pi_{ba}(\xit)=\gamma_0
(T_1)$, and define $\gamma|_{(-\infty,T_1]\cup J}$ by
$\gamma(t)=\pih(\exp((t-T_1)\scHg^b)(\xit))$ for $t\leq T_1$,
$\gamma(t)=\gamma_0(t)$ for $t\in J$. Directly from
Definition~\eqref{Def:gen-br-bichar}, $\gamma$ is a generalized broken
bicharacteristic with the desired properties.
\end{rem}

Since Lipschitz functions are differentiable almost everywhere, and they
are equal to the integral of their a.e.\ defined derivative, we can analyze
$\pi$-invariant functions such as $\tau$ in more detail.

In fact, the Hamilton vector field
$\scHg^b$ applied to $-\tau$ gives $2|\mu_b|^2\geq 0$, so we deduce that
$-\tau$ is monotone increasing, i.e.\ $\tau$ is monotone decreasing
along broken bicharacteristics. In fact, writing $F=-\tau\circ\gamma$,
note that $\xit=\pih^{-1}_b(\xi)$ being in $\Sigma_b(\lambda)$ gives
\begin{equation}
\lambda-|\xit_b|^2=\lambda-\tau^2-|\mu_b|^2\in\pspec(H^b).
\end{equation}
Thus, $\lambda-\tau^2(\gamma(t_0))\nin\pspec(H^b)$ for any $b$
with $C_b\supset C_a$
implies that
\begin{equation}
D_\pm F(t_0)\geq c_0>0,
\end{equation}
so there exist $c>0$ and $\delta>0$ such that
\begin{equation}
|\tau(\gamma(t))-\tau(\gamma(t_0))|\geq c|t-t_0|\Mif |t-t_0|<\delta,
\end{equation}
so the bicharacteristic $\gamma$ is not constant.
Conversely, if there exists $b$ with $C_b\supset C_a$ such that
$\lambda-\tau^2(\gamma(t_0))\in\pspec(H^b)$, then the constant curve
\begin{equation}
\gamma_0:\Real\to\scdt \Xb,\ \gamma_0(t)=\gamma(t_0)\ \text{for all}\ t\in
\Real,
\end{equation}
is a broken bicharacteristic since for $\xit\in\pih_b^{-1}(\gamma(t_0))$,
we have $\mu_b=0$, so $\scHg^b(\xit)$ vanishes.

Thus, we define the
$b$-radial sets
\begin{equation}
R_\pm^b(\lambda)=\{\xi_b\in\Sigma_b(\lambda):\ \mu_b=0,\ \pm\tau\geq 0\},
\end{equation}
and their image by
\begin{equation}
R_\pm(\lambda)=\cup_b\pi_b(R_\pm^b(\lambda)).
\end{equation}
We cannot expect a non-trivial propagation theorem at points in
$R_\pm(\lambda)$
(since one of the broken bicharacteristics through them is stationary), but
we can expect such results elsewhere. We emphasize that there can be (and
usually there are, though not if there are no bound states in the
subsystems) non-constant bicharacteristics through $R_\pm(\lambda)$; it is
because there exists a constant bicharacteristic that we cannot expect an
interesting propagation result.

Note that if over an interval $J\subset I$, the image of
$\gamma$ is disjoint from $R_+(\lambda)\cup R_-(\lambda)$, then $\gamma|_J$
can be reparameterized using $\tau$ as a parameter instead of $t$;
Lipschitz functions on the
reparameterized curve remain Lipschitz.

Another interesting $\pi$-invariant function is $\eta=
\eta_a=z_a\cdot\nu_a$ defined
near $\sct_p(C_a;\Xb)$ for some $p\in C'_a$. Since $\scHg^b\eta$ over
$\Sigma_b(\lambda)$ is
\begin{equation}\label{eq:scHg^b-eta_a-3}
\scHg^b\eta_a=\frac{x_a}{x}(2\tau_a\eta_a+2|\nu_{ab}|^2),
\end{equation}
where $\mu_b=(\mu_a,\nu_{ab})$ is the decomposition corresponding to that
of $X_b$ with respect to $X_a$ and its orthocomplement, so
\begin{equation}\label{eq:nu-mu-a-b}
|\mu_b|^2=|\mu_a|^2+|\nu_{ab}|^2,
\end{equation}
we can obtain a monotonicity result for $F=\eta_a\circ\gamma$,
defined for $t$ near $t_0$. First note that $F(t_0)=0$ (as $z=0$ at
$\gamma(t_0)\in\sct_{C'_a}(C_a;\Xb)$). Next, $F(t)\geq 0$ for $t\geq t_0$.
Indeed, if $F(t_1)<0$, let $t'_0=\sup\{t:\ F(t)=0,\ t<t_1\}$,
so $t_0\leq t_0'<t_1$; since
$F$ is continuous, $F(t)<0$ for $t\in(t'_0,t_1]$. But by
\eqref{eq:scHg^b-eta_a-3}, $D_\pm F(t)\geq cF(t)$ for some $c>0$
whenever $F\leq 0$, hence on $[t_0',t_1]$, so $D_\pm (e^{-c.}F)\geq 0$
there, so $e^{-c.}F$ is non-decreasing, so $F(t_1)\geq F(t'_0)=0$,
a contradiction. Thus, $F(t)\geq 0$ for $t\geq t_0$, so by
\eqref{eq:scHg^b-eta_a-3}, $D_\pm F(t)\geq -c'F(t)$ for some $c'>0$, so
$D_\pm (e^{c'.}F)(t)\geq 0$, hence $e^{c't}F(t)$ is non-decreasing.
Thus, with $\eta_a(t)=\eta_a(\gamma(t))=F(t)$, there
exists $\delta>0$ such that either $\eta_a$ is
identically $0$ on $[t_0,t_0+\delta)$, or $\eta_a(t)>0$ for
$t\in(t_0,t_0+\delta)$.
Since $\scHg^b |z_a|^2=4(x_a/x)\eta_a$ (as $\eta_a$ is $\pi$-invariant),
we deduce that
there exists $\delta>0$ such that
either $z_a=0$ for all $t$ in $(t_0,t_0+\delta)$, or
$z_a\neq 0$ for all $t$ in $(t_0,t_0+\delta)$. If the former holds,
then $\gamma|_{(t_0,t_0+\delta)}$ is a differentiable
curve in $\sct_{C'_a}\Xb_a$ with tangent vector given by $\scHg^a$,
by Remark~\ref{rem:tgt-smoothness}.
In other words, we have proved
the following lemma.

\begin{lemma}\label{lemma:br-bichar-tgt-normal}
Suppose $\gamma$ is a broken bicharacteristic.
Then there exists $\delta>0$ such that either
$\gamma(t)\in\sct_{C'_a}(C_a;\Xb)$ for all $t\in(t_0,t_0+\delta)$, or
$\gamma(t)\nin\sct_{C'_a}(C_a;\Xb)$ for all $t\in(t_0,t_0+\delta)$. If the
former holds, then $\gamma|_{(t_0,t_0+\delta)}$ is an integral curve
of $\scHg^a$.
A similar result holds for $(t_0-\delta,t_0)$. In the former case we
say that $\gamma$ is forward/backward tangential at $t_0$, in the latter
that $\gamma$ is forward/backward normal at $t_0$.
\end{lemma}

\begin{rem}
Note that tangential/normal behavior is the property of $\gamma$, not just of
its value at $t_0$, unlike in the situation when there are no bound
states in any subsystems. Indeed, in that case we only need to
consider $\scHg f(\xit_0)$ for $\xit_0\in\pih_0(\xi_0)$ in
Definition~\ref{Def:gen-br-bichar},
so in particular $\scHg \eta_a(\xit_0)=|\nu_a|^2=\lambda-\tau_a^2-|\mu_a|^2$,
so $\gamma$ is tangential if and only
if $\tau_a^2+|\mu_a|^2=\lambda$ at $\gamma(t_0)\in\sct_{C'_a}\Xb_a$.
This remark also shows that our definition of broken bicharacteristics
agrees with that of \cite{Vasy:Propagation-Many}, which in turn
corresponds to the analogous definition introduced by Lebeau in the
study of the wave equation \cite{Lebeau:Propagation}. In particular,
our general definition is equivalent to the one discussed in the
introduction if $\Lambda_1=\{0\}$.
\end{rem}

We can give a more geometric description of the broken
bicharacteristics, provided that the set of thresholds, $\Lambda_1$,
is discrete, or if $H$ is a four-body Hamiltonian, $\lambda\nin\Lambda_1$.
Since this is improvement is not important for the
propagation estimates, we only state the result below, and give the
proof in the Appendix. Note that if there are no bound states
in any subsystems, then $\Lambda_1=\{0\}$ is certainly discrete, so
in particular the proof applies in the setting of
\cite{Vasy:Propagation-Many} (though it was already proved there, with the
proof based on Lebeau's results).

\begin{thm}\label{thm:geom-br-bichar}
Suppose that $\Lambda_1$ is discrete and
$\gamma:\Real\to\dot\Sigma(\lambda)$
is a continuous curve. Then $\gamma$
is a generalized
broken bicharacteristic of $H-\lambda$ if and only if there exist
$t_0<t_1<t_2<\ldots<t_k$ such that $\gamma|_{[t_j,t_{j+1}]}$,
as well as $\gamma|_{(-\infty,t_0]}$ and $\gamma_{[t_k,+\infty)}$, are
the projections of integral curves
of the Hamilton vector field $\scHg^a$ for some $a$. Similar results
hold if the interval of definition, $\Real$, is replaced by any interval.
\end{thm}

We end this section by analyzing the behavior of generalized broken
bicharacteristics under the assumption that $\lambda\nin\Lambda_1$;
this will be useful when studying the resolvent.

\begin{lemma}\label{lemma:t-to-infty}
Suppose $H$ is a many-body Hamiltonian, $\lambda\nin\Lambda_1$, and
$\gamma:I\to\dot\Sigma(\lambda)$ is a generalized broken
bicharacteristic where $I$ is a closed unbounded interval. Then
\begin{equation}
\overline{\gamma(I)}\subset\gamma(I)\cup
R_+(\lambda)\cup R_-(\lambda).
\end{equation}
Moreover, if $I$ is not bounded above,
$\tau(+\infty)=\lim_{t\to+\infty}\tau(\gamma(t))$
exists, $\lambda-\tau(+\infty)^2\in\Lambda_1$, and there
exists $\xi\in R_+(\lambda)\cup R_-(\lambda)$ and a sequence
$t_j\to+\infty$ such that $\xi=\lim_{j\to\infty}\gamma(t_j)$ in
$\dot\Sigma(\lambda)$ (hence in particular $\tau(\xi)=
\tau(+\infty)$). Similar results hold as $t\to-\infty$ if $I$ is
not bounded below.
\end{lemma}

\begin{proof}
For the sake of definiteness we take $I=[t_0,+\infty)$, all other cases
are very similar.

Since $\tau\circ\gamma$ is monotone decreasing, and $\tau$ is a bounded
function on $\dot\Sigma(\lambda)$, $\tau(+\infty)=
\lim_{t\to+\infty}\tau(t)$ exists
(we wrote $\tau(t)=\tau(\gamma(t))$).

Now, let $t_n\to +\infty$ be any sequence such that $t_n\geq t_0$
for all $n$. Let $\gamma_n:[0,1]\to\dot\Sigma
(\lambda)$ given by $\gamma_n(t)=\gamma(t+t_n)$.
By
Proposition~\ref{prop:Lebeau-compactness}, $\{\gamma_n\}$ has a subsequence,
which we write as $\gamma'_n=\gamma_{k(n)}$, which converges uniformly to
a generalized broken bicharacteristic $\gamma':[0,1]\to\dot\Sigma(\lambda)$.
Thus, for $t\in[0,1]$, $\tau(\gamma'(t))=\lim_{n\to\infty}\tau
(\gamma'_n(t))=\lim_{n\to\infty}\gamma(t+t_{k(n)})=\tau(+\infty)$,
so $\tau$ is constant on $\gamma'$. Since $\tau\circ\gamma'$ is Lipschitz,
it is equal to the integral of its a.e.\ defined derivative, $-2|\mu(t)|^2$,
so we conclude that $|\mu(t)|^2=0$ a.e.\ along $\gamma'$. But then
the arclength of the projection of $\gamma'$ to $\Sn$, which is the integral
of $|\mu(t)|$, is also zero, so $\gamma'$ is a constant curve, and
hence, as shown above, there exists $\xi\in R_+(\lambda)\cup
R_-(\lambda)$ such that $\gamma'(t)=\xi$ for all $t\in[0,1]$.
Thus,
\begin{equation}
\overline{\gamma([t_0,+\infty))}\cap(R_+(\lambda)\cup R_-(\lambda))
\neq\emptyset.
\end{equation}
Also, $\tau(+\infty)=\tau(\xi)$, and $\xi\in R_+(\lambda)\cup R_-(\lambda)$,
so $\lambda-\tau(+\infty)^2\in\Lambda_1$.

We still need to show that if $\{t_n\}$ is a sequence in $[t_0,+\infty)$,
and $\{\gamma(t_n)\}$ converges in $\dot\Sigma(\lambda)$,
say $\xi_0=\lim_{n\to\infty}\gamma(t_n)$, then either $\xi\in\gamma([t_0,
+\infty))$, or $\xi\in R_+(\lambda)\cup R_-(\lambda)$. Now, if
$t_n$ has a subsequence converging to some $T\in[t_0,+\infty)$,
then $\gamma(t_n)$ converges to $\gamma(T)$, so we are done. We may thus
assume that $t_n\to +\infty$ as $n\to\infty$. By the argument of
the previous paragraph, there is a subsequence $t_{k(n)}$ such that
$\gamma'_n=\gamma(t+t_{k(n)})$ converges uniformly (over $[0,1]$)
to a generalized
broken bicharacteristic $\gamma'$, which is a constant curve, $\gamma'(t)=\xi$
for $t\in[0,1]$,
in $R_+(\lambda)\cup R_-(\lambda)$. In particular, $\gamma(t_{k(n)})$
converges to $\xi$. But on the other hand,
$\xi_0=\lim_{n\to\infty}\gamma(t_n)$, so $\xi_0=\xi\in R_+(\lambda)\cup
R_-(\lambda)$. This proves that 
\begin{equation}
\overline{\gamma([t_0,+\infty))}\subset\gamma([t_0,+\infty))\cup R_+(\lambda)
\cup R_-(\lambda).
\end{equation}
\end{proof}

\begin{rem}
Note that under the assumptions of Theorem~\ref{thm:geom-br-bichar},
$\gamma|_{[T,+\infty)}$ is an integral curve of $\scHg^a$ for a sufficiently
large $T$, so
$\xi=\lim_{t\to+\infty}\gamma(t)$ exists in $\dot\Sigma(\lambda)$.
Moreover, if $\xi\in R_+(\lambda)$ (i.e.\ it is incoming), then
$\gamma|_{[T,+\infty)}$ is a constant curve, so $\xi\in\gamma([t_0,+\infty))$.
\end{rem}

\section{Positive commutators}\label{sec:commutators}
Our approach to prove propagation of singularities along generalized
broken bicharacteristics
relies on positive commutator estimates. Just as in \cite{Vasy:Propagation-2}
and \cite{Vasy:Propagation-Many}, we need to commute operators
arising by quantization of $\pi$-invariant functions
$q\in\Cinf(\sct_{\bXb}\Xb)$, $\Xb=\Snp$ as usual. That $q$ is $\pi$-invariant
means that at $C'_a$, the regular part of $C_a$, given by $x_a=0$,
$z_a=0$, $q(y_a,z_a,\xi_a,\xi^a)|_{z_a=0}$ is independent of $\xi^a$,
i.e.\ $q(y_a,0,\xi_a,\xi^a)$ is independent of $\xi_a$. This means that
(except in some special cases) $q$ is not symbolic at cotangent
infinity, i.e.\ as $\xi\to\infty$, since typically $\xi_a$ derivatives
will not give any decay in the $\xi^a$ direction. Recall here
that a (non-polyhomogeneous) scattering symbol $r$ of, say, multiorder $(0,0)$,
$r\in\bcon^{(0,0)}(\Snp\times\Snp)$, is an element of
$\Cinf(\Rn_w\times\Rn_\xi)$
that satisfies (product-type) symbol
estimates in both sets of variables, i.e.\ for all $\alpha,\beta\in\Nat^n$
there exists $C_{\alpha\beta}>0$ such that
\begin{equation}
|D_w^\alpha D_\xi^\beta r(w,\xi)|\leq C_{\alpha\beta}\langle w\rangle
^{-|\alpha|}\langle\xi\rangle^{-|\beta|}.
\end{equation}
However, in our arguments
the failure of $q$ to be in $\bcon^{(0,0)}(\Snp\times\Snp)$
will never cause a problem, only an inconvenience, since we always localize in
the spectrum of the elliptic operator $H$ by composing our operators
with $\psi_0(H)\in\PsiSc^{-\infty,0}(\Xb,\calC)$ where $\psi_0\in\Cinf_c
(\Real)$. Such operators $\psi_0(H)$
are smoothing -- their amplitudes decay rapidly in $\xi$.

So fix $\psi_0\in\Cinf_c(\Real;[0,1])$
which is identically $1$ in a neighborhood
of a fixed $\lambda$.
Thus,
$\psi_0(H)\in\PsiSc^{-\infty,0}(\Xb,\calC)$, so it is smoothing.
At the symbol level,
$\psi_0(H)$ is locally
the right quantization of some
\begin{equation}\label{eq:comm-3}
p\in\Cinf([\Snp;\calC]\times\Snp)
\end{equation}
which vanishes to infinite order at $[\Snp;\calC]\times\partial\Snp$, i.e.\ it
is Schwartz in $\xi$.

We are then interested in the following class of symbols $q$. We assume
that $q\in\Cinf(\Rn_w\times\Rn_\xi)$ and that for
every multiindex $\alpha$, $\beta\in\Nat^n$ there exist constants
$C_{\alpha,\beta}$ and $m_{\alpha,\beta}$ such that
\begin{equation}\label{eq:comm-4-a}
|(D_w^\alpha
D_\xi^\beta q)(w,\xi)|\leq C_{\alpha,\beta}\langle w\rangle^{-|\alpha|}
\langle\xi\rangle^{m_{\alpha,\beta}}.
\end{equation}
This implies, in particular, that
\begin{equation}\label{eq:comm-4-ap}
q\in\bcon^0(\Snp\times\Rn),
\end{equation}
i.e.\ that $q$ is a 0th order symbol in $w$, though it may blow up
polynomially in $\xi$. Indeed, in the compactified notation,
\eqref{eq:comm-4-a} becomes that for every $P\in\Diffb(\Snp)$,
acting in the base ($w$) variables, and for every
$\beta\in\Nat^n$ there exist $C_{P,\beta}$ and $m_{P,\beta}$ such that
\begin{equation}\label{eq:comm-4-b}
|(P D_\xi^\beta) q| \leq C_{P,\beta}
\langle\xi\rangle^{m_{P,\beta}}.
\end{equation}

It is convenient to require in addition
that $q$ be polyhomogeneous on $\Snp\times\Rn$:
\begin{equation}\label{eq:comm-5}
q\in\Cinf(\Snp\times\Rn);
\end{equation}
this stronger statement (implying \eqref{eq:comm-4-ap})
automatically holds for the $\pi$-invariant
symbols we are interested in.

We next introduce the product symbol
\begin{equation}\label{eq:comm-7}
a(w,w',\xi)=q(w,\xi)p(w',\xi),
\end{equation}
where $\psi_0(H)$ is given locally by the right quantization of $p$.
The main point is

\begin{lemma}\label{lemma:comm-1}
The symbol $a$ defined by \eqref{eq:comm-7} is in $\Cinf(\Snp
\times[\Snp;\calC]
\times\Snp)$ and it vanishes with all derivatives at
$\Snp\times
[\Snp;\calC]\times\partial\Snp$. Hence, it defines an operator
$A\in\PsiSc^{-\infty,0}(\Snp,\calC)$ by the oscillatory integral
\begin{equation}\label{eq:Psop-37}
Au(w)=(2\pi)^{-n}\int e^{i(w-w')\cdot\xi}a(w,w',\xi)\,u(w')\,dw'\,d\xi.
\end{equation}
If in addition $q$ is $\pi$-invariant, then for all $\xi\in\sct(\Ct_a;\Xb)$
$\Ah_a(\xi)=q(\xi)\widehat{\psi_0(H)}_a(\xi)$. In particular,
$\widehat{[A,H]}_{a,0}$ vanishes for all $a$, so $[A,H]\in
\PsiSc^{-\infty,1}(\Snp,\calC)$.
\end{lemma}

The last statements can be seen directly, but they follow particularly
easily from the forthcoming discussion in which we quantize $q$ itself
(without the factor $p$) to act on oscillatory functions, and from the
oscillatory testing definition of the indicial operators that we have adopted.

Our immediate goal is now to prove that under certain conditions, the
most important of which is that in a part of phase space, say $U\subset\dot
\Sigma(\lambda)$,
$\scHg^b q$ is negative for all $b$,
\begin{equation}
i\psi(H)[A^*A,H]\psi(H)\geq B^*B+E+F
\end{equation}
where $B,E\in\PsiSc^{-\infty,0}(\Snp,\calC)$ (so they are bounded on
$L^2_{\scl}(\Snp)$), $F\in\PsiSc^{-\infty,1}(\Snp,\calC)$ (so $F$ is compact
on $L^2_{\scl}(\Snp)$), $\WFScp(E)$ disjoint from $U$ and
$\psi\in\Cinf_c(\Real;[0,1])$ has sufficiently small support near $\lambda$.
In view of $B^*B$, this says that the commutator $i[A^*A,H]$
is microlocally positive at $U$. The precise version of this statement
is Proposition~\ref{prop:comm-7} at the end of this section.

Due to a square root
construction in the ps.d.o. algebra, see
\cite[Corollary~9.7]{Vasy:Propagation-Many},
this result follows immediately if we prove that all indicial
operators of $i\psi(H)[A^*A,H]\psi(H)$ over $U$ are positive.
There are two closely related ways to proceed prove the positivity of the
indicial operators. A general approach is to follow the proof
of the Mourre estimate due to Froese and Herbst \cite{FroMourre}, but
replace their global statements by the appropriate indicial operator ones.
This approach is necessary in more general geometric settings.
However, in the Euclidean setting it turns out that all one needs is the
Mourre estimate in {\em all} the proper subsystems. We follow this path.
Thus, we will need to evaluate some commutators quite explicitly, and use
the Mourre estimate in the normal variables and the standard Hamilton
vector field commutator estimate in the tangential variables.

Although the evaluation of the commutator appears (and is)
rather obvious if we write
$A=Q\psi_0(H)$, think of $Q$ as a scattering ps.d.o.\ obtained by
quantizing the symbol $q$, and use
\begin{equation}\label{eq:comm-8}
[A,H]=[Q\psi_0(H),H]=[Q,H]\psi_0(H),
\end{equation}
this
argument is rather formal since $q$ is not a scattering symbol.
Here we follow \cite{Vasy:Propagation-2} to make precise sense of this.

So let $q$ be as in Lemma~\ref{lemma:comm-1},
and suppose that $f\in\Cinf(\Xb;\Real)$, $u=e^{if/x}v$,
$v\in\Cinf([\Xb;\calC])$. Thus, $u\in\cap_{r=0}^\infty H^{r,s}(\Rn)$ for
$s<-n/2$. Using our indicial
operator definition, $A\in\PsiSc^{-\infty,0}(\Xb,\calC)$, the above calculation,
\eqref{eq:comm-8},
makes sense and gives the correct result if $Q$, defined by
\begin{equation}
Qu=(2\pi)^{-n}\int_{\Rn}e^{i w\cdot\xi}q(w,\xi)\hat u(\xi)\,d\xi,
\end{equation}
acts on such oscillatory functions.
But $u\in\cap_{r=0}^\infty H^{r,s}(\Rn)$, so
$\hat u\in\cap_{r=0}^{\infty} H^{s,r}(\Rn)$, so the integral makes sense as
a distributional pairing. The same holds if we replace $q(w,\xi)$ by
$q^{(k)}=q(w,\xi)\langle\xi\rangle^{-2k}$, and in addition, by choosing
$k$ large enough, $q^{(k)}$ satisfies any fixed number of scattering
symbol estimates as a scattering symbol of multiorder, say, $(0,0)$.
In addition,
$Qu=Q^{(k)}(\Delta+1)^k u$, and $(\Delta+1)^k:e^{if/x}\Cinf([\Xb;\calC])
\to e^{if/x}\Cinf([\Xb;\calC])$.
Since any fixed symbol estimate
of $e^{-if/x}Pu$,
$u\in e^{if/x}\Cinf([\Xb;\calC])$, requires only a finite number of scattering
symbol estimates on $P$, we conclude that
\begin{equation}
Q:e^{if/x}\Cinf([\Xb;\calC])\to e^{if/x}\Cinf([\Xb;\calC]).
\end{equation}
In addition, to calculate $Ru$, $R\in\Psisc^{*,0}(\Xb)$, modulo $x^l
e^{if/x}\Cinf([\Xb;\calC])$ requires only a finite number of scattering
symbol estimates on the amplitude $r$ of $R$,
we conclude that in all indicial operator calculations,
as well as calculations of lower order terms, we can work as if
$Q$ were in $\Psisc^{*,0}(\Xb)$.
We also note that if $R\in\Psisc^{*,0}(\Xb)$ has sc-principal symbol
$r(y_a,z_a,\xi_a,\xi^a)$, its $a$-indicial operator, $\Rh_{a,0}(y_a,\xi_a)$
is the
translation-invariant ps.d.o.\ $I(r)\in\Psisc^{*,0}(\Xb^a)$ on $X^a$ given by
\begin{equation}\label{eq:comm-19}
I(r(y_a,0,\xi_a,.))v=\Fr_{X^a}^{-1}r(y_a,0,\xi_a,\xi^a)\Fr_{X^a}v,\ v\in
\temp(X^a);
\end{equation}
see \cite{Vasy:Propagation-2, Vasy:Propagation-Many}.
Here $\Fr_{X^a}$ denotes the Fourier transform on $X^a$.
Correspondingly, when applied to $u\in e^{if/x}v$, $v\in\Cinf([\Xb;\calC])$,
$e^{-if/x}Qu$ is given by $I(q(y_a,0,\xi_a,.))v$
modulo $x \Cinf([\Xb;\calC])$, at $C'_a$. But $q$ is $\pi$-invariant,
so $q(y_a,0,\xi_a,.)$ is independent of $\xi^a$, so $I(q(y_a,0,\xi_a,.))$
is just multiplication by $q_\pi(y_a,0,\xi_a)$, which we will
simply write as $q(y_a,0,\xi_a)$ as in Section~\ref{sec:br-bichar}.
Simlarly, $\scHg^a q|_{z=0}$ is also independent of $\xi^a$, since $q|_{z=0}$
is and $\scHg^a=2|w_a|^{-1}\xi_a\cdot\pa_{w_a}$ is tangent to $\sct_{C_a}\Snp$,
so $I((\scHg^a q)(y_a,0,\xi_a,.))$ is just multiplication by
$(\scHg^a q)_\pi(y_a,0,\xi_a,.)=(\scHg^a q)(y_a,0,\xi_a)$.

While the preceeding discussion was quite general, we will be
interested in rather special choices of $q$, which will make the
indicial operator calculation simpler. In particular, these choices will
allow us to use the standard Mourre estimate on $X^a$, rather than
essentially repeating its proof, as presented by Froese and Herbst
\cite{FroMourre}, see also \cite{Derezinski-Gerard:Scattering}, to obtain
a corresponding estimate for more general operators in the normal variables.
Thus, we will make a first order assumption on $q$ at the collision
planes at infinity, i.e.\ at the $C_a$, as well,
in addition to its $\pi$-invariance. Namely, we assume that for all
clusters $a$, $q$ has the structure
\begin{equation}\label{eq:d_z-q-16}
d_z q|_{z=0}=q_1\, \xi^a\cdot dz,\quad q_1\ \pi-\text{invariant},
\end{equation}
so for all $a$,
\begin{equation}\label{eq:d_z-q-16a}
q=q_0+q_1\eta_a+\sum_{i,j}z_i z_j q_{ij},
\end{equation}
with $q_0,q_1\in\Cinf(\sct\Xb_a)$, $q_{ij}\in\Cinf(\sct\Snp)$.

We use this as follows. For $u\in e^{if/x}\Cinf([\Xb;\calC])$ we can write
\begin{equation}\label{eq:comm-23}
[A,H]u=AHu-HAu=Q\psi_0(H)Hu-HQ\psi_0(H)u=(QH-HQ)\psi_0(H)u.
\end{equation}
As indicated above, during calculation the indicial operator of
$[A,H]$, we can thus work as if $Q$ were a scattering
ps.d.o.
Recall that our many-body Hamiltonians are of the form
\begin{equation}\label{eq:comm-16}
H=H^a+\Delta_{X_a}+I_a,\ I_a\in x\Diffsc^2(U),
\end{equation}
$U$ a neighborhood of $C'_a$, disjoint from $C_b$ such that $C_b\not\supset
C_a$.
Using \eqref{eq:comm-16} and \eqref{eq:comm-23}, we write
\begin{equation}
[A,H]u=[Q,\Delta_{X_a}]\psi_0(H)u+[Q,I_a]\psi_0(H)u+[Q,H^a]\psi_0(H)u.
\end{equation}
Since the principal symbol in the scattering calculus is given by the
Poisson bracket, for $\ut= e^{if/x}\vt$, $\vt\in\Cinf([\Xb;\calC])$,
$p\in\ff([\Xb;C_a])$ with $\beta[\Xb;C_a](p)\in C'_a$,
\begin{equation}\label{eq:comm-24}
ix^{-1}e^{-if/x}[Q,\Delta_{X_a}]\ut(p)
=I((-\scHg^a q)(\xi))\vt(p)=-(\scHg^a q_0)(\xi)
\vt(p),\quad df(p)=\xi,
\end{equation}
here we used that as $\scHg^a q=\scHg^a q_0$ is $\pi$-invariant,
$I(\scHg^a q)$ is a multiplication operator.
In addition, as $I_a\in x\DiffSc^2(\Xb,\calC)$, we see that a priori
$QI_a, I_a Q:e^{if/x}\Cinf([\Xb;\calC])\to x e^{if/x}\Cinf([\Xb;\calC])$
with indicial operator given by the composition of the two indicial
operators. Since that of $Q$ is $I(q)=q\,\Id$,
a multiple of the identity, we conclude
that the $(a,1)$-indicial operator of $[Q,I_a]$ vanishes, so
in fact $[Q,I_a]:e^{if/x}\Cinf([\Xb;\calC])\to x^2 e^{if/x}\Cinf([\Xb;\calC])$.

The only remaining term is $[Q,H^a]\psi_0(H)$.
Note that $z_j=|w_a|^{-1}(w^a)_j$ as usual, so $\eta_a=|w_a|^{-1}w^a\cdot
\xi^a$. Thus, with $Q_i$, $Q_{ij}$ denoting the quantization of $q_i$ and
$q_{ij}$ respectively,
the same way as $Q$ is the quantization of $q$,
\begin{equation}
Q=Q_0+|w^a|^{-1} Q_1 w^a\cdot D_{w^a}+\sum_{i,j}z_i z_j Q_{ij}.
\end{equation}
Here $Q_0$, $Q_1$ depend on the $\sct\Xb_a$ variables only, so they
commute with $H^a$. On the other hand, $[z_i,H^a]\in\PsiSc^{1,1}(\Snp,\calC)$
(since all indicial operators of $z_i$ are multiples of the indentity),
and $Q_{ij}:e^{if/x}\Cinf([\Xb;\calC])\to e^{if/x}\Cinf([\Xb;\calC])$,
so the same holds for $[Q_{ij},H^a]$. Thus,
\begin{equation}
[z_i z_j Q_{ij},H^a]:e^{if/x}\Cinf([\Xb;\calC])\to
|z|^2 e^{if/x}\Cinf([\Xb;\calC]).
\end{equation}
This proves that for oscillatory functions $u$,
\begin{equation}
[Q,H^a]\psi_0(H)u=|w_a|^{-1}q_1 [w^a\cdot D_{w^a},H^a]\psi_0(H)u+u',
\quad u'\in |z|^2 e^{if/x}\Cinf([\Xb;\calC]),
\end{equation}
where $q_1$ is evaluated at $df$ as in \eqref{eq:comm-24}.
Thus, we conclude that
\begin{equation}\begin{split}\label{eq:comm-H-h_a}
i\widehat{[A^*A,H]}_{a,1}(y_a,\xi_a)=&-2q
\scHg^a q_0(y_a,0,\xi_a)\psi_0(H^a+|\xi_a|^2)^2\\
&\quad+i[B(y_a,\xi_a)+B(y_a,\xi_a)^*,H^a]
\end{split}\end{equation}
where
\begin{equation}\label{eq:comm-B-def}
B(y_a,\xi_a)=\psi_0(H^a+|\xi_a|^2)q q_1
\sum_j (w^a)_j D_{(w^a)_j} \psi_0(H^a+|\xi_a|^2)\in\Psisc^{-\infty,-1}(\Xb^a)
\end{equation}
is a ($(y_a,\xi_a)$-dependent) multiple of the generator of dilations on $X^a$.
Although we worked at $C'_a$ for the sake of convenience,
$i\widehat{[A^*A,H]}_{a,1}$ is a well-defined element of
$\Cinf(\sct(\Ct_a;\Xb);\PsiSc^{-\infty,0}(\Xb^a,\calC^a))$; in particular
it is continuous on $\sct(\Ct_a;\Xb)$, in which $\sct_{C'_a}(\Ct_a;\Xb)$ is
dense, so \eqref{eq:comm-H-h_a} in fact gives $i\widehat{[A^*A,H]}_{a,1}$
on $\sct(\Ct_a;\Xb)$.

We also remark that if $\psi\in\Cinf_c(\Real;[0,1])$, $\psi_0\equiv 1$ on
$\supp\psi$, then $\psi(H)=\psi_0(H)\psi(H)$, so $A\psi(H)=Q\psi_0(H)\psi(H)
=Q\psi(H)$ on oscillatory functions. Thus,
\begin{equation}\begin{split}\label{eq:comm-H-h_a-8}
i\widehat{\psi(H)[A^*A,H]\psi(H)}_{a,1}(y_a,\xi_a)=&-2q
\scHg^a q_0(y_a,0,\xi_a)\psi(H^a+|\xi_a|^2)^2\\
&\quad+i[B_\psi(y_a,\xi_a)+B_\psi(y_a,\xi_a)^*,H^a]
\end{split}\end{equation}
where
\begin{equation}\label{eq:comm-B-def-8}
B(y_a,\xi_a)=\psi(H^a+|\xi_a|^2)q q_1
\sum_j (w^a)_j D_{(w^a)_j} \psi(H^a+|\xi_a|^2)\in\Psisc^{-\infty,-1}(\Xb^a).
\end{equation}

This also proves the following lemma.

\begin{lemma}\label{lemma:comm-4}
(see also \cite[Lemma~9.4]{Vasy:Propagation-Many})
Let $q$ and $A$ be as in Lemma~\ref{lemma:comm-1}, $q$ being
$\pi$-invariant. Suppose also that $q$ is of the form
\eqref{eq:d_z-q-16a}. For every seminorm
in $\PsiSc^{-\infty,0}(\Xb^a,\calC^a)$
and for every $l\in\Nat$
there exist
$C>0$ and $m\in\Nat$ such that for every $a$ and every
$\xi\in\sct_p(\Ct_a;\Xb)$, $p\in\Ct_a$,
the seminorm
of $\widehat{[A,H]}_{a,1}(\xi)$
in $\PsiSc^{-\infty,0}(\Xb^a,\calC^a)$ is bounded by
\begin{equation}
C(|q(\xi)|+|q_1(\xi)|).
\end{equation}
\end{lemma}

\begin{rem}
Similar conclusions hold for every seminorm in
$\Cinf(\sct(\Ct_a;\Xb),\PsiSc^{-\infty,0}(\Xb^a,\calC^a))$,
which can be seen directly
from our calculations in the preceeding proof.
\end{rem}

In view of \eqref{eq:comm-B-def},
the Mourre estimate is immediately applicable; we first recall it.
We present it essentially
as stated in \cite{Derezinski-Gerard:Scattering}, statement
$H_3(a)$, Section~6.4.
Let
\begin{equation}\begin{split}
&d^a(\sigma)=\inf\{\sigma-\sigma':\ \sigma'\leq\sigma,\ \sigma'\in\Lambda'_a\}
\geq 0,\Mif \sigma\geq\inf\Lambda'_a,\\
&d^a(\sigma)=0,\Mif \sigma<\inf\Lambda'_a,
\end{split}\end{equation}
\begin{equation}
d^{a,\kappa}(\sigma)=\inf\{d^a(\sigma'):\ \sigma'\in[\sigma-\kappa,\sigma+
\kappa]\}.
\end{equation}
Then given $\ep>0$, $\sigma_0\in\Real$,
there exists $\delta>0$ such that for all $\sigma\leq\sigma_0$ and
for all $\psi_1\in\Cinf_c(\Real;
[0,1])$ supported in $(\sigma-\delta,\sigma+\delta)$,
\begin{equation}\begin{split}\label{eq:Mourre-1}
-i[&\psi(H^a)(w^a\cdot D_{w^a}+D_{w^a}\cdot w^a)
 \psi(H^a+|\xi_a|^2),H^a]/2\\
&\geq 2(d^{a,\delta}(\sigma)-\ep)\psi_1(H^a)^2.
\end{split}\end{equation}
Note that this definition of $d^a(\sigma)$ differs slightly from
that of \cite{Derezinski-Gerard:Scattering}, namely our $d^a(\sigma)$
is bounded from above by theirs if $\sigma\geq\inf\Lambda_a$,
and \eqref{eq:Mourre-1} is trivial if $\sigma<\inf\Lambda_a$.

We apply this with $\sigma=\lambda-|\xi_a|^2$, $\sigma_0=\lambda+1$.
Thus, given $\ep>0$
there exists $\delta>0$ such that for all $\psi\in\Cinf_c(\Real;
[0,1])$ supported in $(\lambda-\delta,\lambda+\delta)$,
\begin{equation}\begin{split}
-i[&\psi(H^a+|\xi_a|^2)(w^a\cdot D_{w^a}+D_{w^a}\cdot w^a)
 \psi(H^a+|\xi_a|^2),H^a]/2\\
&\geq 2(d^{a,\delta}(\lambda-|\xi_a|^2)-\ep)\psi(H^a+|\xi_a|^2)^2.
\end{split}\end{equation}
Note that if $\sigma\geq\inf\Lambda'_a$,
\begin{equation}\begin{split}\label{eq:d^a-eta^a}
2 d^a(\sigma)&=2\inf\{|\xi^a_b|^2:\ \sigma-|\xi^a_b|^2\in\pspec(H^b),
\ C_b\supset C_a\}\\
&=\inf\{H_g^b(\xi^a\cdot w^a):\ \sigma-|\xi^a_b|^2\in\pspec(H^b),
\ C_b\supset C_a\}\\
&=\inf\{\scHg^b\eta_a(y_a,0,\xi_a,\xi^a_b):\ \sigma-|\xi^a_b|^2\in\pspec(H^b),
\ C_b\supset C_a\}\geq 0;
\end{split}\end{equation}
here $\scHg^b\eta_a=|\xi_a^b|^2$
is independent of $y_a$ and $\xi_a$, so the appearance
of these on the last line is irrelevant. For the sake of
convenience, in the next paragarphs, we write $\inf\emptyset=0$,
then \eqref{eq:d^a-eta^a} holds for all $\sigma\in\Real$. We remark that
$\xi^a\cdot w^a$ is the principal symbol of
$(w^a\cdot D_{w^a}+D_{w^a}\cdot w^a)/2$ in the scattering calculus (here
the part at infinity is the one that matters),
which explains how it arises in the commutator estimates.
Thus, with $I=[\lambda-\delta,\lambda+\delta]$,
\begin{equation}\begin{split}\label{eq:d^ad-H^b}
2&d^{a,\delta}(\lambda-|\xi_a|^2)\\
&=\inf\{\scHg^b\eta_a(\xit):\ \lambda'-|\xi_a|^2-|\xi^a_b|^2\in\pspec(H^b),
\ C_b\supset C_a,\ \lambda'\in[\lambda-\delta,\lambda+\delta]\}\\
&=\inf\{\scHg^b\eta_a(\xit):\ \pi_b(\xit)=\xi,\ C_b\supset C_a,
\ \exists\sigma\in\pspec(H^b),\ |\lambda-|\xit|^2-\sigma|\leq \delta\}\\
&=\inf\{\scHg^b\eta_a(\xit):\ \ C_b\supset C_a,
\ \xit\in\pih_{b,I}^{-1}(\{\xi\})\}.
\end{split}\end{equation}

In our estimates, in the part of phase space where we wish to prove that
the commutator is positive,
we will have
\begin{equation}
q\geq 0,\ q_1\leq 0.
\end{equation}
Note that due to
\eqref{eq:d_z-q-16},
for $b$ with $C_b\supset C_a$,
\begin{equation}\label{eq:scHg^b-q-z=0}
\scHg^b q|_{z=0}=\scHg^a q_0+q_1\scHg^b\eta_a.
\end{equation}
The Mourre estimate thus shows that
\begin{equation}
\label{eq:Mourre-8}
i\psi(H^a+|\xi_a|^2)[B+B^*,H^a]\psi(H^a+|\xi_a|^2)
\geq -4 qq_1 (d^{a,\delta}(\lambda-|\xi_a|^2)-\ep)\psi(H^a+|\xi_a|^2)^2.
\end{equation}
Substituting this into \eqref{eq:comm-H-h_a} yields
\begin{equation}\begin{split}\label{eq:comm-H-h_a-2p}
i&\psi(H^a+|\xi_a|^2)\widehat{[A^*A,H]}_{a,1}(\xi)\psi(H^a+|\xi_a|^2)\\
&\geq -2q(\scHg^a q_0+2q_1 d^{a,\delta}(\lambda-|\xi_a|^2)-2\ep q_1)
\psi(H^a+|\xi_a|^2)^2,
\end{split}\end{equation}
with
\begin{equation}\begin{split}\label{eq:comm-H-h_a-2pp}
&-\scHg^a q_0-2q_1 d^{a,\delta}(\lambda-|\xi_a|^2)\\
&\qquad=\inf\{-\scHg^b q(\xit):\ \ C_b\supset C_a,
\ \xit\in\pih_{b,I}^{-1}(\{\xi\})\},
\end{split}\end{equation}
where we used
\eqref{eq:d^ad-H^b}, \eqref{eq:scHg^b-q-z=0}. In \eqref{eq:comm-H-h_a-2pp},
$\inf\emptyset$ must be understood
as $-\scHg^a q_0$, but due to the factors of $\psi(H^a+|\xi_a|^2)$ in
\eqref{eq:comm-H-h_a}, that equation holds for any value
replacing $\scHg^a q_0+2q_1 d^{a,\delta}(\lambda-|\xi_a|^2)-2\ep q_1$
if $\lambda-|\xi_a|^2<\inf\Lambda'_a-\delta$ (which is exactly when
$\inf\emptyset$ arises above).

In order to estimate the commutator on some set $U\subset\scdt\Snp$,
it is convenient to rewrite this as follows.
Suppose that
$b\in\Cinf(\sct\Snp)$ is $\pi$-invariant, $\delta_0>0$, $I=[\lambda-\delta_0,
\lambda+\delta_0]$,
for all clusters $b$,
\begin{equation}
\xit\in\pih_{b,I}^{-1}(U) \Rightarrow \scHg^b q(\xit)\leq -b(\xit)^2,
\end{equation}
and there exists $C>0$ such that for all $\xi\in U$
\begin{equation}
|q_1(\xi)|\leq Cb(\xi)^2.
\end{equation}
Then for all $\ep'>0$ there exists $\delta'>0$ such that
if $\supp\psi\subset(\lambda-\delta',\lambda+\delta')$ then
\begin{equation}\label{eq:comm-H-h_a-2}
i\psi(H^a+|\xi_a|^2)\widehat{[A^*A,H]}_{a,1}\psi(H^a+|\xi_a|^2)
\geq (2-\ep')qb^2\psi(H^a+|\xi_a|^2)^2.
\end{equation}
We have thus proved the following proposition.

\begin{prop}\label{prop:comm-1}
(cf.\ \cite[Proposition~9.6]{Vasy:Propagation-Many})
Suppose that $H$ is a many-body Hamiltonian
and $\lambda\in\Real$.
Suppose also that $q,b\in\Cinf(\sct \Snp;\Real)$ are $\pi$-invariant,
for all clusters $a$
\begin{equation}\label{eq:comm-33}
q=q_0+q_1 \eta_a+O(|z_a|^2),\ q_0,q_1\ \pi-\text{invariant},
\end{equation}
they satisfy
the bounds \eqref{eq:comm-4-b}, $q,b,-q_1\geq 0$,
and that there
exist $\delta>0$, $C>0$, $C'>0$, $U\subset\dot\Sigma(I)$, where
$I=[\lambda-\delta,\lambda+\delta]$,
such that for every index $a$ and for all $\xi\in\sct_{C_a}
\Xb_a$,
\begin{equation}\label{eq:comm-35}
\xi\in\pih_{a,I}^{-1}(U)\Rightarrow
\scHg^a q(\xi)\leq -b(\xi)^2
\end{equation}
and
\begin{equation}\label{eq:comm-37}
\xi\in\pih_{a,I}^{-1}(U)\Rightarrow
q(\xi)\leq Cb(\xi)^2\Mand|q_1(\xi)|\leq C'
b(\xi)^2.
\end{equation}
Let $A\in\PsiSc^{-\infty,0}(\Xb,\calC)$ be as in Lemma~\ref{lemma:comm-1}.
For any $\ep'>0$
there exists $\delta'>0$
such that
if $\psi\in\Cinf_c(\Real)$ is supported in $(\lambda-\delta',\lambda+\delta')$
and $\xi\in\sct(\Ct_a;\Xb)$ for some $a$ with $\pit(\xi)\in U$, then
\begin{equation}\label{eq:ind-est-1}
i\widehat{(\psi(H)[A^*A,H]\psi(H))}_{a,1}(\xi)
\geq (2-\ep')b^2 q\psi(\Hh_{a,0}(\xi))^2.
\end{equation}
\end{prop}

We can now apply \cite[Corollary~9.7]{Vasy:Propagation-Many}
to obtain estimates on the original operators. Note that the proof of
this corollary consists of a square root construction,
\cite[Proposition~8.3]{Vasy:Propagation-Many}, for which the existence
of $L^2$ eigenvalues in subsystems is irrelevant, and an indicial
operator estimate, which was given in Proposition~9.6 in
\cite{Vasy:Propagation-Many}; here its place is taken by the
preceeding result, Proposition~\ref{prop:comm-1}. Recall that
the square root construction is via a ps.d.o.\ version of the functional
calculus, and that \eqref{eq:comm-37-p} ensures that this can be done
within the ps.d.o.\ calculus. Also, \eqref{eq:comm-37-p} might appear
weaker than its analogue in \cite{Vasy:Propagation-Many}, but due to
\eqref{eq:comm-33}, the other estimates of Corollary~9.7 of that paper
follow automatically.
Note that $\pi$-invariant (continuous)
functions on $\sct\Snp$ can be regarded as (continuous) functions on
$\scdt\Snp$, so it makes sense to talk about the support of these functions
in $\scdt\Snp$.

\begin{prop}\label{prop:comm-7}
(cf.\ \cite[Corollary~9.7]{Vasy:Propagation-Many})
Suppose that the assumptions of Proposition~\ref{prop:comm-1} are
satisfied and let $C$ be as in \eqref{eq:comm-37}, and $U$ open.
Suppose in addition that for any cluster $a$ and any
differential operator $Q$ on
$\sct(\Ct_a;\Xb)$ there exists a constant
$C_Q$ such that
\begin{equation}\begin{split}\label{eq:comm-37-p}
&\pit(\xi)\in U\Mand b(\xi)\neq 0\\
&\Rightarrow
|Q(b^{-2} q)(\xi)|\leq C_Q\Mand|Q(b^{-2}q_1)(\xi)|\leq C_Q.
\end{split}\end{equation}
Then for any $\ep'>0$, $M>0$, and for any
$K\subset\scdt \Xb$ compact with $K\subset U$
there exists $\delta'>0$, $\psi\in\Cinf_c(\Real)$ is supported in
$(\lambda-\delta',\lambda+\delta')$, $\psi\equiv 1$ near $\lambda$,
$B,E\in\PsiSc^{-\infty,0}(\Xb,\calC)$, $F\in\PsiSc^{-\infty,1}(\Xb,\calC)$
with $\WFScp(B),\WFScp(E),\WFScp(F)\subset\dot\Sigma(I)$,
\begin{equation}\label{eq:comm-51}
\WFScp(E)\cap K=\emptyset,\ \WFScp(F)\subset\supp q,\ \Bh_{a,0}(\xi)=
b(\xi)q(\xi)^{1/2}
\psi(\Hh_{a,0}(\xi)),\ \pit(\xi)\in K,
\end{equation}
such that
\begin{equation}\label{eq:comm-53}
i\psi(H)x^{-1/2}[A^*A,H]x^{-1/2}\psi(H)-M\psi(H)A^*A\psi(H)
\geq (2-\ep'-MC)B^*B+E+F.
\end{equation}
\end{prop}

\begin{rem}
Multiplying both sides of \eqref{eq:comm-53} by $\psi_1(H)$ such that
$\psi\equiv 1$ on $\supp\psi_1$ shows that $\psi$ can be replaced by
any such $\psi_1$.
\end{rem}

\section{Propagation of singularities}\label{sec:propagation}
For our positive commutator estimates it is convenient to consider
two scenarios separately, though we present them parallel to each other.
Recall that we are interested in proving propagation estimates
in $S=\dot\Sigma(\lambda)\setminus(R_+(\lambda)\cup R_-(\lambda))$.
Note that $\xib\in S\cap\sct_{C'_a}\Xb_a$, $\xib=(\bar y_a,\taub_a,\mub_a)$,
means that for all $b$
with $C_a\subset C_b$, $\xib_b\in\pih_b^{-1}(\xib)$, $\mu_b(\xib_b)\neq 0$.
The two possibilities are:

\begin{enumerate}
\item
$\mub_a\neq 0$. Then $\scHg^a$ does not vanish at $\xib$,
so there is tangential propagation along $C_a$.

\item
$\mub_a=0$. Then there is no tangential propagation along $C_a$.
(However, note that as $\xib\in S$
this automatically gives $\xib\nin\Sigma_a
(\lambda)$, so certainly the particles cannot be in an $a$-bound state.)
In this case, by \eqref{eq:scHg^b-eta_a},
\begin{equation}
\scHg^b(z_a\cdot\nu_a)=(x_a/x)(2\tau_a(z_a\cdot\nu_a)+2|\nu_{ab}|^2)>0
\end{equation}
(near $\pih_b^{-1}(\xib)$) since $\mu_b\neq 0$ there, so $\nu_{ab}\neq 0$
by \eqref{eq:nu-mu-a-b}. Thus, we can
use $z_a\cdot\nu_a$ as a parameter along broken bicharacteristics.
\end{enumerate}

The main technical result on the propagation of singularities is the following.

\begin{prop}\label{prop:prop}
Suppose that $H$ is a many-body Hamiltonian.
Let $u\in\dist(\Xb)$,
$\lambda\nin\Lambda_1$.
Let
$\xib\in S\cap\sct_{C'_a}\Xb_a$, $\xib=(\bar y_a,\taub_a,\mub_a)$, and
suppose that $\xib\nin\WFSc((H-\lambda)u)$.
Let $\eta_a=z_a\cdot \nu_a$.
Then there exists a neighborhood $U_0$ of $\xib$ with the following
properties.

\begin{enumerate}
\item
If $\mub_a\neq 0$ and there exists a neighborhood $U\subset U_0$ of $\xib$ in
$\dot\Sigma(\lambda)$ and if there exists $s_0>0$
such that
$\exp(-s_0\scHg^a)(\xib)\nin\WFSc(u)$ and for all $s$ satisfying
$0\leq s\leq s_0$, $\exp(-s\scHg^a)(\xib)\in U$ (i.e.\ the
integral curve segment of $\scHg^a$ is completely inside $U$),
and in addition
\begin{equation}\begin{split}\label{eq:prop-i}
\xi\in U\Mand \eta_a(\xi)<0\Rightarrow\xi\nin\WFSc(u),
\end{split}\end{equation}
then $\xib\nin\WFSc(u)$.

\item
If $\mub_a=0$ and
there exists a neighborhood $U\subset U_0$ of $\xib$ in $\dot\Sigma(\lambda)$
such that
\begin{equation}\begin{split}\label{eq:prop-ii}
\xi\in U\Mand \eta_a(\xi)<0\Rightarrow\xi\nin\WFSc(u),
\end{split}\end{equation}
then $\xib\nin\WFSc(u)$.
\end{enumerate}
\end{prop}

\begin{rem}
Note that $\eta_a(\xi)<0$ implies $z_a\neq 0$, so $\xi\nin\sct_{C_a'}\Xb_a$.
\end{rem}

\begin{proof}
The main step in the proof is the construction of an operator which has
a microlocally positive commutator with $H$ near $\xi_0$. In fact, we
construct the symbol of this operator. This symbol will not be a
scattering symbol, i.e.\ it will not be in $\Cinf(\Snp\times\Snp)$,
only due to its behavior as $\nu\to\infty$ corresponding to its
$\pi$-invariance. This
will be accommodated by composing its quantization with a cutoff in the
spectrum of $H$, $\psi(H)$, $\psi\in\Cinf_c(\Real)$ supported near $\lambda$,
as discussed in Lemma~\ref{lemma:comm-1}. This approach simply
extends the one taken in
\cite{Vasy:Propagation-2}-\cite{Vasy:Propagation-Many}, though the
actual construction is different due to the more complicated geometry.

We introduce some notation and then fix $U_0$. We define $c_0>0$ by
\begin{equation}
c_0=\half\inf\{-\scHg^b\tau(\xib_b):\ \xib_b\in\pih_b^{-1}(\xib),
\ C_b\supset C_a\}>0.
\end{equation}
Note that at $C_a$, where $|w_a|/|w|=x/x_a=1$,
\begin{equation}
\tau=\tau_a-\eta_a,\ \text{since}\ \tau=-\frac{w\cdot\xi}{|w|},
\ \tau_a=-\frac{w^a\cdot\xi^a}{|w_a|},\ \eta_a=\frac{w_a\cdot\xi_a}{|w_a|},
\end{equation}
and if $\mub_a=0$, then $\scHg\tau_a(\xib_b)=0$ when $\pi(\xib_b)=\xib$.
Thus, if $\mub_a=0$, then
\begin{equation}
c_0=\half\inf\{\scHg^b(z_a\cdot\nu_a)
(\xib_b):\ \xib_b\in\pih_b^{-1}(\xib)\}>0.
\end{equation}
Now fix $I=[\lambda-\delta_0,\lambda+\delta_0]$ and
$U_0$ to a neighborhood of $\xib$ in $\dot\Sigma(I)$
such that $U_0\subset\bigcup\{\sct_{C'_b}\Xb_b:\ C_a\subset C_b\}$,
$U_0\cap\WFSc((H-\lambda)u)=\emptyset$ (this is possible since
$\WFSc((H-\lambda)u)$ is closed), $U_0$ is inside a fixed system of
local coordinates, and
\begin{equation}\begin{split}\label{eq:c_0-est}
&\text{if}\ \mub_a\neq 0\Mthen
 \xi\in U_0,\ \xib\in\pih_{b,I}^{-1}(\xi),\ C_b\supset C_a\Rightarrow
-\scHg^b\tau(\xi_b)>3c_0/2\\
&\text{if}\ \mub_a= 0\Mthen
 \xi\in U_0,\ \xi_b\in\pih_{b,I}^{-1}(\xi),\ C_b\supset C_a\Rightarrow
\scHg^b(z_a\cdot\nu_a)(\xi_b)>3c_0/2.
\end{split}\end{equation}

If $\mub_a\neq 0$, we make an additional definition. Namely, we let
$\delta>0$ be given by
\begin{equation}
-\delta=\taub-\tau(\exp(-s_0\scHg^a)(\xib)).
\end{equation}
We will use $\ep\in(0,1)$
as a small parameter that microlocalizes in a neighborhood
of $\xib$ if $\mub_a=0$, and in a neighborhood of
$\{\exp(-s\scHg^a)(\xib):\ s\in[0,s_0]\}$ if $\mub_a\neq 0$; if $\mub_a\neq 0$
we will always take $\ep<\delta$.

Employing an iterative argument as usual, we may assume
that $\xib\nin\WFSc^{*,l}(u)$ and we need
to show that $\xib\nin\WFSc^{*,l+1/2}(u)$. (We can start the induction
with an $l$ such that $u\in\Hsc^{*,l}(\Xb)$.)

Our positive commutator estimates at a point $\xib_a=(\bar y_a,\taub,
\mub_a)$ will arise by considering functions
\begin{equation}\begin{split}
&\phi=\taub-\tau+\ep^{-1}(|z_a|^2+
\omega),\Mif\mub_a\neq 0,\\
&\phi=z_a\cdot\nu_a+\frac{\beta}{\ep}(|z_a|^2+\omega)
,\Mif\mub_a= 0,
\end{split}\end{equation}
where $\omega$ localizes in the tangential variables $y_a$, $\tau_a$, $\mu_a$.
In fact, if $\mub_a\neq 0$, then $\scHg^a(\xib_a)\neq 0$, so we construct
$\omega$ to be a `quadratic distance' (in $\sct_{C_a}\Xb_a$) from
the $\scHg^a$ integral curve through $\xib_a$ constructed so that
$\scHg^a\omega=0$. That is, we define $\omega$ on a hypersurface
through $\xib$ that is transversal to $\scHg^a$, e.g.\ $\tau_a=\taub_a$,
to be given by a positive
definite quadratic form in some local coordinates centered at $\xib$,
e.g.\ $(y_a-\bar y_a)^2+(\mu_a-\mub_a)^2$, and
extend it to $\sct_{C_a}\Xb_a$
to be constant along the $\scHg^a$ integral curves;
cf.\ \cite[Section~7]{Vasy:Propagation-Many}.
We can then extend $\omega$ to a function on
$\sct_{\bXb}\Xb$ defined near $C_a$ as discussed before \eqref{eq:pi_a^e-scHg},
and then
\begin{equation}
\scHg\omega=0,\ \text{so}\ \scHg^b\omega=0\ \text{for all}\ b.
\end{equation}
On the other hand, if $\mub_a=0$, we can take
\begin{equation}
\omega_0=(\xi_a-\xib_a)^2+(y_a-\bar y_a)^2;
\end{equation}
now
\begin{equation}
\scHg^a \omega_0=2(y_a-\bar y_a)\cdot\scHg^a y_a=4\mu_a\cdot(y_a-
\bar y_a),
\end{equation}
so $|\scHg^a\omega_0|\leq C'\omega_0$. Thus, for $C$ sufficiently large, we
see that for
\begin{equation}
\omega=e^{-C\tau}\omega_0,\ \scHg^a \omega\geq 0
\end{equation}
since $-\scHg^a\tau>0$ on $S$, hence near $\xib$.
Since $\omega$ depends on the tangential variables only, we conclude in
either case that for all $b$,
\begin{equation}
\scHg^b\omega\geq 0.
\end{equation}
Moreover,
$\scHg^b|z_a|^2=4z_{ab}\cdot\nu_{ab}$ under the decomposition $\mu_b=(\mu_a,
\nu_{ab})$, so $\scHg^b|z_a|^2\leq C_1|z_a|$.

Now suppose that $\mub_a\neq 0$,
\begin{equation}\label{eq:case-i-cond-1}
\phi\leq 2\ep,\ \taub-\tau\geq -\delta-\ep.
\end{equation}
Since $\omega\geq 0$ and $|z_a|^2\geq 0$,
the first of these inequalities implies that
$\taub-\tau\leq 2\ep$, and the combination of these two gives
$\phi-(\taub-\tau)\leq \delta+3\ep\leq 4\delta$ as $\ep<\delta$.
Then we
conclude that
\begin{equation}\label{eq:supp-q-a}
-\delta-\ep\leq\taub-\tau\leq 2\ep,\ |z_a|\leq 2(\ep\delta)^{1/2},
\ \omega^{1/2}\leq 2(\delta\ep)^{1/2}.
\end{equation}
Note that for all $b$ with $C_b\supset C_a$, $z_{ab}\cdot\nu_{ab}=z_a
\cdot\nu_a$ at $\sct_{C_b}\Xb_b$. Thus, if
\begin{equation}
z_{a}\cdot\nu_{a}\geq -C_3 \ep  c_0
\end{equation}
we see that
\begin{equation}
\scHg^b|z_a|^2=4z_{a}\cdot\nu_{a}\geq -4 C_3\ep  c_0.
\end{equation}
Thus, using \eqref{eq:c_0-est} as well, we deduce that
there exists $C_4>0$ (independent of $\ep>0$) such that
\begin{equation}\label{eq:case-i-pos-est-1}
\text{if}\ z_a\cdot\nu_a\geq-C_4 \ep,\ \text{then}
\ \scHg^b\phi\geq -\scHg^b\tau-c_0/2>c_0.
\end{equation}
Note that here $\ep\in (0,\delta)$ is arbitrary.

Similarly, suppose that $\mub_a= 0$,
\begin{equation}\label{eq:case-ii-cond-1}
\phi\leq 2\ep,\ z_a\cdot\nu_a\geq -2\ep.
\end{equation}
Then we
conclude that
\begin{equation}\label{eq:supp-q-a-2}
|z_a\cdot\nu_a|\leq 2\ep,\ |z_a|\leq 2\ep/\sqrt{\beta},
\ \omega^{1/2}\leq 2\ep/\sqrt{\beta}.
\end{equation}
Thus, from \eqref{eq:c_0-est},
\begin{equation}
\scHg^b\phi\geq \scHg^b(z_a\cdot\nu_a)-2\sqrt{\beta} C_1
>c_0\Mif\beta=(c_0/4C_1)^2.
\end{equation}

The positive commutator estimate then arises by considering the following
$\pi$-invariant symbol $q$ and quantizing it as in
Lemma~\ref{lemma:comm-1}.
Let $\chi_0\in\Cinf(\Real)$ be equal to $0$ on $(-\infty,0]$ and
$\chi_0(t)=\exp(-1/t)$ for $t>0$. Thus, $\chi_0'(t)=t^{-2}\chi_0(t)$, $t>0$,
and $\chi'_0(t)=0$, $t\leq 0$.
Let $\chi_1\in\Cinf(\Real)$ be $0$
on $(-\infty,0]$, $1$ on $[1,\infty)$, with $\chi_1'\geq 0$
and $\chi_1(t)=\exp(-1/t)$
on some small interval
$(0,t_0)$, $t_0>0$.
Furthermore, for $A_0>0$ large, to be determined, let
\begin{equation}\label{eq:prop-22}
q=\chi_0(A_0^{-1}(2-\phi/\ep))\chi_1((\taub-\tau+\delta)/
\ep+1)
\end{equation}
if $\mub_a\neq 0$, and let
\begin{equation}\label{eq:prop-22b}
q=\chi_0(A_0^{-1}(2-\phi/\ep))\chi_1((z_a\cdot\nu_a)/
\ep+2)
\end{equation}
if $\mub_a=0$.
Thus, $q(\xit)=\chi_0(2/A_0)>0$, and on $\supp q$ we have
\begin{equation}
\phi\leq 2\ep\Mand \taub-\tau\geq -\delta-\ep\ \text{in case (i)},
\end{equation}
and
\begin{equation}
\phi\leq 2\ep\Mand z_a\cdot\nu_a\geq -2\ep\ \text{in case (ii)},
\end{equation}
which are exactly \eqref{eq:case-i-cond-1} and \eqref{eq:case-ii-cond-1},
so $\supp q$ is a subset of \eqref{eq:supp-q-a} and \eqref{eq:supp-q-a-2}
respectively.
We also see that
as $\ep$ decreases, so does $\supp q=\supp q_\ep$,
in fact, if $0<\ep'<\ep$ then $q_\ep>0$ on $\supp q_{\ep'}$.
Note that
in case (i),
by reducing $\ep$, we can make $q$ supported in an arbitrary small
neighborhood of a compact backward bicharacteristic segment through
$\xib$, and in case (ii), by reducing $\ep$,
we can make $q$ supported in an arbitrary small
neighborhood of $\xib$.

We at once obtain positivity estimates for $\scHg^c q$.
The following argument works
similarly for both $\mub_a=0$ and $\mub_a\neq 0$; we consider the slightly
more complicated case
$\mub_a\neq 0$.
Thus, if $\mub_a\neq 0$, then
\begin{equation}\begin{split}\label{eq:scHg q-calc}
\scHg^c q=-A_0^{-1}&\ep^{-1}\chi'_0(A_0^{-1}(2-\phi/\ep))
\chi_1((\taub-\tau+\delta)/\ep+1)\scHg^c \phi\\
&-\ep^{-1}
\chi_0(A_0^{-1}(2-\phi/\ep))\chi_1'((\taub-\tau+\delta)/\ep+1)
\scHg^c \tau.
\end{split}\end{equation}
We break up the first term by using a cutoff
that ensures that the hypothesis in \eqref{eq:case-i-pos-est-1} is
satisfied.
Thus, let
\begin{equation}
\chit=\chi_1(4z_{a}\cdot\nu_{a}/(C_4\ep)+2),
\end{equation}
so
\begin{equation}\label{eq:prop-26}
\text{on}\ \supp \chit,\ 
z_{a}\cdot\nu_{a}\geq-C_4 \ep/2,
\end{equation}
\begin{equation}
\text{and on}\ \supp(1-\chit),\ z_{a}\cdot\nu_{a}\leq-C_4 \ep/4.
\end{equation}
Then
\begin{equation}\label{eq:scHg q-pos}
\scHg^c q=-\bt^2_c+e_c
\end{equation}
with
\begin{equation}
\bt^2_c=A_0^{-1}\ep^{-1}\chi'_0(A_0^{-1}(2-\phi/\ep))
\chi_1((\taub-\tau+\delta)/\ep+1)\chit \scHg^c \phi.
\end{equation}
Hence, with
\begin{equation}
b^2=c_0 A_0^{-1}\ep^{-1}\chi'_0(A_0^{-1}(2-\phi/\ep))
\chi_1((\taub-\tau+\delta)/\ep+1)\chit,
\end{equation}
using \eqref{eq:case-i-pos-est-1} and \eqref{eq:prop-26}, we deduce that
\begin{equation}
\scHg^c q\leq -b^2+e_c.
\end{equation}
Moreover,
\begin{equation}\label{eq:prop-33}
b^2\geq (c_0 A_0/16) q
\end{equation}
since $\phi\geq \taub-\tau\geq -2\delta$ on $\supp q$, so
\begin{equation}\begin{split}\label{eq:prop-34}
\chi'_0(A_0^{-1}(2-\phi/\ep))
&=A_0^2(2-\phi/\ep)^{-2}\chi_0(A_0^{-1}(2-\phi/\ep))\\
&\geq (A_0^2/16)
\chi_0(A_0^{-1}(2-\phi/\ep)).
\end{split}\end{equation}
On the other hand, $e_c$ is supported where either
\begin{equation}\label{eq:supp-e-fine}
-\delta-\epsilon\leq\taub-\tau\leq-\delta,
\ \omega^{1/2}\leq 2(\ep\delta)^{1/2},\ |z_a|\leq 2(\ep\delta)^{1/2}
\end{equation}
so near the backward direction along $\scHg^a$
bicharacteristic through $\xib$, or
\begin{equation}\label{eq:supp-e-fine-2}
z_{a}\cdot\nu_{a}\leq -C_4\ep/4,\Mand -\delta-\ep\leq
\taub-\tau\leq 2\ep,
\ \omega^{1/2}\leq 2(\ep\delta)^{1/2},\ |z_a|\leq 2(\ep\delta)^{1/2}.
\end{equation}
But by our assumption
$\exp(-s_0\scHg^a)(\xib)\nin\WFSc(u)$, $s_0>0$,
so the same holds for a sufficiently
small neighborhood of $\exp(-s_0\scHg^a)(\xib)$ as $\WFSc(u)$ is closed.
By choosing $\ep>0$ sufficiently small, we can thus make sure that
the region defined by \eqref{eq:supp-e-fine} is disjoint from $\WFSc(u)$.
Moreover, by further reducing $\ep>0$ if necessary and using our
second assumption, we can also make sure that the region
\eqref{eq:supp-e-fine-2} is also disjoint from $\WFSc(u)$, so that
$\supp e_c$ is disjoint from $\WFSc(u)$ for all $c$.
Moreover, by \eqref{eq:supp-q-a},
for $\ep>0$ sufficiently small, we deduce from the
inductive hypothesis that $\supp q$ (hence $\supp b$) is disjoint
from $\WFSc^{*,l}(u)\cap\dot\Sigma(I)$.

Moreover, with $\partial$ denoting a partial derivative with respect to
one of $(y,z,\tau,\mu,\nu)$,
\begin{equation}\begin{split}
\partial q=-A_0^{-1}&\ep^{-1}\chi'_0(A_0^{-1}(2-\phi/\ep))
\chi_1((\taub-\tau+\delta)/\ep+1)\partial \phi\\
&-\ep^{-1}
\chi_0(A_0^{-1}(2-\phi/\ep))\chi_1'((\taub-\tau+\delta)/\ep+1)
\partial \tau.
\end{split}\end{equation}
At any $C_b$ with $p\in C_b$, defined by $x=0$, $z''=0$, as
above, $\phi$ is independent of $\nu''$ at $z''=0$ so outside $\supp e_c$
\begin{equation}
|\partial_{\nu''}^\beta dq|\leq C_\beta b^2\ \text{at}\ z''=0.
\end{equation}
In fact, outside $\supp e_c$, but in the set where $b$ is positive,
\begin{equation}\label{eq:prop-39}
b^{-2}\partial q=c_0^{-1}\partial\phi,
\end{equation}
so the uniform bounds of \eqref{eq:comm-37-p} also follow.
In addition, at any cluster $b$, $|z_a|^2=|z_{ab}|^2+|z_a^b|^2$,
$\eta_a=z_{ab}\cdot\nu_{ab}+z_a^b\cdot\nu_a^b=z_{ab}\cdot\nu_{ab}+\eta_b$,
and $z_{ab}$, $\nu_{ab}$ are $b$-tangential variables, so
$\phi$, hence $q$ has the form \eqref{eq:comm-33} around each $C_b$.

Let $\psit\in\Cinf_c(\Real)$ be identically $1$ near $0$ and supported
close to $0$.
We also define
\begin{equation}
\qt=\psit(x)q.
\end{equation}
Thus, $\qt\in\Cinf(\sct \Xb)$ is a $\pi$-invariant function satisfying
\eqref{eq:comm-4-b}. Let $A$ be the operator given by Lemma~\ref{lemma:comm-1}
with $\qt$ in place of $q$,
so in particular its indicial operators are $q(\xi)\psi_0(\Hh_{b,0}(\xi))$.
Note that \eqref{eq:comm-37} holds with $C=16 c_0^{-1}A_0^{-1}$. So suppose
that $M>0$ and $\ep'>0$. Choose $A_0$ so large that $MC<\ep'$,
and let $U$ be the complement of
$\cup_c \pi_c(\supp e_c)$ in $\dot\Sigma(I)$, and let $K=\WFSc(u)\cap
\dot\Sigma(I)$; so $K\subset U$ by our choice of $\ep>0$.
By Proposition~\ref{prop:comm-7} we deduce that
there exists $\delta'>0$, $\psi\in\Cinf_c(\Real)$ is supported in
$(\lambda-\delta',\lambda+\delta')$, $\psi\equiv 1$ near $\lambda$,
$B,E\in\PsiSc^{-\infty,0}(\Xb,\calC)$, $F\in\PsiSc^{-\infty,1}(\Xb,\calC)$
with $\WFScp(B),\WFScp(E),\WFScp(F)\subset\dot\Sigma(I)$,
\begin{equation}\label{eq:comm-51a}
\WFScp(E)\cap K=\emptyset,\ \WFScp(F)\subset\supp q,\ \Bh_{a,0}(\xi)=
b(\xi)q(\xi)^{1/2}
\psi(\Hh_{a,0}(\xi)),\ \pit(\xi)\in K,
\end{equation}
such that
\begin{equation}\label{eq:comm-53a}
i\psi(H)x^{-1/2}[A^*A,H]x^{-1/2}\psi(H)-M\psi(H)A^*A\psi(H)
\geq (2-2\ep')B^*B+E+F.
\end{equation}

Let
\begin{equation}
\Lambda_r=x^{-l-1/2}(1+r/x)^{-1},\quad r\in(0,1),
\end{equation}
so $\Lambda_r\in\PsiSc^{0,-l+1/2}(\Xb,\calC)$ for $r\in(0,1)$ and it is
uniformly bounded in $\PsiSc^{0,-l-1/2}(\Xb,\calC)$. The last statement
follows from $(1+r/x)^{-1}$ being uniformly bounded as a 0th order
symbol, i.e.\ from $(x\partial_x)^k(1+r/x)^{-1}\leq C_k$ uniformly
($C_k$ independent of $r$).
We also define
\begin{equation}\label{eq:prop-57}
A_r=A\Lambda_r x^{-1/2}\psi(H),\ B_r=B\Lambda_r,\ E_r=\Lambda_r E\Lambda_r.
\end{equation}
Then, with $\psi_0\in\Cinf_c(\Real;[0,1])$ identically $1$ near $\supp\psi$,
\begin{equation}\begin{split}
ix^{l+1/2}&[A^*_rA_r,H]x^{l+1/2}\\
&=i(1+r/x)^{-1}\psi(H)x^{-1/2}[A^*A,H]x^{-1/2}\psi(H)(1+r/x)^{-1}\\
&\qquad+i\psi(H)A^*x^{l+1/2}[\Lambda_r x^{-1/2},H](1+r/x)^{-1}x^{-1/2}\psi_0(H)
A\psi(H)\\
&\qquad
+i\psi(H)A^*\psi_0(H)
x^{-1/2}(1+r/x)^{-1}[\Lambda_r x^{-1/2},H]x^{l+1/2}A\psi(H)
+H_r,
\end{split}\end{equation}
where $H_r$ is uniformly bounded in $\PsiSc^{-\infty,1}(\Xb,\calC)$.
Note that $H_r$ arises by commuting $A$, powers of $x$ and $\Lambda_r$ through
other operators, but as the indicial operators of $A$ and $x$
are a multiple of the identity, $A$, $x$ and $\Lambda_r$ commute with these
operators to top order, and in case of $\Lambda_r$, the commutator
is uniformly bounded as an operator of one lower order.
Then, multiplying \eqref{eq:comm-53a} by $(1+r/x)^{-1}$ from the left
and right and rearranging the terms we obtain the following estimate of
bounded self-adjoint operators on $L^2_\scl(\Xb)$:
\begin{equation}\label{eq:prop-61}\begin{split}
ix^{l+1/2}[A^*_rA_r&,H]x^{l+1/2}
-\psi(H)A^*(G_r^*+G_r)A\psi(H)-M\psi(H)A^*A\psi(H)\\
&\geq x^{l+1/2}((2-\epsilon')B_r^*B_r+E_r+F_r)x^{l+1/2}
\end{split}\end{equation}
where
\begin{equation}
G_r=i\psi_0(H)x^{-1/2}(1+r/x)^{-1}[\Lambda_r x^{-1/2},H]x^{l+1/2},
\end{equation}
and $F_r\in\PsiSc^{-\infty,-2l+1}(\Xb,\calC)$ is
uniformly bounded in $\PsiSc^{-\infty,-2l}(\Xb,\calC)$ as $r\to 0$.
Now, $G_r\in\PsiSc^{-\infty,2}(\Xb,\calC)$ is uniformly bounded
in $\PsiSc^{0,0}(\Xb,\calC)$, hence as a bounded operator on $L^2_\scl(\Xb)$.
Thus, if $M>0$ is chosen sufficiently large, then $G_r+G_r^*\geq- M$ for
all $r\in(0,1)$, so
\begin{equation}
\psi(H)A^*(G_r+G_r^*+M)A\psi(H)\geq 0.
\end{equation}
Adding this to \eqref{eq:prop-61} shows that
\begin{equation}\label{eq:prop-61a}
ix^{l+1/2}[A^*_rA_r,H]x^{l+1/2}
\geq x^{l+1/2}((2-\epsilon')B_r^*B_r+E_r+F_r)x^{l+1/2}.
\end{equation}
The point of the commutator calculation is that in
$L^2_\scl(\Xb)$
\begin{equation}\begin{split}
\langle u,&[A_r^*A_r,H]u\rangle\\
&=\langle u,A_r^*A_r(H-\lambda)u\rangle
-\langle u,(H-\lambda)A_r^*A_r u\rangle\\
&=2i\im\langle u,A_r^*A_r(H-\lambda)u\rangle;
\end{split}\end{equation}
the pairing makes sense for $r>0$ since $A_r\in\PsiSc^{-\infty,-l}(\Xb,\calC)$.
Now apply \eqref{eq:prop-61a} to $x^{-l-1/2}u$ and pair it with
$x^{-l-1/2}u$ in $L^2_{\scl}(\Xb)$. Then for $r>0$
\begin{equation}\label{eq:pos-comm-est}
\|B_r u\|^2\leq|\langle u,E_r u\rangle|+|\langle u,F_ru\rangle|
+2|\langle u,A^*_rA_r(H-\lambda)u\rangle|.
\end{equation}

Letting $r\to 0$ now keeps the right hand side of \eqref{eq:pos-comm-est}
bounded.
In fact,
$A_r(H-\lambda)u\in\dCinf(\Xb)$ remains bounded in $\dCinf(\Xb)$ as $r\to 0$.
Similarly, by \eqref{eq:comm-51a},
$E_r u$ remains bounded in $\dCinf(\Xb)$
as $r\to 0$ since $\WFSc(u)\cap \WFScp(E)=\emptyset$.
Also, $F_r$ is bounded in $\bop(\Hsc^{m,l}(\Xb),\Hsc^{-m,-l}(\Xb))$, so
$\langle u,F_r u\rangle$ stays bounded by \eqref{eq:comm-51a} as well
since $\supp q\cap\WFSc^{*,l}(u)\cap\dot\Sigma(I)=\emptyset$.
These estimates
show that $B_r u$ is uniformly bounded in $L^2_\scl(\Xb)$.
Since $(1+r/x)^{-1}\to \Id$ strongly
on $\bop(\Hsc^{m',l'}(\Xb),\Hsc^{m',l'}(\Xb))$, we conclude that
$Bx^{-l-1/2}u\in L^2_\scl(\Xb)$. By \eqref{eq:comm-51a}
and Proposition~\ref{prop:WF-3} this implies that for every $m$,
\begin{equation}
\xib\nin\WFSc^{m,l+1/2}(u);
\end{equation}
in fact that $\xi\nin\WFSc^{m,l+1/2}(u)$ for all $\xi\in K$ for which
$q>0$.
This is exactly the iterative step we wanted to prove. In the next step
we decrease $\ep>0$
slightly to ensure that $\WFScp(F)\subset \supp\qt$
is disjoint from $\WFSc^{m,l+1/2}(u)$, just as H\"ormander decreases
$t$ in the proof of \cite[Proposition~24.5.1]{Hor}.
\end{proof}

Based on this proposition, we
proceed inductively, assuming that propagation has been proved at $C'_b$ with
$C_a\subsetneq C_b$. Noting that $\eta_a=z_{a}\cdot\nu_{a}<0$ implies that a
backward broken bicharacteristic through $\xi$ must stay away from $C_a$,
we can use
Lebeau's argument, as presented in \cite[Theorem~10.6]{Vasy:Propagation-Many},
to prove propagation along broken bicharacteristics. The basic idea
is that $\xib\in
\WFSc(u)\cap\sct_{C'_a}\Xb_a$ implies that either every point on the
non-constant backward
$\scHg^a$-bicharacteristic is in $\WFSc(u)$ (in a neighborhood of $\xib$),
or there are points in $\WFSc(u)\setminus\sct_{C'_a}\Xb_a$
arbitrarily close to $\xib$ such
that there is a backward broken bicharacteristic through these points which
is completely in $\WFSc(u)$ (by the inductive hypothesis). Using the
compactness of the set of broken bicharacteristics, we can extract
a sequence of backward broken bicharacteristics which converge to
a backward broken bicharacteristic through $\xib$, which will prove the
propagation statement.

\begin{thm*}(Theorem~\ref{thm:prop-sing})
Let $u\in\dist(\Xb)$,
$\lambda\nin\Lambda_1$. Then
\begin{equation}
\WFSc(u)\setminus\WFSc((H-\lambda)u)
\end{equation}
is a union of maximally extended
generalized broken bicharacteristics of $H-\lambda$ in
$\dot\Sigma(\lambda)\setminus\WFSc((H-\lambda)u)$.
\end{thm*}

\begin{proof}
As usual, broken bicharacteristic means {\em generalized} broken
bicharacteristic in this proof.

We only need to prove that for every cluster $a$, if
\begin{equation}\label{eq:prop-91}
\xib\in\WFSc(u)\setminus\WFSc((H-\lambda)u)\Mand \xib\in\sct_{C'_a}
\Xb_a
\end{equation}
then
\begin{equation}\begin{split}\label{eq:prop-92}
\exists\ &\text{broken bicharacteristic}
\ \gamma:[-\ep_0,\ep_0]\to\dot\Sigma(\lambda),\ \ep_0>0,\Mst\\
&\gamma(0)=\xib,
\ \gamma(t)\in\WFSc(u)\setminus\WFSc((H-\lambda)u)\Mfor t\in[-\ep_0,\ep_0].
\end{split}\end{equation}
In fact, if \eqref{eq:prop-91}$\Rightarrow$\eqref{eq:prop-92} holds
for all $a$ with $C_c\subset C_a$, let
\begin{equation}\begin{split}
\calR=&\{\text{broken bicharacteristics}\\
&\ \gamma:(\alpha,\beta)
\to (\WFSc(u)\setminus\WFSc((H-\lambda)u))\cap\bigcup
\{\sct_{C'_a}\Xb_a:\ C_c\subset C_a\},\\
&\ \alpha<0<\beta,\ \gamma(0)=\xib\},
\end{split}\end{equation}
and put the natural partial order on $\calR$, so $\gamma\leq\gamma'$
if the domains satisfy $(\alpha,\beta)\subset(\alpha',\beta')$ and
$\gamma=\gamma'|_{(\alpha,\beta)}$. Then $\calR$ is not empty
(due to \eqref{eq:prop-92})
and every non-empty
totally ordered subset
of $\calR$ has an upper bound,
so an application of Zorn's lemma gives
a maximal broken bicharacteristic of $H-\lambda$
in the intersection of
$\WFSc(u)\setminus\WFSc((H-\lambda)u)$ with $\cup_{C_c\subset C_a}
\sct_{C'_a}\Xb_a$ which passes through $\xib$. A similar maximal
statement holds if we replace $C_c\subset C_a$ by $C_c\subsetneq C_a$.

Indeed, it suffices to show that
for any $a$, if
\begin{equation}\label{eq:prop-103}
\xib\in\WFSc(u)\setminus\WFSc((H-\lambda)u)\Mand\xib\in\sct_{C'_a}
(C_a;\Xb)
\end{equation}
then
\begin{equation}\label{eq:prop-104}\begin{split}
&\text{there exists a broken bicharacteristic}
\ \gamma:[-\ep_0,0]\to\dot\Sigma(\lambda),\ \ep_0>0,\\
&\qquad\qquad \gamma(0)=\xib,
\ \gamma(t)\in\WFSc(u)\setminus\WFSc((H-\lambda)u),\ t\in[-\ep_0,0],
\end{split}\end{equation}
for the existence of a broken bicharacteristic
on $[0,\ep_0]$ can be demonstrated similarly
by replacing the forward propagation
estimates by backward ones, and, directly from
Definition~\ref{Def:gen-br-bichar}, piecing together the two
broken bicharacteristics gives one defined on $[-\ep_0,\ep_0]$.

We proceed to prove that \eqref{eq:prop-103} implies
\eqref{eq:prop-104} by induction on $a$. This is
certainly true for $a=0$ by Proposition~\ref{prop:prop}: there are
no normal variables $z_a$, $\nu_a$, so $\eta_a=0$ in the notation
of that Proposition, showing that a segment of
the backward bicharacteristic through
$\xib$ must be in $\WFSc(u)$. Of course, this is simply Melrose's
propagation theorem \cite[Proposition~7]{RBMSpec}.

In addition, if $a$ is arbitrary and $\xib\in R_+(\lambda)\cup R_-(\lambda)$,
then the constant curve $\gamma$ through $\xib$ is a broken bicharacteristic,
so  \eqref{eq:prop-103}$\Rightarrow$\eqref{eq:prop-104} holds with
this $\gamma$.

So suppose that \eqref{eq:prop-103}$\Rightarrow$\eqref{eq:prop-104}
has been proved for all $b$ with
$C_a\subsetneq C_b$ and that $\xib\in\dot\Sigma(\lambda)\cap\sct_{C'_a}
(C_a;\Xb)$ satisfies \eqref{eq:prop-103}.
As noted in the first paragraph, we thus know
that
the intersection of
$\WFSc(u)\setminus\WFSc((H-\lambda)u)$ with $\cup_{C_a\subsetneq C_b}
\sct_{C'_b}\Xb_b$ is a union of maximally extended
broken bicharacteristics of $H-\lambda$.
We use the notation of the proof of
Proposition~\ref{prop:prop} below.
Let $U_0$ be a neighborhood of $\xib=(\bar y_a,\taub_a,\mub_a)$ in
$\dot\Sigma(\lambda)$ as in Proposition~\ref{prop:prop};
we may assume that
$U_0\cap\WFSc((H-\lambda)u)=\emptyset$.
By Proposition~\ref{prop:prop}, either every point on the
non-constant backward
$\scHg^a$-bicharacteristic segment through $\xib$ is in $\WFSc(u)$
(in the neighborhood $U_0$ of $\xib$), in which case we have proved
\eqref{eq:prop-104}, so we are done,
or there is a sequence of points
$\xi_n\in U_0$ such that $\xi_n\in\WFSc(u)$,
$\xi_n\to\xib$ as $n\to\infty$, and
$\eta_a(\xi_n)<0$ for all $n$.
Since $\eta_a(\xi_n)<0$, $\xi_n\nin\sct_{C'_a}\Xb_a$. Thus,
by the inductive hypothesis
there exist broken bicharacteristics
$\gamma_n:(-\ep_n,0]\to\dot\Sigma(\lambda)$, $\ep_n>0$,
with $\gamma_n(0)=\xi_n$, $\gamma_n(t)\in\WFSc(u)$ for all $t\in(-\ep_n,0]$,
and $\gamma_n$ is maximal with this property in $U_0\cap
\cup_{C_a\subsetneq C_b}\sct_{C'_b}\Xb_b\cap\WFSc(u)$.
That is, if $\gammat_n$ is an extension of $\gamma_n$, then $\gammat_n(-\ep_n)
\nin\cup_{C_a\subsetneq C_b}\sct_{C'_b}\Xb_b$, so
$\gammat_n(-\ep_n)\in\sct_{C'_a}\Xb_a$. By Remark~\ref{rem:ext-to-Real}
we can fix some $\ep_0>0$ and extend (or restrict)
each $\gamma_n$ to a broken
bicharacteristic defined on $[-\ep_0,0]$,
which we keep denoting by $\gamma_n$;
by the previous remark
$\gamma_n(-\ep_n)\in\sct_{C'_a}\Xb_a$. 
But $e^{Ct}\eta_a$
is non-decreasing (for sufficiently large $C>0$) along
broken bicharacteristics by the argument that preceeds
Lemma~\ref{lemma:br-bichar-tgt-normal}, so we conclude that
\begin{equation}
\eta_a(\gamma_n(t))\leq e^{C(t-t_0)}
\eta_a(\gamma_n(0))<0
\end{equation}
for $t\in[-\ep_0,0]$, hence $z_a(\gamma_n(t))\neq 0$,
so $\ep_n\geq \ep_0$, $\gamma_n(t)\in\WFSc(u)$
for all $t\in[-\ep_0,0]$. By Proposition~\ref{prop:Lebeau-compactness},
there
is a subsequence of $\gamma_n$ converging uniformly to a broken
bicharacteristic $\gamma:[-\ep_0,0]\to\dot\Sigma(\lambda)$.
Since $\WFSc(u)$ is closed,
$\gamma:[-\ep_0,0]\to\WFSc(u)$. In particular,
$\gamma(0)=\xib$ and $\gamma(t)\in\WFSc(u)$ for all $t\in[-\ep_0,0]$,
providing the inductive step.
\end{proof}

\begin{rem}
As expected,
the theorem does not provide any interesting information at the radial sets
$R_+(\lambda)\cup R_-(\lambda)$.
\end{rem}

This theorem essentially says that the propagation of quantum particles, modulo
`trivial' (i.e.\ Schwartz) terms can be understood as a series of
collisions in which the total energy and the external momentum are preserved.

\section{The resolvent}\label{sec:resolvent}
Before we can turn Theorem~\ref{thm:prop-sing} into a result on the
wave front relation of the S-matrix, we need to analyze the resolvent.
More precisely, we need to understand the boundary values
\begin{equation}
R(\lambda\pm i0)=(H-(\lambda\pm i0))^{-1}
\end{equation}
of the resolvent at the real axis in a microlocal sense.
To do so, we also need estimates at the radial sets $R_\pm(\lambda)$.
Since some of the Hamilton vector fields $\scHg^a$
of the metric $g$ vanish at
$R_+(\lambda)\cup R_-(\lambda)$, the estimates must utilize
the weights $x^{-l-1}$ themselves. In this sense they are delicate,
but on the other hand they only involve $x$ and its sc-microlocal
dual variable $\tau$, so they do not need to reflect the
geometry of $\calC$. The best known positive commutator estimate is
the Mourre estimate, originally proved by Perry, Sigal and Simon
in Euclidean many-body scattering \cite{Perry-Sigal-Simon:Spectral},
in which one takes $q=x^{-1}\tau$ with the notation
of Section~\ref{sec:commutators}. Since it is easy to analyze the commutator
of powers of $x$ with $H$ (in particular, they commute with $V$),
the functional calculus allows one to obtain microlocal estimates
from these, as was done by G\'erard, Isozaki and Skibsted
\cite{GerComm, GIS:N-body}.
Thus, nearly all the technical results in this section can be proved, for
example, by using the Mourre estimate and Theorem~\ref{thm:prop-sing}.
In particular, apart from the propagation statements, they are well-known
in our Euclidean many-body setting.

We first recall the limiting absorption principle.

\begin{thm}\label{thm:lim-absorb}
Suppose that $H$ is a many-body Hamiltonian, and
$\lambda\nin\Lambda$.
Let $f\in\dCinf(\Xb)$, $u_t^\pm=R(\lambda\pm it)f$, $t>0$. Then
$u_t^\pm$ has a limit $u_\pm=R(\lambda\pm i0)f$ in $\Hsc^{m,l}(\Xb)$, $l<-1/2$,
as $t\to 0$. In addition,
\begin{equation}\label{eq:lim-absorb}
\WFSc(u_\pm)\subset \Phi_\pm(R_\mp(\lambda)).
\end{equation}
\end{thm}

\begin{rem}
The first statement in the theorem also holds if we merely assume
$f\in\Hsc^{m,l'}(\Xb)$ with $l'>1/2$, but then $\WFSc(u_\pm)$ has to be replaced
by the filtered wave front set $\WFSc^{m,l'-1}(u_\pm)$. Moreover,
$R(\lambda\pm i0)$ give continuous operators from $\Hsc^{m,l'}(\Xb)$ to
$\Hsc^{m+2,l}(\Xb)$.
\end{rem}

\begin{proof}
The existence of $u_\pm$ in $\Hsc^{m,l}(\Xb)$, $l<-1/2$ follows from
the Mourre estimate, as presented in \cite{Perry-Sigal-Simon:Spectral}.
That $\xi\in\WFSc(u_+)$ implies $\tau(\xi)<0$ (in fact, $\tau(\xi)
\leq -\sqrt{a(\lambda)}$ where $a(\lambda)=\inf\{\lambda-\sigma:
\sigma\in\Lambda,\ \sigma<\lambda\}$), follows from the work of
G\'erard, Isozaki and Skibsted \cite{GerComm}; see
\cite{Vasy:Propagation-Many} to see how the proof would proceed with
our notation. By Theorem~\ref{thm:prop-sing},
$\WFSc(u_+)\setminus\WFSc(f)=\WFSc(u_+)$ is a union of maximally
extended generalized broken bicharacteristics. So suppose
that $\xi\in\WFSc(u_+)$, and let $\gamma:\Real\to\dot\Sigma(\lambda)$
be a generalized broken
bicharacteristic in $\WFSc(u_+)$ with $\gamma(t_0)=\xi$. Then, by
Lemma~\ref{lemma:t-to-infty},
$\tau(-\infty)=\lim_{t\to -\infty}\tau(\gamma(t))$ exists. If
$\tau(-\infty)\geq 0$ then $\tau(\gamma(t))>-\sqrt{a(\lambda)}$ for
large negative $t$, contradicting that $\tau\leq-\sqrt{a(\lambda)}$
on $\WFSc(u_+)$. Thus, $\tau(-\infty)<0$, and hence
by Lemma~\ref{lemma:t-to-infty}, $\overline{\gamma((-\infty,t_0])}
\cap R_-(\lambda)\neq\emptyset$, so $\xi\in\Phi_+(R_-(\lambda))$.
We thus conclude that \eqref{eq:lim-absorb} holds.
\end{proof}

A pairing argument immediately shows $R(\lambda\pm i0)v$ also exists
for distributions $v\in\dist(\Xb)$ with wave front set disjoint from the
incoming and outgoing radial set respectively. Combining it with the
propagation theorem, Theorem~\ref{thm:prop-sing}, we
can deduce the following result.

\begin{thm*}[Theorem~\ref{thm:res-WF}]
Suppose that $H$ is a many-body Hamiltonian, and
$\lambda\nin\Lambda$.
Suppose also that $v\in\dist(\Xb)$ and $\WFSc(v)\cap\Phi_-(R_+(\lambda))
=\emptyset$.
Let $u_t^+=R(\lambda+ it) v$, $t>0$. Then $u^+_t$
has a limit
$u_+=R(\lambda+ i0)v$ in $\dist(\Xb)$ as $t\to 0$ and
$\WFSc(u_+)\cap \Phi_-(R_+(\lambda))=\emptyset$. Moreover,
$\WFSc(u_+)\subset \Phi_+(R_-(\lambda))\cup\Phi_+(\WFSc(v))$.
The result also holds with $R_+(\lambda)$ and $R_-(\lambda)$ interchanged,
$R(\lambda+it)$ replaced by $R(\lambda-it)$,
$\Phi_+$ and $\Phi_-$ interchanged.
\end{thm*}

\begin{proof}
As mentioned above, the first part follows from the self-adjointness of $H$,
so that for $t>0$, $v\in\dist(\Xb)$, $f\in\dCinf(\Xb)$, we have
$v(R(\lambda+it)f)=R(\lambda+it)v(f)$; recall that the distributional
pairing is the real pairing, not the complex (i.e.\ $L^2$) one. The
wave front statement of Theorem~\ref{thm:lim-absorb} and the assumption on
$v$ show the existence of the limit $u_+=R(\lambda+i0)v$ in $\dist(\Xb)$
and that in addition $\WFSc^{m,l}(u_+)\cap \Phi_-(R_+(\lambda))=\emptyset$ for
every
$l<-1/2$.
The statement
$\WFSc(u_+)\cap \Phi_-(R_+(\lambda))=\emptyset$ follows from a microlocalized
version of the Mourre estimate due to G\'erard, Isozaki and Skibsted
\cite{GIS:N-body}; see \cite{Hassell-Vasy:Symbolic} or
\cite{Vasy:Propagation-Many} for a detailed
argument. The final part of the conclusion follows from
Theorem~\ref{thm:prop-sing}, much as in the previous proof.
\end{proof}

\section{Poisson operators}\label{sec:Poisson}
The propagation of singularities theorem, especially in its form
for the resolvent, Theorem~\ref{thm:res-WF},
has immediate consequences for the wave front relation of
all S-matrices.
To see this, recall the definition of the Poisson operator from
\eqref{eq:Poisson-def}. First, let
\begin{equation}
I_a=\sum_{C_a\not\subset C_b}V_b,\ \It_a=(I_a/x_a)|_{C_a}\in\Cinf
(C_a),
\end{equation}
so $I_a$ is $\Cinf$ near $C_a'$ with simple vanishing at $C_a'$, so
$I_a/x_a$ is $\Cinf$ there.
Then, for $\lambda\in(\epsilon_\alpha,\infty)\setminus\Lambda$
and $g\in\Cinf_c(C_a')$,
there is a unique $u\in\dist(\Xb)$
such that $(H-\lambda)u=0$, and $u$ has the form
\begin{equation}\begin{split}\label{eq:intro-21-S}
u=&e^{-i\sqrt{\lambda-\ep_\alpha}/x}x^{\dim C_a/2}
x^{i\It_a/2\sqrt{\lambda-\ep_\alpha}}\vt_-+R(\lambda+i0)f,\\
&\vt_-=((\pi^a)^*\psi_\alpha)(\pi_a^*g)+x\vt'_-,
\ f\in\dCinf(\Xb),
\end{split}\end{equation}
and $\vt'_-$ is conormal to the boundary of $[\Xb;C_a]$, with
infinite order vanishing off the front face (the inverse image of
$C_a$ under the blow-down map), and has a full asymptotic expansion
at the front face of the form
\begin{equation}\label{eq:full-asymp-exp}
\sum_{k=0}^\infty\sum_{s=0}^{2k+2}x^k(\log x)^s a_j(y,Z),
\end{equation}
$a_j$ smooth in $y$ and Schwartz in $Z$ where we used the notation
of \eqref{eq:ff-coords}-\eqref{eq:ff-coords-2}
for the coordinates in the interior of the front face of $[\Xb;C_a]$.
(Uniqueness holds
if we merely require that $f\in\Hsc^{0,1/2+\ep'}(\Xb)$, $\ep'>0$,
due to Isozaki's results \cite{IsoUniq, IsoRad} applied to the
difference of two solutions of the form \eqref{eq:intro-21-S},
see \cite{Vasy:Scattering}; existence is a direct consequence of the
following parametrix construction.)
The Poisson operator $P_{\alpha,+}(\lambda)$ is the map
\begin{equation}\label{eq:P-def-16}
P_{\alpha,+}(\lambda):\Cinf_c(C_a')\to\dist(\Xb)\ \text{defined by}
\ P_{\alpha,+}(\lambda)g=u.
\end{equation}
The scattering matrix is then given by the following formula, which is
related to the representation formula of Isozaki and Kitada
\cite{Isozaki-Kitada:Scattering}:
\begin{equation}\label{eq:equiv-7-S}
S_{\beta\alpha}(\lambda)=\frac{1}{2i\sqrt{\lambda-\ep_\beta}}
((H-\lambda)\Pt_{\beta,-}(\lambda))^*P_{\alpha,+}(\lambda),
\end{equation}
$\lambda>\max(\epsilon_\alpha,\epsilon_\beta)$, $\lambda\nin\Lambda$
where $\Pt_{\beta,-}(\lambda)$ is a microlocalized version of $P_{\beta,-}
(\lambda)$, microlocalized near the outgoing radial set (see the
following paragraphs for details).

This formula reduces the understanding
of the structure of the S-matrices to that of the Poisson operators.
These in turn can be described in two steps. First, one constructs
the Poisson operators approximately in the incoming region of phase
space, and then one uses the propagation results for the resolvent
to obtain the wave front relation of the Poisson operators.

A very good parametrix, which we here denote by $\Pt_{\alpha,+}(\lambda)$,
for $P_{\alpha,+}(\lambda)$ in the incoming region
of phase space,
i.e.\ where $\tau$ is close to $\sqrt{\lambda-\ep_\alpha}$, has
been constructed by Skibsted \cite{Skibsted:Smoothness} in the short-range
and by Bommier \cite{Bommier:Proprietes} in the long-range setting using
a WKB-type construction under the assumption that $\ep_\alpha\in\spec_d(H^a)$,
i.e.\ $\ep_\alpha$ is below the thresholds of $H^a$.
Thus, writing their construction explicitly for the class of potentials
that we are considering and using the notation
of \eqref{eq:ff-coords}-\eqref{eq:ff-coords-2}
for the coordinates in the interior of the front face of $[\Xb;C_a]$,
they construct an operator
\begin{equation}
\Pt_{\alpha,+}(\lambda):\dist_c(C'_a)\to\dist(\Xb)
\end{equation}
by constructing its kernel
\begin{equation}\label{eq:Pt-kernel-structure}
\Kf_{\alpha,+}=e^{-i\sqrt{\lambda-\ep_\alpha}\cos\distance(y,y')/x}
x^{i\It_a(y')/2\sqrt{\lambda-\ep_\alpha}}a_+(x,y,Z,y'),
\end{equation}
$a_+\in\Cinf([\Xb;C_a]\times C_a)$
vanishes to infinite order at the `main face' (the lift of $C_0=\bXb$
to $[\Xb;C_a]$; this means that $a_+$ is Schwartz in the fibers of the
blow-down map, which are $X^a$, i.e.\ it is Schwartz in $Z$),
and $a_+$ is supported where $y$ is near $y'$, which in
turn implies (by looking at the phase of the exponential)
that $\tau$ is near $\sqrt{\lambda-\ep_\alpha}$. Note that
the phase is a multiple of $\cos\distance(y,y')=y\cdot y'$, where
we consider $C_a$ as the unit sphere in $X_a$, $y,y'\in C_a$. In terms
of Euclidean variables $w_a=y/x$ on $X^a$,
\begin{equation}
\Kf_{\alpha,+}=e^{-i\sqrt{\lambda-\ep_\alpha}w_a\cdot y'}
|w_a|^{-i\It_a(y')/2\sqrt{\lambda-\ep_\alpha}} \at_+(w_a,y',Z),
\end{equation}
with $\at_+$ polyhomogeneous, of degree $0$ in $w_a$, smooth in $y'$,
Schwartz in $Z$.

The construction of $\Kf_{\alpha,+}$, i.e.\ finding an oscillatory function of
the form \eqref{eq:Pt-kernel-structure} satisfying that
$(H-\lambda)\Kf_{\alpha,+}$ is in $\dCinf(\Xb\times C'_a)$ at least near
$y=y'$, where $H-\lambda$ is applied to $\Kf_{\alpha,+}$ in the left
factor $\Xb$,
is essentially given by a non-commutative WKB-type procedure.
Since one already has the correct phase function, one only needs
to solve transport equations starting at the diagonal
$y=y'$, which, however, are operator-valued.
To remove the errors in these equations,
$H^a-\ep_\alpha$ has to be inverted on the
orthocomplement of its $L^2$ null space. Since $H^a-\ep_\alpha$ has
a parametrix in $\PsiSc^{-2,0}(\Xb^a,\calC^a)$ (because $\ep_\alpha$
is below the set of thresholds), this is possible with the resulting
(generalized) inverse being in  $\PsiSc^{-2,0}(\Xb^a,\calC^a)$. This
gives a kernel $K_0$ of the form \eqref{eq:Pt-kernel-structure},
defined near $y=y'$, $x=0$, and satisfying $(H-\lambda)K_0$ is in
$\dCinf(\Xb\times C'_a)$ in the region where the transport equations have
been solved (see \cite{Skibsted:Smoothness, Bommier:Proprietes} for more
details).
Multiplying it by a cut-off function
$\phi\in\Cinf(C'_a\times C'_a)$ supported near the diagonal, identically
$1$ in a smaller neighborhood of the diagonal, we obtain a kernel $K$
which is still of the form \eqref{eq:Pt-kernel-structure}, now globally
defined, and $(H-\lambda)K$ is in $\dCinf(\Xb\times C'_a)$ near
$y=y'$. Finally, we take $\psi\in\Cinf_c(\Real)$ identically $1$ near
$\lambda$, and apply $\psi(H)\in\PsiSc^{-\infty,0}(\Xb,\calC)$ to
$K$ in the left factor $\Xb$, i.e.\ we compose the operator $P$ given
by the kernel $K$ with $\psi(H)$, so that
$\Pt_{\alpha,+}(\lambda)=\psi(H)P$. Writing
out the composition explicitly, similarly to how one shows that $\PsiSc(\Xb,
\calC)$ preserves oscillatory functions, see
\cite[Section~4]{Vasy:Propagation-Many},
shows that
$K^\flat_{\alpha,+}$ is indeed of the required form,
\eqref{eq:Pt-kernel-structure}. (Recall that this is based on conjugating
an Sc-ps.d.o.\ by an oscillatory function, observing
that the resulting kernel is that of an Sc-ps.d.o, and using that Sc-ps.d.o's
map smooth functions on $\Xb$ to smooth functions,
which is quite straightforward to see by an explicit integral representation
of the kernel.) In addition,
$(H-\lambda)\Pt_{\alpha,+}(\lambda)g=\psi(H)(H-\lambda)Pg$ for all
$g\in\dist_c(C'_a)$, which will allow us to use that
$(H-\lambda)K$ is in $\dCinf(\Xb\times C'_a)$ near $y=y'$. Also,
$(P-\Pt_{\alpha,+}(\lambda))g=(\Id-\psi(H))Pg=\psit(H)(H-\lambda)Pg$
where $\psit(s)=(1-\psi(s))/(s-\lambda)$ so $\psit\in S^{-1}_{\phg}(\Real)$,
$\psit(H)\in\PsiSc^{-2,0}(\Xb,\calC)$. In view of the corresponding
property of $(H-\lambda)K$, $K^\flat_{\alpha,+}-K$ vanishes to infinite
order near $y=y'$, $x=0$, i.e.\ $K^\flat_{\alpha,+}$ and $K$ are `microlocally
equal' there.

We note that in fact the construction of the kernel can be continued
up to $C_{a,\sing}$; we just need to cut off $a_+$ shortly beforehand.
In particular, if $a$ is a 2-cluster, then the construction works
until reaching the outgoing radial set.
We also remark that if the $V_b$ are Schwartz, then a product decomposition
immediately gives $\Pt_{\alpha,+}(\lambda)$, without a need to remove
errors, so in that case one can allow $\ep_\alpha$ to be arbitrary.

The main properties of $\Pt_{\alpha,+}(\lambda)$ are that
for any $g\in\dist_c(C'_a)$, $h\in\Cinf_c(C'_a)$, $\ep'>0$,
\begin{equation}\label{eq:Pt-C'_a}
\WFSc(\Pt_{\alpha,+}(\lambda)g)\subset\sct_{C'_a}(C_a;\Xb),
\end{equation}
\begin{equation}\begin{split}\label{eq:Pt-smooth}
\Pt_{\alpha,+}(\lambda)h=e^{-i\sqrt{\lambda-\ep_\alpha}/x}x^{\dim C_a/2}
x^{i\It_a/2\sqrt{\lambda-\ep_\alpha}}
&(((\pi^a)^*\psi_\alpha)((\pi_a^*h)+xv'),\\
&v'\ \text{as in}\ \eqref{eq:full-asymp-exp},
\end{split}\end{equation}
\begin{equation}\label{eq:Pt-smooth-error}
(H-\lambda)\Pt_{\alpha,+}(\lambda)h\in\dCinf(\Xb),
\end{equation}
\begin{equation}\begin{split}\label{eq:Pt-WF}
\WFSc(\Pt_{\alpha,+}(\lambda)g)\subset&
\{(y,\tau,\mu):\  y\in\supp g,\ \tau=\sqrt{\lambda-\ep_\alpha},\ \mu=0\}\\
&\quad\cup\{\xi\in\dot\Sigma(\lambda)\setminus (R_+(\lambda)\cup R_-(\lambda))
:\ \exists\zeta\in\WF(g),\ \xi\in\gamma_{\alpha,-}(\zeta)\},
\end{split}\end{equation}
\begin{equation}\begin{split}\label{eq:Pt-error}
\WFSc((H&-\lambda)\Pt_{\alpha,+}(\lambda)g)\\
&\quad\subset\{\xi\in\dot\Sigma(\lambda)
\setminus(R_+(\lambda)\cup R_-(\lambda)):
\ \exists\zeta\in\WF(g),\ \xi\in\gamma_{\alpha,-}(\zeta)\}.
\end{split}\end{equation}
See \cite[Appendix~A]{Vasy:Propagation-2} for a discussion of these mapping
properties if $\alpha$ is the free cluster; the general case is similar
since we are working in the regular part of $C_a$, and $a_+$ is rapidly
decreasing in $X^a$.
However, since the notion of wave front set is more
complicated than at the free cluster, we briefly indicate how to prove
these statements. Before doing this, we discuss the consequences of
\eqref{eq:Pt-C'_a}-\eqref{eq:Pt-error} for the actual Poisson
operator $P_{\alpha,+}(\lambda)$.

The Poisson operator, first defined on $\Cinf_c(C'_a)$ by \eqref{eq:P-def-16},
is given by
\begin{equation}\label{eq:P-Pt-rel}
P_{\alpha,+}(\lambda)=\Pt_{\alpha,+}(\lambda)-R(\lambda+i0)(H-\lambda)
\Pt_{\alpha,+}(\lambda).
\end{equation}
Note that for $g\in\Cinf_c(C'_a)$, $(H-\lambda)
\Pt_{\alpha,+}(\lambda)g\in\dCinf(\Xb)$ by \eqref{eq:Pt-smooth-error},
so the resolvent can certainly
be applied to it. For $g\in\Cinf_c(C'_a)$, the right hand side
of \eqref{eq:P-Pt-rel}
is of the form \eqref{eq:intro-21-S}
by \eqref{eq:Pt-smooth}, so it indeed yields
the Poisson operator by its definition, \eqref{eq:P-def-16}.
Similar results hold for the Poisson operator $P_{\beta,-}(\lambda)$
with outgoing initial data. Theorem~\ref{thm:res-WF}, \eqref{eq:Pt-error} and
\eqref{eq:P-Pt-rel} immediately show the following proposition.

\begin{prop*}(Proposition~\ref{prop:Poisson-extend})
Suppose that $H$ is a many-body Hamiltonian, $\lambda\nin\Lambda$,
and \eqref{eq:param-Poisson} holds.
Suppose also that $g\in\Cinf_c(C'_a)$. Then
\begin{equation}\label{eq:Poisson-WF-smooth-p}
\WFSc(P_{\alpha,+}(\lambda)g)
\setminus R_+(\lambda)\subset\Phi_+(R_-(\lambda)).
\end{equation}

In addition, $P_{\alpha,+}(\lambda)$ extends by continuity from
$\Cinf_c(C'_a)$ to distributions
$g\in\dist_c(C'_a)$ with $(R_+(\lambda)\times\WF(g))\cap\calR_{\alpha,+}
=\emptyset$ (i.e.\ $\calR_{\alpha,+}(\WF(g))\cap R_+(\lambda)=\emptyset$).
If $g$ is such a distribution, then
\begin{equation}
\WFSc(P_{\alpha,+}(\lambda)g)\setminus R_+(\lambda)
\subset\Phi_+(R_-(\lambda))\cup
\calR_{\alpha,+}(\WF(g)).
\end{equation}
\end{prop*}

We now return to the proof of \eqref{eq:Pt-C'_a}-\eqref{eq:Pt-error}.
First note that \eqref{eq:Pt-C'_a} follows from the fact that $\supp a_+$
is disjoint from $C_{a,\sing}\times D\times\Xb^a$ for any
$D\subset C'_a$ compact.
Next, \eqref{eq:Pt-smooth}-\eqref{eq:Pt-smooth-error}
are direct consequences of the stationary
phase lemma. The
$\log x$ factors appear in the lower order terms (i.e.\ in $v'$) due to the
derivatives falling on $x^{i\It_a(y')/2\sqrt{\lambda-\ep_\alpha}}$
as stated in \eqref{eq:full-asymp-exp};
see \cite[Sections~11-12]{Vasy:Propagation-Many}
for a detailed description. Note
that in addition to the stationary phase lemma,
\eqref{eq:Pt-smooth-error} uses that
$(H-\lambda)\Pt_{\alpha,+}(\lambda)h=\psi(H)(H-\lambda)Ph$, and the kernel
of $(H-\lambda)P$ is rapidly decreasing near $y=y'$, where the phase
is stationary. The stationary phase lemma also shows \eqref{eq:Pt-WF}
with $h\in\Cinf_c(C'_a)$ in place of $g$ (with the second set on
the right hand side empty), since all terms in the stationary phase
expansion arise by differentiating $h$ and evaluating the result at $y$,
so we obtain rapid decay away from $\supp h$.

Concerning \eqref{eq:Pt-WF} (for $g\in\dist_c(C'_a)$), we are done
if we show that for any $\xi$ that is not
in the right hand side of \eqref{eq:Pt-WF}, we can write
$u=Bv+f$ with $B\in\PsiSc^{-\infty,0}(\Xb,\calC)$, $v\in\dist(X)$, $f\in
\dCinf(\Xb)$, $\xi\nin\WFScp(B)$. For the sake of simplicity, and
since we only need this weaker statement here,
we only state the proof with
$\{(y,\sqrt{\lambda-\ep_\alpha},0):\ y\in\supp g\}$ in the right hand side
of \eqref{eq:Pt-WF} replaced by
$\{(y,\sqrt{\lambda-\ep_\alpha},0):\ y\in C'_a\}$; the improved
statement only requires a more careful consideration of
supports in \eqref{eq:qc-relation-0}.

To find such $B$, $v$, $f$, we consider
$q\in\Cinf_c(\sct \Xb_a)$ as in
Section~\ref{sec:commutators}, $\rho\in\Cinf(\Xb)$ supported near
$C_a$, identically $1$ in a small neighborhood of $C_a$, and take
$\psi_0\in\Cinf_c(\Real)$ is identically $1$
on $\supp\psi$ (recall that $\Pt_{\alpha,+}(\lambda)=\psi(H)P$),
and let $A=\rho Q\psi_0(H)\in\PsiSc^{-\infty,0}(\Xb,\calC)$
as in Section~\ref{sec:commutators}.
Then, as discussed above (with $\psi(H)$ in place of $A$), the kernel $K'$ of
$A\Pt_{\alpha,+}(\lambda)$ has the form 
\eqref{eq:Pt-kernel-structure}, i.e.
\begin{equation}
K'=e^{-i\sqrt{\lambda-\ep_\alpha}y\cdot y'/x}
x^{i\It_a(y')/2\sqrt{\lambda-\ep_\alpha}}a'(x,y,Z,y'),
\end{equation}
$a'\in\Cinf([\Xb;C_a]\times C_a)$, and in addition
\begin{equation}\begin{split}\label{eq:csupp-a'}
(\yb,d_{(y/x)}&(-\sqrt{\lambda-\ep_\alpha}\,y\cdot y'/x)|_{(0,\yb,\yb')})
\nin\supp q\\
&\Rightarrow a'\ \text{vanishes to infinite order near}
\ (\yb,\yb')\in C_a\times C'_a.
\end{split}\end{equation}
Suppose also that
\begin{equation}\label{eq:qc-relation-0}
y\in C'_a\Rightarrow (y,\sqrt{\lambda-\ep_\alpha},0)\nin\supp q,
\end{equation}
so by \eqref{eq:csupp-a'}, the diagonal $y=y'$ is disjoint from $\supp a'$.
This means that, with $\Fr_{X_a}$ denoting the Fourier transform in $X_a$,
$\Fr_{X_a} A\Pt_{\alpha,+}(\lambda)$
is a Fourier integral operator in the usual sense (on $X_a^*\times C'_a$
where $X_a^*$ is the dual of $X_a$, which we identify with $X_a$ for
the sake of convenience), with values in Schwartz functions on $X^a$,
and with canonical relation in
\begin{equation}\label{eq:can-rel-nzero}
(T^*X_a\setminus 0)\times(T^*C'_a\setminus 0)
\end{equation}
due to the previous remark.
We also
take $C\in\Psi^0(C'_a)$ with compactly supported kernel and amplitude
(full symbol) $c$, and suppose that
\begin{equation}\begin{split}\label{eq:qc-relation}
(y,\xi)&\in\supp q\Mand (y',\zeta)\in\csupp c\cap S^*C'_a\\
&\Rightarrow
(y,\xi)\nin \overline{\gamma_{\alpha,-}(y',\zeta)}=
\gamma_{\alpha,-}(y',\zeta)\cup\{(y',\sqrt{\lambda-\ep_\alpha},0\}.
\end{split}\end{equation}
(N.B. $(y,\xi)\in\supp q$ automatically
implies that $(y,\xi)$ is not in the second set
on the right hand side, i.e.\ it is not incoming,
by \eqref{eq:qc-relation-0}.)
We claim that under these conditions
$A\Pt_{\alpha,+}(\lambda)C$
has a smooth rapidly decreasing kernel $K_{APC}\in\dCinf(\Xb\times C'_a)$.
This can be seen from the standard FIO calculus (with values in
Schwartz functions on $X^a$ if stated this way) if we
perform Fourier transform in $X_a$. Namely, as remarked above
$\Fr_{X_a} A\Pt_{\alpha,+}(\lambda)$
is a Fourier integral operator, with canonical relation
satisfying \eqref{eq:can-rel-nzero}.
Thus, $\Fr_{X_a} A\Pt_{\alpha,+}(\lambda)C$ is
also an FIO.
But by \eqref{eq:csupp-a'} and \eqref{eq:qc-relation},
$(y',\zeta)\in\WF'(C)$ implies that for no $(\xi,y')\in T^*X_a$ is
$(\xi,y',y,\zeta)$ in
the canonical relation of $\Fr_{X_a} A\Pt_{\alpha,+}(\lambda)$,
hence
smoothness of the kernel of $\Fr_{X_a} A\Pt_{\alpha,+}(\lambda)C$ follows.
Rapid decay at infinity of the kernel of $\Fr_{X_a} A\Pt_{\alpha,+}(\lambda)C$
being automatic, this shows that
$K_{APC}\in\dCinf(\Xb\times C'_a)$.
This can also be seen more
directly, writing out the composition as
an oscillatory integral, in the uncompactified notation, with
$w_a=y/x$,
\begin{equation}\begin{split}
K_{APC}(w_a,Z,y'')=
&\int e^{-i\sqrt{\lambda-\ep_\alpha}w_a\cdot y'}
e^{i(y'-y'')\cdot\zeta}|w_a|^{-i\It_a(y')/2\sqrt{\lambda-\ep_\alpha}}\\
&\qquad a'(w_a,y',Z)c(y',\zeta)\,dy'\,d\zeta,
\end{split}\end{equation}
where $w_a\cdot y'$ stands for the
$X_a$-scalar product of $w_a\in X_a$ and $y'\in C_a$,
we see that the phase of the $y'$ integral is only stationary if
$(y,\xi_a)\in \overline{\gamma_{\alpha,-}(y',\zeta)}$,
but $a'c$ is rapidly decreasing there.

Now, given any $g\in\dist_c(C'_a)$, and any $U\subset S^*C'_a$ open such
that $\WF(g)\subset U$, we can find $C\in\Psi^0(C'_a)$ with
compactly supported kernel such that $\WFp(C)\subset U$, $\WFp(\Id-C)\cap
\WF(g)=\emptyset$.
Thus, $(\Id-C)g\in\Cinf_c(C'_a)$. Choosing $q$ as above,
i.e.\ satisfying \eqref{eq:qc-relation-0} and \eqref{eq:qc-relation},
we see that
\begin{equation}\label{eq:Pt-WF-16}
Au=A\Pt_{\alpha,+}(\lambda)Cg+A\Pt_{\alpha,+}(\lambda)(\Id-C)g\in\dCinf(\Xb).
\end{equation}
As $(\Id-\psi_0(H))u=0$, this shows that
\begin{equation}
u=\psi_0(H)u=Q\psi_0(H)u+(\psi_0(H)-Q\psi_0(H))u=Au+(\psi_0(H)-A)u,
\end{equation}
so if $q\equiv 1$ near $(y,\xi)$, so $(y,\xi)
\nin\WFScp(B)$, $B=\psi_0(H)-A\in\PsiSc^{-\infty,0}(\Xb,\calC)$,
then $(y,\xi)\nin\WFSc(u)$. But $U$, hence $\supp c$,
can be chosen arbitrarily small
around $\WF(g)$; correspondingly for any $(y,\xi)$ such that
$(y,\xi)\nin \overline{\gamma_{\alpha,-}(y',\zeta)}$
for any $(y',\zeta)\in\WF(g)$, and $(y,\xi)$ is not incoming,
we deduce that $(y,\xi)\nin\WFSc(u)$.
The estimate for $\WFSc((H-\lambda)u)$ works the same way; now
$(H-\lambda)u=\psi_0(H)(H-\lambda)Pg$, with $(H-\lambda)P$ of the same form
as $P$ but with smaller (cone) support. This gives the improved estimate.

Note if we only required finitely many of the
symbol estimates on, say, the symbol whose left quantization is $B$,
a conclusion analogous to \eqref{eq:Pt-WF}
would easily follow for any operator with
kernel given by an oscillatory function as in \eqref{eq:Pt-kernel-structure},
i.e.\ we could allow $a_+$ to be arbitrary, and there would be
no need for the factor $\psi(H)$ in the definition of
$\Pt_{\alpha,+}(\lambda)$. (Of course, \eqref{eq:Pt-error}
could not hold then.)
Indeed, $v=(\Delta+1)^k u$ is of the same form as $u$,
so $u=(\Delta+1)^{-k}v$. Choosing $q\in\Cinf_c(\sct \Xb_a)$, $C\in
\Psi^0(C'_a)$ as above, so $Q\in\Psisc^{-\infty,0}(\Xb_a)$, regarded
as an operator on $\Xb$, we deduce just as above that
$Q\Pt_{\alpha,+}(\lambda)C$ has a smooth rapidly decreasing kernel,
so \eqref{eq:Pt-WF-16} holds with $Q$ in place of $A$. Let $\rho\in\Cinf(\Xb)$
be a cutoff function supported near $C_a$ as above.
Then
$u=(\Id-\rho Q)(\Delta+1)^{-k}v+f$, $f=\rho Qu\in\dCinf(\Xb)$, so if
$B=(\Id-\rho Q)(\Delta+1)^{-k}$ were in $\PsiSc^{-\infty,0}(\Xb,\calC)$,
we would be able to deduce that
$\xi\nin\WFScp(B)$ implies $\xi\nin\WFSc(u)$. The problem
is that $q$ is not a symbol in $(\xi_a,\xi^a)$ jointly since it is independent
of the latter. However, note that
$k$ can be chosen arbitrarily large, so it can be arranged that
$B$ satisfies an arbitrary, but finite, number of the symbol estimates
corresponding to $\PsiSc^{-\infty,0}(\Xb,\calC)$, while keeping $v$ in
a fixed weighted Sobolev space. While this is sufficient for any
ps.d.o.\ argument based on weighted Sobolev space
estimates, etc., to go through, it is not sufficient
for ensuring that $B$ is in $\PsiSc^{-\infty,0}(\Xb,\calC)$; that was
the reason for the somewhat more complicate argument of the previous
paragraphs.

\section{Scattering matrices}\label{sec:S-matrix}
In view of the pairing formula \eqref{eq:equiv-7-S} for the S-matrix,
which we
restate below,
Proposition~\ref{prop:Poisson-extend} allows us to describe the wave front
relation of the S-matrices rather directly. Recall that
the S-matrix from an incoming cluster $\alpha$ to
an outgoing cluster $\beta$ is given by
\begin{equation}\label{eq:equiv-7-S-1}
S_{\beta\alpha}(\lambda)=\frac{1}{2i\sqrt{\lambda-\ep_\beta}}
((H-\lambda)\Pt_{\beta,-}(\lambda))^*P_{\alpha,+}(\lambda),
\end{equation}
$\lambda>\max(\epsilon_\alpha,\epsilon_\beta)$, $\lambda\nin\Lambda$.
This formula defines $S_{\beta\alpha}(\lambda)$ as a map
$\Cinf_c(C'_a)\to\dist(C'_b)$. Namely if $g\in\Cinf_c(C'_a)$ then
$S_{\beta\alpha}(\lambda)g$ is given by the formula
\begin{equation}\label{eq:S-matrix-pair-2}
\langle S_{\beta\alpha}(\lambda)g,h\rangle=\frac{1}{2i\sqrt{\lambda-\ep_\beta}}
\langle P_{\alpha,+}(\lambda)g,(H-\lambda)\Pt_{\beta,-}(\lambda)h\rangle,
\quad h\in\Cinf_c(C'_b),
\end{equation}
where the pairing makes sense since, by \eqref{eq:Pt-error} with
$\Pt_{\beta,-}$ in place of $\Pt_{\alpha,+}$,
$(H-\lambda)\Pt_{\beta,-}(\lambda)h\in\Sch(\Rn)$. (Note that $\langle\cdot,
\cdot\rangle$ denotes the complex $L^2$ pairing, linear in the first
variable, not the real distributional
pairing; the two are related by taking the complex conjugate of the
second term.)

Before stating the full results on the wave front relation of
$S_{\beta\alpha}(\lambda)$,
we first address the question whether $S_{\beta\alpha}(\lambda)$ maps
$\Cinf_c(C'_a)$ to $\Cinf(C'_b)$ and extends to a map
$\dist_c(C'_a)\to\dist(C'_b)$.
This is closely related to the question whether
the wave front set of the kernel of $S_{\beta\alpha}(\lambda)$ is
a subset of $T^*(C'_b\times C'_a)\setminus ((0\times T^*C'_a)\cup
(T^* C'_b\times 0))$, as we discussed in Section~\ref{sec:results}.
By a pairing argument we can answer the first question affirmatively
if for all $g\in\Cinf_c(C'_a)$, $h\in\dist_c(C'_b)$,
the two terms being paired in
\eqref{eq:S-matrix-pair-2}, that is
$P_{\alpha,+}(\lambda)g$ and $(H-\lambda)\Pt_{\beta,-}(\lambda)h$,
have disjoint wave front sets. Indeed, suppose that
$K_1=\WFSc(P_{\alpha,+}(\lambda)g)$
and $K_2=\WFSc((H-\lambda)\Pt_{\beta,-}(\lambda)h)$ are disjoint.
Being wave front sets, they are closed, and as $\psi(H)P_{\alpha,+}(\lambda)g
=P_{\alpha,+}(\lambda)g$ if $\psi\equiv 1$ near $\lambda$, $\psi\in
\Cinf_c(\Real)$, we see that $K_1$ is compact, and $K_2$ is also
compact by \eqref{eq:Pt-error} applied for $\Pt_{\beta,-}(\lambda)$.
As discussed in Section~\ref{sec:indicial}, there exists
$q\in\Cinf(\sct\Xb)$, $q$ $\pi$-invariant, such that $q\equiv 1$ on
$K_1$, $q\equiv 0$ on $K_2$, and $q$ satisfies estimates \eqref{eq:comm-4-a},
so we can define an operator $A\in\PsiSc^{-\infty,0}(\Xb,\calC)$
as in Lemma~\ref{lemma:comm-1} which is of the form $Q\psi(H)$,
$\psi$ as above, and $Q$ acts on oscillatory functions. Then
$\WFScp(\psi(H)-A)\cap K_1=\emptyset$ by construction,
so $(\psi(H)-A)P_{\alpha,+}(\lambda)g\in\dCinf(\Xb)$. In addition,
$\WFScp(A^*)\cap K_2=\emptyset$, so $A^*(H-\lambda)\Pt_{\beta,-}(\lambda)h
\in\dCinf(\Xb)$ as well. Then
\begin{equation}\begin{split}
&\langle P_{\alpha,+}(\lambda)g,(H-\lambda)\Pt_{\beta,-}(\lambda)h\rangle\\
&=\langle (\psi(H)-A)P_{\alpha,+}(\lambda)g,
(H-\lambda)\Pt_{\beta,-}(\lambda)h\rangle+
\langle P_{\alpha,+}(\lambda)g,A^*(H-\lambda)\Pt_{\beta,-}(\lambda)h\rangle,
\end{split}\end{equation}
with the equality a priori valid for $g$, $h$ smooth, for then the second
factor is Schwartz, shows that the pairing \eqref{eq:S-matrix-pair-2} makes
sense (i.e.\ extends by continuity from smooth data)
if $K_1$ and $K_2$ are disjoint.

This is certainly the case if $\alpha$ and
$\beta$ are either the free channel or 2-cluster channels
satisfying \eqref{eq:param-Poisson}
(including the possibility that one is free, the other is
such a two cluster channel).
In fact, if $\beta=0$ is the free channel, $\alpha$ is arbitrary
(but satisfies \eqref{eq:param-Poisson}), $\lambda>0$ (so with the
notation of Theorem~\ref{thm:lim-absorb} $a(\lambda)=\lambda$) then
on $\WFSc((H-\lambda)\Pt_{\beta,-}(\lambda)h)$, $-\sqrt{\lambda}<\tau
<\sqrt{\lambda}$ by \eqref{eq:Pt-error}, while
$\WFSc(P_{\alpha,+}(\lambda)g)
\subset \Phi_+(R_-(\lambda))$ by Proposition~\ref{prop:Poisson-extend},
so as $\tau$ is a decreasing function along generalized broken
bicharacteristics, $\tau\leq-\sqrt{a(\lambda)}=-\sqrt{\lambda}$ on it,
so these two sets
are indeed disjoint. On the other hand, if $\beta$ is a 2-cluster channel
satisfying \eqref{eq:param-Poisson}, $\alpha$ arbitrary,
then one can arrange (by continuing
the kernel construction until close to the $-$ outgoing, i.e.\ $+$
incoming, set) that
on $\WFSc((H-\lambda)\Pt_{\beta,-}(\lambda)h)$, $\tau$ is arbitrarily
close to $\sqrt{\lambda-\ep_\beta}$, which is again disjoint from
the set $\tau\leq-\sqrt{a(\lambda)}$, hence from
$\WFSc(P_{\alpha,+}(\lambda)g)\subset \Phi_+(R_-(\lambda))$.
This suffices to prove that $S_{\beta\alpha}(\lambda):\Cinf(C'_a)
\to\Cinf(C'_b)$ in this case. By a similar argument, or by duality,
$S_{\beta\alpha}(\lambda)$ extends to a map $\dist(C'_a)\to\dist(C'_b)$
if either $\alpha$ is a 2-cluster channel or the free channel.
In other cases Proposition~\ref{prop:Poisson-extend} and
\eqref{eq:Pt-error} do not allow us to conclude that the wave front sets
of the two terms are disjoint, and correspondingly we cannot expect that
$S_{\beta\alpha}(\lambda)$ either preserves smoothness or extends to
distributional data. We summarize these results in the following corollary.

\begin{cor}
Suppose that $H$ is a many-body Hamiltonian, $\lambda\nin\Lambda$,
and \eqref{eq:param-Poisson} holds. If $\beta$ is a two-cluster channel
or the free cluster then $S_{\beta\alpha}(\lambda)$ preserves
smoothness: $S_{\beta\alpha}(\lambda):\Cinf_c(C'_a)\to
\Cinf_c(C'_b)$. If $\alpha$ is a two-cluster channel
or the free cluster then $S_{\beta\alpha}(\lambda)$ extends
to a map on distributions: $S_{\beta\alpha}(\lambda):\dist_c(C'_a)\to
\dist_c(C'_b)$
\end{cor}

The pairing argument of course allows us to describe the wave front
relation of the S-matrices in general. If $g\in\dist_c(C'_a)$
satisfies $\calR_{\alpha,+}(\WF(g))\cap R_+(\lambda)=\emptyset$, then
$P_{\alpha,+}(\lambda)g$ is defined, and its wave front set is given
by Proposition~\ref{prop:Poisson-extend}. Similarly, if $h\in
\dist_c(C'_b)$, then $\Pt_{\beta,-}(\lambda)h$ is always defined
and the wave front set of its `error as a generalized
eigenfunction' is given by \eqref{eq:Pt-error}.
Thus, under the additional assumption that these two wave front sets
are disjoint, which holds if
\begin{equation}\label{eq:g-h-WF-pair}
\{\gamma_{\beta,+}(\zeta):\zeta\in \WF(h)\}
\cap(\calR_{\alpha,+}(\WF(g))\cup\Phi_+(R_-(\lambda)))
=\emptyset,
\end{equation}
we see that $\langle S_{\beta\alpha}(\lambda)g,h\rangle$ is defined via
\eqref{eq:S-matrix-pair-2}. Since $h$ is arbitrary, except that
its wave front set satisfies \eqref{eq:g-h-WF-pair}, this allows us
to conclude that if $K\subset S^*C'_a$ is such that
\begin{equation}\label{eq:g-K-WF-pair}
\{\gamma_{\beta,+}(\zeta):\zeta\in K\}
\cap(\calR_{\alpha,+}(\WF(g))\cup\Phi_+(R_-(\lambda)))
=\emptyset,
\end{equation}
then $\WF(S_{\beta\alpha}(\lambda)g)\cap K=\emptyset$.
But by the definition of $\calR_{\beta\alpha}$ and $\calR_{\beta,-}$,
\eqref{eq:g-K-WF-pair} is equivalent to
\begin{equation}
\calR_{\beta\alpha}(\WF(g))\cap K=\emptyset\Mand
\calR_{\beta,-}^{-1}(R_-(\lambda))\cap K
=\emptyset.
\end{equation}
(See the proof of Theorem~12.4 in \cite{Vasy:Propagation-Many} for a
more detailed version of the pairing argument.)
We have thus completed the proof of
Theorem~\ref{thm:sc-matrix}.

\begin{thm*}(Theorem~\ref{thm:sc-matrix})
Suppose that $H$ is a many-body Hamiltonian, $\lambda\nin\Lambda$,
and \eqref{eq:param-Poisson} holds.
Then $S_{\beta\alpha}(\lambda)$ extends by continuity from
$\Cinf_c(C'_a)$ to distributions
$g\in\dist_c(C'_a)$ with $\calR_{\alpha,+}(\WF(g))\cap R_+(\lambda)=\emptyset$.
If $g$ is such a distribution, then
\begin{equation}
\WF(S_{\beta\alpha}(\lambda)g)\subset\calR^{-1}_{\beta,-}(R_-(\lambda))\cup
\calR_{\beta\alpha}(\WF(g)).
\end{equation}
\end{thm*}

\appendix

\section{Geometric description of broken bicharacteristics}
\label{app:geometric}
The purpose of the appendix is
to give a more geometric description of the generalized broken
bicharacteristics, provided that the set of thresholds, $\Lambda_1$,
is discrete, that is to prove Theorem~\ref{thm:geom-br-bichar}.
We also analyze four-body scattering, where we can also obtain a similar
simple geometric description for generalized broken bicharacteristics.

It is convenient to break up the analysis of the structure of the
generalized broken bicharacteristics into a local and a global step.
Although until now we used the expression `broken bicharacteristic'
as a synonym for the cumbersome term `generalized broken
bicharacteristic', we now give a stronger definition for the former
{\em which is only valid in this Appendix.}

\begin{Def}
Suppose $I$ is compact, $\gamma:I\to\dot\Sigma(\lambda)$ is a
generalized broken bicharacteristic. We say that it is a broken
bicharacteristic if there exist $t_0<t_1<\ldots<t_k$, $I=[t_0,t_k]$,
such that for each $j$,
$\gamma|_{[t_j,t_{j+1}]}$ is the projection of an integral curve
of $\scHg$ to $\dot\Sigma(\lambda)$. In general, if $I$ is not compact,
we say that $\gamma$ is a broken bicharacteristic if $\gamma|_J$ is
a broken bicharacteristic for all compact intervals $J\subset I$.
\end{Def}

\begin{Def}\hfill

\begin{enumerate}
\item
An $N$-body Hamiltonian has a locally broken bicharacteristic relation
if for all $I$ compact, and for all generalized broken bicharacteristics
$\gamma:I\to\dot\Sigma(\lambda)$,
there exist $t_0<t_1<\ldots<t_k$, $I=[t_0,t_k]$,
such that for each $j$,
$\gamma|_{[t_j,t_{j+1}]}$ is the projection of an integral curve
of $\scHg$ to $\dot\Sigma(\lambda)$.

\item
We say that an $N$-body Hamiltonian has a globally broken bicharacteristic
relation if for all generalized broken bicharacteristics
$\gamma:\Real\to\dot\Sigma(\lambda)$, there exist
$t_0<t_1<\ldots<t_k$,
such that $\gamma|_{(-\infty,t_0]}$, $\gamma_{[t_k,+\infty)}$, as well
as for each $j$,
$\gamma|_{[t_j,t_{j+1}]}$, are the projections of integral curves
of $\scHg$ to $\dot\Sigma(\lambda)$.
\end{enumerate}
\end{Def}

When we talk about the length $\ell(\gamma)$
of a generalized broken bicharacteristic
$\gamma$, we mean the length of its projection to $\Sn$ (recall that the
projection is Lipschitz).

\begin{lemma}\label{lemma:local-if-short}
Suppose that $H$ is a many-body Hamiltonian. Then there exists $l>0$
such that for every generalized broken bicharacteristic $\gamma:[\alpha,\beta]
\to\dot\Sigma(\lambda)$ of length $\ell(\gamma)\leq l$
there exists a cluster $a$
such that the image of $\gamma$ lies in $\cup_{C_a\subset C_b}
\sct_{C'_b}\Xb_b$.
\end{lemma}

\begin{proof}
For each point $p\in\Sn$ there exists a unique $a$ such that $p\in C'_a$.
Also, there exists $l_p>0$ such that $\exp_p(B_{l_p}(0))$ lies
in $\cup_{C_a\subset C_b}\sct_{C'_b}\Xb_b$, where $B_{l_p}(0)$ denotes
the open ball of radius $l_p$ in $T_p\Sn$ with respect to the standard metric.
Now, $\{\exp_p(B_{l_p/2}(0)):\ p\in\Sn\}$ is an open cover of the
compact set $\Sn$, so it has a finite subcover corresponding to some
points, say, $p_j$, $j=1,\ldots,m$. Let $l=\min\{l_{p_j}/2:1\leq j\leq m\}$.
If $\gamma$ is as above, let $j$ be such that $\gamma(\alpha)
\in \exp_{p_j}(B_{l_{p_j/2}}(0))$. For any $t\in[\alpha,\beta]$, the
distance of $\gamma(t)$ and $\gamma(\alpha)$ is bounded by the length
of $\gamma$, hence by $l\leq l_{p_j}/2$, so the image of $\gamma$ lies in
$\exp_{p_j}(B_{l_{p_j}}(0))$ which in turn lies in
$\cup_{C_a\subset C_b}\sct_{C'_b}\Xb_b$ for some $a$.
\end{proof}

We next give an upper bound for the length of a broken bicharacteristic
if $\Lambda_1$ is discrete. Note that $\Lambda_1$ is bounded, so under
this assumption it is finite.

\begin{lemma}\label{lemma:discrete-arclength}
Suppose that $\Lambda_1$ is discrete, and let $C_1$ be the number
of its elements.
Suppose that $\gamma:I\to\dot\Sigma(\lambda)$ is a broken bicharacteristic, $I$
an interval. Then $\ell(\gamma)\leq C_1\pi$.
\end{lemma}

\begin{proof}
We recall from \cite{RBMZw}
the explicit arclength parametrization, $s=s(t)$,
of the integral curves of $\scHg$ with kinetic energy $\sigma>0$.
In terms of this parameterization for $\gamma|_{[t_j,t'_j]}$,
$\tau(s)=\sqrt{\sigma_j}\,\cos(s-s_0)$
where $s$
varies in a subinterval of $(s_0,s_0+\pi)$, and $[t_j,t'_j]$ is
such that $\gamma|_{[t_j,t'_j]}$ is an integral curve of $\scHg$
with kinetic energy $\sigma_j$.
Since $\tau\circ\gamma$ (which we usually just write as $\tau$) is
monotone decreasing, this shows that the total length of the segments
of $\gamma$ which are integral curves of $\scHg$ with any given kinetic
energy $\sigma>0$, is at most $\pi$. If $\sigma=0$, the bicharacteristic
segment is constant. Since $\sigma$ must be such that
$\lambda-\sigma\in\Lambda_1$, there are $C_1$ possible values of $\sigma$,
which proves our estimate.
\end{proof}

Next, we note that in three-body scattering the uniform upper bound $\pi$
holds without any assumptions on the structure of $\Lambda_1$,

\begin{lemma}\label{lemma:3-body-pi}
Suppose that $H$ is a three-body Hamiltonian.
Suppose that $\gamma:I\to\Real$ is a broken bicharacteristic, $I$
an interval. Then $\ell(\gamma)\leq \pi$.
\end{lemma}

\begin{proof}
In three-body scattering kinetic energy is constant along
generalized broken bicharacteristics (essentially
because there are no positive energy bound states), so the proof
of the previous lemma applies and gives the desired conclusion.
\end{proof}

Our strategy to analyze generalized broken bicharacteristics is to
divide them into pieces of length $\leq l$, and reduce the analysis to
that in a subsystem, by virtue of Lemma~\ref{lemma:local-if-short}.

\begin{lemma}\label{lemma:loc-breaks-bounds}
Suppose $H$ is an $N$-body hamiltonian
all of its proper subsystems
have a globally broken bicharacteristic relation.
Then every generalized broken bicharacteristic of $H$ of length $\leq l$,
$l$ as in
Lemma~\ref{lemma:local-if-short},
is a broken bicharacteristic. In addition,
if in all of the proper subsystems
the maximum number of breaks in a broken bicharacteristic (defined over
$\Real$ and of {\em arbitrary energy}) is at most $M_{N-1}$, then
every generalized broken bicharacteristic of $H$ of length $\leq l$
has at most $2M_{N-1}+2$ breaks.
\end{lemma}

\begin{proof}
Let $\gamma:[\alpha,\beta]\to\dot\Sigma(\lambda)$ be
a generalized
broken bicharacteristic 
of length $\ell(\gamma)\leq l$, so the image of $\gamma$ lies in a region
$\cup_{C_a\subset C_b}
\sct_{C'_b}\Xb_b$ for some $a$.
Thus, we can use local coordinates around $a$ for
describing $\gamma$. By Lemma~\ref{lemma:br-bichar-tgt-normal} and
the argument preceeding it (showing that $\eta_a$ cannot change
sign more than once),
there exist some points $\alpha',\beta'$ such that on $[\alpha,\alpha')$
and on $(\beta',\beta]$ the image of $\gamma$ is disjoint from $\sct
_{C'_a}\Xb_a$, and on $[\alpha',\beta']$, $\gamma$ is an integral curve
of $\scHg^a$, where some intervals may be empty or reduce to a point.
Consider the interval $(\beta',\beta]$ for the sake of definiteness,
and assume that it is non-empty.
Let $S(t)$ be given by
\begin{equation}
S(t)=-\int^{\beta}_t
(1+|z_a(t')|^2)^{1/2}|z_a(t')|^{-1}\,dt',
\end{equation}
so $S(t)$ is the solution of the ODE
\begin{equation}
dS/dt=(1+|z_a(t)|^2)^{1/2}/|z_a(t)|,\ t\in(\beta',\beta],\ S(\beta)=0,
\end{equation}
where we wrote $|z_a(t)|=|z_a(\gamma(t))|$. Thus, $S$ is $\calC^1$ and its
derivative is positive, so the same holds for its inverse function, $S^{-1}$,
defined on an interval $J$.
We denote by
$\dot\Sigma_{H^a}(\lambda-|\xib_a|^2)$ the characteristic variety of the
(proper) subsystem Hamiltonian $H^a$ at energy $\lambda-|\xib_a|^2$.
Now let $\gammat:J\to\dot\Sigma_{H^a}(\lambda-|\xib_a|^2)$,
given by
\begin{equation}\label{eq:gl->loc-5}
\gammat(s)=(z_a(S^{-1}(s))/|z_a(S^{-1}(s))|,\nu_a(S^{-1}(s))),
\end{equation}
so in terms of Euclidean coordinates,
\begin{equation}\label{eq:gl->loc-7}
\gammat(s)=(w^a(S^{-1}(s))/|w^a(S^{-1}(s))|,\xi^a(S^{-1}(s))),
\end{equation}
It is straightforward to check that $\gammat$ is a generalized
broken bicharacteristic
of $H^a$; the change of parameters accounts for the change in normalization
of the rescaled Hamilton vector fields. In fact, $\scHg=\langle w\rangle H_g$,
while its analog in the subsystem is $\langle w^a\rangle H_{g^a}$, and the
quotient $\langle w^a\rangle/\langle w\rangle$ is $|z_a|/(1+|z_a|^2)^{1/2}$.
Thus, by the hypothesis, $\gammat$ is a broken bicharacteristic.
Since $\xi_a$ is constant along generalized broken bicharacteristics,
we conclude that on $(\beta',\beta]$, $\gamma|_{(\beta',\beta]}$
is a broken bicharacteristic
if and only if $\gammat$ is, so by the hypothesis it is a broken
bicharacteristic, and it has as many
breaks as $\gammat$, so if we assume uniform bounds in the proper
subsystems, at most $M_{N-1}$.

A similar estimate holds for
$[\alpha,\alpha')$, so on $[\alpha,\beta]$, $\gamma$ is a broken
bicharacteristic, and if we also assume uniform bounds in the proper
subsystems, it can have at
most $2M_{N-1}+2$ breaks, proving the lemma.
\end{proof}

\begin{cor}\label{cor:glob->loc}
Suppose $H$ is an $N$-body hamiltonian, and
all of its proper subsystems
have a globally broken bicharacteristic relation.
Then $H$ has a locally broken bicharacteristic relation.
\end{cor}

\begin{proof}
If $\gamma:I\to\dot\Sigma(\lambda)$ is a generalized broken bicharacteristic,
$I$ compact, then its total arclength is finite, so dividing it up into
segments of length $\leq l$ and applying the previous lemma
proves the proposition.
\end{proof}

\begin{prop}\label{prop:N-1->N}
Suppose $H$ is an $N$-body hamiltonian with $\Lambda_1$ discrete,
all of its proper subsystems
have a globally broken bicharacteristic relation,
and in all of these subsystems
the maximum number of breaks in a broken bicharacteristic (defined over
$\Real$ and {\em of any energy}) is at most $M_{N-1}$.
Then $H$ has a globally broken bicharacteristic
relation and there is a constant $M_N$ such that
every broken bicharacteristic of $H$ {\em of any energy}
has at most $M_N$ breaks.
\end{prop}

\begin{proof}
Due to the previous corollary, every generalized broken
bicharacteristic of $H$ is a broken bicharacteristic.
By Lemma~\ref{lemma:discrete-arclength}, its length $\ell(\gamma)$
satisfies the inequality $\ell(\gamma)\leq C_1\pi$.
Dividing it up into pieces
of length $\leq l$, of which (we can arrange that)
there are at most $C_1\pi/l+1$, shows that the total number of breaks is
at most $(C_1\pi/l+1)(2M_{N-1}+3)$ which is independent of $\gamma$,
proving the proposition.
\end{proof}

An inductive argument completes the proof of Theorem~\ref{thm:geom-br-bichar},
which we restate.

\begin{thm}
Suppose that $H$ is a many-body Hamiltonian and $\Lambda_1$ is
discrete. Then $H$ has a globally broken bicharacteristic relation, and there
is a constant $M$ such that every generalized broken bicharacteristic of
$H$ has at most $M$ breaks.
\end{thm}

We next discuss what happens if $\Lambda_1$ is not discrete.

In the following we let $C_0>0$ be a constant such that
for $s,s'\in(0,\pi)$,
\begin{equation}
|\cos s-\cos s'|\geq C_0^{-2} |s-s'|^2.
\end{equation}
Thus, $C_0^2$ is essentially the H\"older constant in the H\"older of order
$\half$ estimate for $\arccos$.

\begin{lemma}\label{lemma:tau-arclength}
Suppose that $\gamma:I\to\Real$ is a broken bicharacteristic, $I$ compact,
$[t_j,t'_j]$, $j=1,\ldots,k$, are subintervals of $I$ such that the
open intervals $(t_j,t'_j)$ are disjoint. Suppose that the kinetic
energy over the interval $[t_j,t'_j]$ is $\sigma_j$, and let
\begin{equation}
\Delta \tau=\sum_{j=1}^k\tau(\gamma(t'_j))-\tau(\gamma(t_j)).
\end{equation}
Let $\ell(\gamma_j)$
be the length of $\gamma_j=\gamma|_{[t_j,t'_j]}$. Then
\begin{equation}\begin{split}
\sum_{j=1}^k\ell(\gamma_j)
&\leq C_0 (\sum_{j=1,\ \sigma_j>0}^k\sigma_j^{-1})^{1/2}|\Delta\tau|\\
&\leq C_0(\min\{\sigma_j:\ \sigma_j\neq 0,\ j=1,\ldots,k\})^{-1/2}
\,k^{1/2}|\Delta\tau|.
\end{split}\end{equation}
\end{lemma}

\begin{proof}
We again use the explicit arclength parametrization, $s=s(t)$,
of the integral curves of $\scHg$ with $\sigma_j>0$. Thus,
$\tau(s)=\sqrt{\sigma_j}\,\cos(s-s_0)$ for
$\gamma_j$ where $s$
varies in a subinterval of $(s_0,s_0+\pi)$. This gives
\begin{equation}
C_0\sigma_j^{-1/2}
|\tau(s(t'_j))-\tau(s(t_j))|^{1/2}\geq  s(t'_j)-s(t_j).
\end{equation}
If $\sigma_j=0$ then both $\tau(s(t'_j))-\tau(s(t_j))$ and
$s(t'_j)-s(t_j)$ vanish.
Summing over $j$ and
applying the Cauchy-Schwartz inequality to the left hand side
proves the lemma.
\end{proof}

Note that in Lemma~\ref{lemma:loc-breaks-bounds}, we really only need
to assume that the estimates $M_{N-1}$
in the subsystems are uniform over
bounded sets of energy to obtain a uniform estimate $M_N$
over bounded sets of energy for $H$.

\begin{cor}\label{cor:glob-bds-not-disc}
Suppose $H$ is an $N$-body hamiltonian
all of its proper subsystems
have a globally broken bicharacteristic relation, and
in all of the proper subsystems
the maximum number of breaks in a broken bicharacteristic (defined over
$\Real$) is at most $M_{N-1}$.
Suppose that $\gamma:I\to\dot\Sigma(\lambda)$ is a generalized broken
bicharacteristic of $H$ (hence a broken
bicharacteristic under these assumptions), and suppose that
there exist $c_0>0$, $m>0$,
such that the kinetic energy assumes at most $m$ values which are
less than $c_0$. Then there is a constant $M_N$ that depends only
on $c_0$, $m$, and $H$, such that $\gamma$ has at most $M_N$ breaks.
\end{cor}

\begin{proof}
We only need to show that there exists $M_N$ as above such that
for every compact interval $J\subset I$, $\gamma|_J$ has at most
$M_N$ breaks. Let $n$ denote the number of breaks in $\gamma|_J$.
Note that $\tau\circ\gamma$ is
monotone decreasing and it is bounded, with a bound given by
$(\lambda-\inf\Lambda_1)^{1/2}$ (since $\lambda-\tau^2-|\mu|^2\in
\Lambda_1$ in $\dot\Sigma(\lambda)$). By Lemma~\ref{lemma:tau-arclength},
the total length of bicharacteristic segments with kinetic energy
at least $c_0$ is at most
$2C_0 c_0^{-1/2}n^{1/2}(\lambda-\inf\Lambda_1)^{1/2}$, while by
the proof of Lemma~\ref{lemma:discrete-arclength}, the total length
of bicharacteristic segments with kinetic energy
less than $c_0$ is at most $m\pi$. Thus,
\begin{equation}
\ell(\gamma|_J)\leq 2C_0 c_0^{-1/2}n^{1/2}(\lambda-\inf\Lambda_1)^{1/2}
+m\pi.
\end{equation}
On the other hand, dividing $\gamma|_J$ into segments of length
$\leq l$, of which we can arrange that there are at most $\ell(\gamma|_J)/l+1$,
and applying Lemma~\ref{lemma:loc-breaks-bounds}, shows that
\begin{equation}\begin{split}
n&\leq (2M_{N-1}+2)(\ell(\gamma|_J)/l+1)\\
&\leq (2M_{N-1}+2)(1+m\pi/l)
+2(2M_{N-1}+2)
C_0 c_0^{-1/2}n^{1/2}(\lambda-\inf\Lambda_1)^{1/2}/l.
\end{split}\end{equation}
Dividing through by $n^{1/2}$ gives the desired estimate for $n$,
uniform in the energy $\lambda$ as long as $\lambda$ stays in a bounded
set.
\end{proof}

\begin{prop}
Suppose that $H$ is as in Corollary~\ref{cor:glob-bds-not-disc},
and $\lambda\nin\Lambda_1$. Then there exists $M_N$ depending only on
$H$ and $\inf\{\lambda-E:\ \lambda-E>0,\ E\in\Lambda_1\}$
such that every generalized broken bicharacteristic
$\gamma:I\to\dot\Sigma(\lambda)$ has at most $M_N$ breaks.
\end{prop}

\begin{proof}
Note first that $\Lambda_1$ is closed, so for $\lambda\nin\Lambda_1$,
$\inf\{\lambda-E:\ \lambda-E>0,\ E\in\Lambda_1\}>0$. The kinetic
energy $\sigma\geq 0$ on any bicharacteristic segment satisfies
$\lambda-\sigma\in\Lambda_1$, so
$\sigma\geq\inf\{\lambda-E:\ \lambda-E>0,\ E\in\Lambda_1\}>0$. Applying
Corollary~\ref{cor:glob-bds-not-disc} with
$c_0=\inf\{\lambda-E:\ \lambda-E>0,\ E\in\Lambda_1\}$, $m=0$, completes
the proof.
\end{proof}

Combining this result with Lemma~\ref{lemma:3-body-pi} yields the
following theorem.

\begin{thm}
Suppose that $H$ is a $4$-body Hamiltonian, and $\lambda\nin\Lambda_1$.
Then there exists $M$ depending only on
$H$ and $\inf\{\lambda-E:\ \lambda-E>0,\ E\in\Lambda_1\}$
such that every generalized broken bicharacteristic
$\gamma:I\to\dot\Sigma(\lambda)$ has at most $M$ breaks.
\end{thm}

While in $4$-body scattering one is mostly interested in
$\lambda\nin\Lambda_1$, making this result useful, in the inductive step,
to analyze $5$-body scattering, the $4$-body thresholds must be understood
as well, as indicated by the proof of Lemma~\ref{lemma:loc-breaks-bounds},
in particular by $\gammat:J\to\dot\Sigma_{H^a}(\lambda-|\xib_a|^2)$
(the notation is the same as in the Lemma), for $\lambda-|\xib_a|^2$
can be a threshold of the subsystem $H^a$.

\section{Lagrangian structure of the forward broken bicharacteristic
relation}\label{app:Lagrangian}
In this section we show that if
the set of thresholds of $H$ is discrete,
then the forward bicharacteristic relation is the finite union of
Lagrangian submanifolds, and it maps the outgoing radial set to a finite
union of Lagrangian
submanifolds as well. Together with Theorem~\ref{thm:sc-matrix},
(a simple modification of) this immediately implies the corresponding
result, stated in Theorem~\ref{thm:sc-matrix-disc},
about the wave front relation of the
various scattering matrices with incoming and outgoing
channels satisfying the assumptions \eqref{eq:param-Poisson}.

First, recall that if $\Lambda_1$ is discrete,
every generalized broken bicharacteristic $\gamma$
consists of a finite number of bicharacteristic segments (with
a uniform bound on the number of these segments). Let $a_j$ be the
cluster and $\alpha_j$ the channel
in which the $j$th segment propagates.
Also, let $c_j$ be the cluster at which $\gamma$ breaks between
propagation along the clusters $a_j$ and $a_{j+1}$, so
$C_{c_j}\subset C_{a_j}\cap C_{a_{j+1}}$ and the `break point' is
in $\sct_{C'_{c_j}}\Xb_{c_j}$, i.e.\ over the regular part of $C_{c_j}$.
Thus, we can associate a string
\begin{equation}\label{eq:br-strings}
a_1,\alpha_1,c_1,a_2,\alpha_2,c_2,\ldots,c_m,a_{m+1},\alpha_{m+1}
\end{equation}
to $\gamma$, where each term has the same meaning as above. 
Since there is a uniform bound on the finite number of breaks that
$\gamma$ has, there are only finitely many such permissible strings.
Moreover, since we are assuming that there is a break at $c_j$, $j=1,\ldots,
m$, we can impose that either $a_j\neq c_j$ or $a_{j+1}\neq c_j$.
(Note that the requirement that the `break points' are
in $\sct_{C'_{c_j}}\Xb_{c_j}$ is only necessary to associate unique
$c_j$'s to $\gamma$; it plays no role in the
following arguments.)

It will be convenient to use the linear structure of the collision planes
in this discussion. Thus, we will think of the relations, just as of the wave
front sets, as conic sets, invariant under a natural $\Real^+$ action {\em in
the base variables} (i.e.\ in configuration space).
Correspondingly, we talk about conic
Lagrangian rather than Legendre submanifolds. Since multiplying a vector
field by a non-zero function simply reparameterizes the integral curves,
the bicharacteristics of $\scHg^a$ and $H_{g_a}$ agree after
reparamaterization, so the former are given by the ($\Real^+$-quotient
induced)
projection of straight lines in
$T^*X_a=X_a\times X_a^*$ whose projection to $X_a^*$ is constant, say $\xi$,
such that $|\xi|_a^2=g_a(\xi)
=\lambda-\ep_{\alpha}$ where $\ep_\alpha$ is the energy
of the bound state the particle is propagating in. We sometimes write
$|\xi|$ for $|\xi|_a$ to simplify the notation.

We show that the forward broken bicharacteristic relation is given by
a finite union of Lagrangians, one for each string as in
\eqref{eq:br-strings}. To see this, we write the part of the
forward broken bicharacteristic relation corresponding to such a string as
the composite of elementary relations, each of which is Lagrangian, and
which intersect transversally in the usual sense.
Thus, we let $\Lambdat=\Lambdat_{c_j,\alpha_j,a_j,c_{j-1}}$ to be the relation
on $T^*X_{c_j}\times T^*X_{c_{j-1}}$ corresponding to a bicharacteristic
emanating from the collision plane $X_{c_{j-1}}$,
propagating along $a_j$ in channel $\alpha_j$, and then hitting
the collision plane $X_{c_{j}}$. That is, $\Lambdat$ is given by the projection
of the end points of bicharacteristic segments of $H_{g_{a_j}}$
whose endpoints lie over $X_{c_{j-1}}$ and $X_{c_{j}}$ respectively,
to $T^*X_{c_j}\times T^*X_{c_{j-1}}$ via the projections
$\pi_{a_j,c_{j-1}}$ and $\pi_{a_j,c_j}$.
We identify
$X_a^*$ with $X_a$ via the metric $g_a$ (induced by $g$); we do so
in particular when talking about
$T^* X_a=X_a\times X_a^*$. Thus,
$(w,\xi,w',\xi')\in\Lambdat$ means that
\begin{equation}\label{eq:Lambdat-def-1}
\xi=\pi_{a_j,c_j}(\xit),\ \xi'=\pi_{a_j,c_{j-1}}(\xit),
\ \xit=\sqrt{\lambda-\ep_{\alpha_j}}\,\frac{w-w'}{|w-w'|_{a_j}}.
\end{equation}
In particular, this shows that $\Lambdat$ is smoothly parameterized by
$w\in X_{c_j}$, $w'\in X_{c_{j-1}}$, $w\neq w'$,
so $\Lambdat$ is a graph over the
base $X_{c_j}\times X_{c_{j-1}}$ of the cotangent bundle (away from
$w=w'$). Note that if $w=w'$, the segment of $\gamma$ is trivial,
hence $\gamma$ is described by a different string, namely one from
which $\alpha_j$ is missing.
Thus, the tangent space of $\Lambdat$ at some point $p\in\Lambdat$,
given by $(w,w')$ as above,
is spanned by
\begin{equation}\begin{split}\label{eq:Lambdat-5}
&v\cdot\pa_w+B_j v\cdot\pa_{\xi}+C_j v\cdot\pa_{\xi'},\quad v\in X_{c_j},\\
&v'\cdot\pa_{w'}+B'_j v'\cdot\pa_{\xi'}+C'_j
v'\cdot\pa_{\xi},\quad v'\in X_{c_{j-1}},
\end{split}\end{equation}
where $B_j\in\End(X_{c_j})$, $B'_j\in\End(X_{c_{j-1}})$ are given by
\begin{equation}\begin{split}\label{eq:Lambdat-8}
&B_j=\frac{\sqrt{\lambda-\ep_{\alpha_j}}}{|w-w'|^3}(|w-w'|^2\,\Id
-(w-\pi_{a_j,c_j}(w'))\otimes (w-\pi_{a_j,c_j}(w'))),\\
&B'_j=-\frac{\sqrt{\lambda-\ep_{\alpha_j}}}{|w-w'|^3}(|w-w'|^2\,\Id
-(w'-\pi_{a_j,c_{j-1}}(w))\otimes (w'-\pi_{a_j,c_{j-1}}(w))).
\end{split}\end{equation}
Here, as in what follows, the center dot $\cdot$, e.g.\ when
writing $v\cdot\pa_{\xi'}$, etc., simply denotes the $\pa_{\xi'}$, etc.,
component of the tangent vector. Note that $\Lambdat$ is a Lagrangian
submanifold of $T^*X_{c_j}\times T^*X_{c_{j-1}}$ with respect
to the usual twisted symplectic form $\omega_{c_{j}}-\omega_{c_{j-1}}$,
where $\omega_c$ denotes the canonical symplectic form on $T^* X_c$
as well as its lift to the product (by an abuse of notation).
As we discuss later, this implies $C'_j=C_j^*$ as well as $B_j$, $B'_j$
self-adjoint. In addition, from \eqref{eq:Lambdat-8},
$|w-\pi_{a_j,c_j}(w')|\leq |w-w'|$ shows that $B_j\geq 0$, and indeed
it is positive definite if $|w-\pi_{a_j,c_j}(w')|\neq |w-w'|$, i.e.\ if
$w'\nin X_{c_j}$. If $w'\in X_{c_j}$, then so is $\xit=\xi$, and we see
that the null-space of $B_j$ is $\Span\{\xi\}$, its range the orthocomplement,
$\Span\{\xi\}^\perp$, of $\Span\{\xi\}$ in $X_{c_j}$.
Similarly, $B'_j\leq 0$, and it is negative definite
if $w\nin X_{c_{j-1}}$.

With this notation, every element of the part of the forward generalized broken
bicharacteristic relation
corresponding
to the string \eqref{eq:br-strings} is in the composite relation
\begin{equation}\label{eq:composite-rel}
\Lambdat_{a_{m+1},\alpha_{m+1},a_{m+1},c_m}\circ
\Lambdat_{c_m,\alpha_m,a_m,c_{m-1}}\circ\ldots\circ
\Lambdat_{c_2,\alpha_2,a_2,c_1}\circ
\Lambdat_{c_1,\alpha_1,a_1,a_1},
\end{equation}
and conversely, every element of the composite relation certainly
corresponds to an element of the forward generalized broken bicharacteristic
relation. For the sake of convenience we set $c_0=a_1$, $c_{m+1}=a_{m+1}$.
We now prove that
\eqref{eq:composite-rel} is a smooth Lagrangian.
As indicated
in the first paragraphs, a simple modification yields
Theorem~\ref{thm:sc-matrix-disc} if we also prove that the composite relation
maps $R_-(\lambda)\cap\sct_{C'_{a_1}}X_{a_1}$ to a finite union of
smooth Lagrangians.
We do not state this mapping property explicitly to avoid an overburdened
notation, but we indicate in Remark~\ref{rem:R_-} how it is proved by a simple
modification of the following argument.

\begin{prop}\label{prop:comp-Lagr}
For each string as in \eqref{eq:br-strings}
the composite relation \eqref{eq:composite-rel}
is a smooth Lagrangian submanifold of $T^*X_{a_{m+1}}\times T^*X_{a_1}$.
\end{prop}

The rest of this appendix is denoted to the proof of this proposition and
to the indication of the minor changes required for the proof of
Theorem~\ref{thm:sc-matrix-disc}.

\begin{proof}
The proof proceeds by induction. Thus, let
\begin{equation}
\Lambda^\sharp_{j}=
\Lambdat_{c_j,\alpha_j,a_j,c_{j-1}}\circ\ldots\circ
\Lambdat_{c_2,\alpha_2,a_2,c_1}\circ
\Lambdat_{c_1,\alpha_1,a_1,a_1};
\end{equation}
we need to show that for each $j\leq m+1$, given that $\Lambda^\sharp_{j-1}$ is
smooth Lagrangian, so is
$\Lambda^\sharp_{j}
=\Lambdat_{c_{j},\alpha_{j},a_{j},c_{j-1}}\circ\Lambda^\sharp_{j-1}$.
By the usual argument, this certainly follows if the intersection
of
$\Lambdat_{c_{j},\alpha_{j},a_{j},c_{j-1}}\times\Lambda^\sharp_{j-1}$
with the partial diagonal $\diag$ in $T^* X_{c_j}\times(T^* X_{c_{j-1}}
\times T^* X_{c_{j-1}})\times T^* X_{c_0}$ is transversal, i.e.\ for
$p\in (\Lambdat\times\Lambda^\sharp_{j-1})\cap\diag$,
\begin{equation}
T_p(\Lambdat\times \Lambda^\sharp_{j-1})
+T_p\diag=T_p(T^* X_{c_j}\times T^* X_{c_{j-1}}
\times T^* X_{c_{j-1}}\times T^* X_{c_0}).
\end{equation}
Here `partial' means that we take the diagonal in the central factor
$T^*X_{c_{j-1}}\times T^* X_{c_{j-1}}$,
and for the sake of simpler notation we wrote
\begin{equation}
\Lambdat=\Lambdat_{c_{j},\alpha_{j},a_{j},c_{j-1}}.
\end{equation}
In this case the tangent space of
the composite relation, $\Lambda^\sharp_j$ is given by the projection of
\begin{equation}
T_p(\Lambdat\times \Lambda^\sharp_{j-1})
\cap T_p\diag
\end{equation}
to the first and last factors, i.e.\ to $T^*X_{c_j}\times T^*X_{c_0}$.
The transversality is equivalent to
\begin{equation}
T(\Lambdat\times T^*X_{c_0})
+T(T^*X_{c_j}\times \Lambda^\sharp_{j-1})
=T(T^* X_{c_j}
\times T^* X_{c_{j-1}}\times T^* X_{c_0})
\end{equation}
over $(\Lambdat\times T^*X_{c_0})\cap(T^*X_{c_j}\times \Lambda^\sharp_{j-1})$,
and then $T\Lambda^\sharp_j$ is given by the projection of
\begin{equation}
T(\Lambdat\times T^*X_{c_0})
\cap T(T^*X_{c_j}\times \Lambda^\sharp_{j-1})
\end{equation}
to the first and third factors.
Transversality certainly follows if we can find a Lagrangian subspace $V'$ of
$TT^*X_{c_{j-1}}$ such that $V'$ is in the range of the
differential of the projection $\Lambda^\sharp_{j-1}
\to T^* X_{c_{j-1}}$, and
\begin{equation}
T\Lambdat+TT^*X_{c_j}\times V'
=TT^*X_{c_j}\times TT^*X_{c_{j-1}},
\end{equation}
and in this case the Lagrangian subspace
\begin{equation}
V=T\Lambdat\circ V'\subset TT^*X_{c_j},
\end{equation}
given by the projection of
$T\Lambdat\cap (TT^*X_{c_j}\times V')$ to $TT^*X_{c_j}$,
is in the range of the differential of the projection $\Lambda^\sharp_{j}
\to T^* X_{c_{j}}$.
Since the previous statements involving $V'$ referred to tangent vectors,
they are essentially equivalent to mapping properties of
$\Lambdat$ on Lagrangian submanifolds
$\Lambda'$
of $T^*X_{c_{j-1}}$. Thus, for the sake of convenience in notation,
we will consider these mapping properties, i.e.\ we will take $V'=T_p\Lambda'$;
it is simple to reinterpret results in the desired form.

To be concrete, for $j=1$ we take
\begin{equation}\label{eq:V_m-def}
V'=V'_1=\Span\{v\cdot\pa_{w'}+B_jv\cdot\pa_{\xi'}:\ v\in X_{c_1}\}.
\end{equation}
Note that
$V$ is certainly in the range of the
differential of the projection $\Lambdat_0=\Lambda^\sharp_0
\to T^* X_{c_1}$
by \eqref{eq:Lambdat-5} (the projection is to the
unprimed variables in the notation of this equation!).
If $a_2\neq c_1$,
an equally good choice is $V'=V'_1=X_{c_1}\subset T^*X_{c_1}$,
which simply corresponds to the Lagrangian given by
plane waves $\xi'_{1}=\text{const}$.

Before continuing, we make some general remarks about Lagrangian
submanifolds $\Lambda$ of $T^*X_a$. Suppose that $\Lambda$ is (locally)
a smooth graph
over $X_a$; i.e.\ that for every $p\in\Lambda$, the bundle
projection $T^*X_a\to X_a$ has a surjective differential at $p$ (which
is thus an isomorphism). Then
$\Lambda$ is (locally)
the image of a smooth bundle map $F:X_a\to T^*X_a$, i.e.\ the
composite of $F$ with the projection to the base is the identity map. Since
$T^*X_a=X_a\times X_a^*$ naturally, $F$ has the form $F(w)=(w,F_2(w))$
where $F_2:X_a\to X_a^*$. Moreover, $T_p T^*X_a$ can be naturally
identified with the vector space $T^*X_a$ itself; so for each $w\in X_a$,
$(F_2)_*|_w=(F_2)_*$ is
a linear map from $X_a$ to $X_a^*$. Let $A=A_w$ be the induced endomorphism
of $X_a$ via the metric identification of $X_a$ and $X_a^*$. That
$\Lambda$ is Lagrangian means that for any two tangent vectors
$V,V'\in T_p\Lambda$, $\omega(V,V')=0$, where $\omega=\omega_a$ is the standard
symplectic
form given by $d\xi\wedge dw$, so writing $V,V'\in T_p T^*X_a=X_a\times X_a^*$
as $V=(v,v^*)$, $V'=(v',(v')^*)$,
$\omega(V,V')=v^*(v')-(v')^*(v)$,
so in our case, with $V=F_* v$, $V'=F_* v'$, the Lagrangian condition
becomes $(F_2)_* v(v')-(F_2)_* v'(v)=0$. Since the metric identification
means that $(F_2)_* v(v')=g_a(Av,v')=Av\cdot v'$,
the Lagrangian condition amounts
to the statement that $A$ is self-adjoint.

A similar discussion also applies to $T^*X_{c_j}\times T^*X_{c_{j-1}}$ with
the twisted symplectic form $\omega_{c_j}-\omega_{c_{j-1}}$, and shows that
for $\Lambdat_{c_j,\alpha_j,a_j,c_{j-1}}$, which is a graph by
\eqref{eq:Lambdat-def-1}, with the notation of \eqref{eq:Lambdat-5},
$B_j$ and $B'_j$ are self-adjoint on $X_{c_j}$ and $X_{c_{j-1}}$
respectively (this of course follows from \eqref{eq:Lambdat-8} as well), while
$C'_j=C^*_j$.

Now we return to transversality.
In fact, we show the following stronger statement. Let $j$ be
arbitrary, $a=a_j$, $d=c_{j-1}$, $c=c_j$, $\alpha=\alpha_j$, and suppose that
$\Lambda'$ is a Lagrangian submanifold of $T^*X_d$
which is (locally) a smooth graph over $X_d$ of the
form $(w',\xi')=(w',F'(w'))$ where for each $w'\in X_d$,
$F'_*|_{w'}$ is given by a positive operator $A'=A'_{w'}$ on
$X_d$ (i.e.\ $A'\geq 0$;
$A'$ is automatically self-adjoint as discussed above), and suppose
in addition that $A'$ is positive definite if $d=a$.
Then we show that the relation $\Lambdat
=\Lambdat_{c\alpha ad}$ intersects $T^*X_c\times\Lambda'$ transversally, and
maps $\Lambda'\in T^*X_d$
to a Lagrangian $\Lambda\subset T^*X_c$ with similar properties,
i.e.\ $\Lambda$
is (locally)
a graph of the form $(w,F(w))$, $A=A_w=F_*|_w$ is positive for every
$w\in X_c$, and
$A$ is positive definite if $c\neq a$. Note
that the transversal intersection property automatically shows that
$\Lambda$ is smooth and Lagrangian. Once this is shown, the desired result
follows by induction starting from $j=1$. Indeed, to keep the induction
going we need to check that $A$ is positive definite if $c_j=a_{j+1}$,
but since this implies $c_j\neq a_j$, i.e.\ $c\neq a$, as discussed before,
we see that our claim, that $A$ is positive definite if $c\neq a$,
ensures that this holds.

We remark that for $j=1$, the Lagrangian subspace
$V_1$ of $T_{F(w'_1)}T^*X_{c_1}$ defined in
\eqref{eq:V_m-def} certainly has the
desired positivity; the operator $A'=A'_1$ is given by
$B_1$, hence it is positive definite unless $c_1=a_1$.
If instead we take $V_{1}
=X_{c_1}$ (corresponding to plane waves), then $A'=A'_1=0$, so it
is still positive (though not positive definite), which suffices
if $a_2\neq c_1$.

We proceed to show the claimed mapping property of $\Lambdat$.
We need to show that for $p=(w,\xi,w',\xi')\in(T^*X_c\times\Lambda')\cap
\Lambdat$
\begin{equation}\label{eq:transv-3}
T_p(T^*X_c\times \Lambda')+T_p\Lambdat=T_p(T^*X_c\times T^*X_d).
\end{equation}
But, due
to the first summand, $\pa_{w}$, $\pa_{\xi}$ are in the left hand side,
hence so is $v'\cdot\pa_{w'}+B'_jv'\cdot \pa_{\xi'}$, $v'\in
X_d$, due to
\eqref{eq:Lambdat-5}. Since the vectors
$v_0\cdot\pa_{w'}+A'v_0\cdot\pa_{\xi'}$,
$v_0\in X_d$, are
also in the left hand side due to the first summand, \eqref{eq:transv-3}
will follow if the block matrix $\begin{bmatrix} I & I\\ B'_j & A'
\end{bmatrix}$ is invertible, i.e.\ if $B'_j-A'$ is invertible.
But $B'_j\leq 0$, in fact negative definite if $a\neq d$, while
$A'\geq 0$, in fact positive definite if $a=d$, we conclude that
$A'-B'_j$ is positive definite, hence invertible. This proves
\eqref{eq:transv-3}.

Next, we need to find
\begin{equation}\label{eq:transv-8}
T_p(T^*X_c\times \Lambda')\cap T_p\Lambdat.
\end{equation}
Dimension counting shows that we will have found the whole intersection if
for every $v\in X_c$ we find an element of the intersection
which is of the form
\begin{equation}
v\cdot\pa_w+Av\cdot \pa_\xi+v_1\cdot\pa_{w'}+v_2\cdot\pa_{\xi'}.
\end{equation}
To do so, fix $v\in X_c$, and consider elements of $T_p\Lambdat$
which are of the form
\begin{equation}\label{eq:transv-16}
v\cdot\pa_w+B_j v\cdot\pa_{\xi}+C_j v\cdot\pa_{\xi'}
+v'\cdot\pa_{w'}+B'_j v'\cdot\pa_{\xi'}+C'_j
v'\cdot\pa_{\xi},\quad v'\in X_d.
\end{equation}
For these to be in the intersection \eqref{eq:transv-8}, we need
that there exist $v_0\in X_d$ such that the projection of
\eqref{eq:transv-16} to $T_{(w',\xi')}T^*X_d$ is equal to $v_0\cdot\pa_{w'}
+A'v_0\cdot\pa_{\xi'}$, i.e.\ such that
\begin{equation}
v'=v_0,\quad C_j v+B'_j v'=A'v_0.
\end{equation}
The existence of such $v'$ and $v_0$ thus again
follows from the invertibility of $\begin{bmatrix} I & -I\\ B'_j & -A'
\end{bmatrix}$, which holds since $A'-B'_j$ is
positive definite, and we see that in particular
$v'=(A'-B'_j)^{-1}C_j v$.
The projection of the corresponding tangent vector
in \eqref{eq:transv-8} to $T_{(w,\xi)}T^*X_c$ is given by
\begin{equation}
v\cdot\pa_w+(B_j v+C'_j v')\cdot\pa_\xi=v\cdot\pa_w
+(B_j+C^*_j(A'-B'_j)^{-1}C_j)v\cdot\pa_\xi
\end{equation}
where we used $C'_j=C_j^*$. Thus, the projection of the intersection
$(T^*X_c\times \Lambda')\cap \Lambdat$ to $T^*X_c$ is a smooth Lagrangian
$\Lambda$ with tangent vectors of the form
\begin{equation}
v\cdot\pa_w+Av\cdot\pa_\xi,\quad A=B_j+C^*_j(A'-B'_j)^{-1}C_j.
\end{equation}
Since $A'-B'_j$ is positive definite, the same holds for its inverse,
so $C^*_j(A'-B'_j)^{-1}C_j\geq 0$. Moreover, $B_j\geq 0$, so
we conclude that $A\geq 0$. If in addition $c\neq a$, then $B_j$ is
positive definite, hence the same holds for $A$.
This provides the necessary inductive step and proves
Proposition~\ref{prop:comp-Lagr}.
\end{proof}

\begin{rem}\label{rem:R_-}
The outgoing radial set, $R_-(\lambda)\cap\sct_{C'_{c_1}}X_{c_1}$ gives
another example of a Lagrangian $\Lambda'$ such that the operator $A'$,
defined above, is positive. Namely, a simple calculation shows
that on $\mu=0$, $\tau=-\sqrt{\lambda-\sigma}$, $\sigma\in\Lambda_1$,
$A'$ is given by
\begin{equation}
A'=\frac{\sqrt{\lambda-\sigma}}{|w|^3}(|w|^2\Id-w\otimes w),
\end{equation}
hence $A'\geq 0$. Thus, if $a_2\neq c_1$, which holds in this case (otherwise
the bicharacteristic were constant!), the proof of
Proposition~\ref{prop:comp-Lagr}
actually shows that the forward broken bicharacteristic relation maps
it to a finite union of smooth Lagrangian submanifolds.
\end{rem}

We now indicate the changes required to prove the S-matrix result,
Theorem~\ref{thm:sc-matrix-disc}.

\begin{proof}(Proof of Theorem~\ref{thm:sc-matrix-disc}.)
With the notation of \eqref{eq:br-strings}, if $m=0$, then $\gamma$ is
an unbroken bicharacteristic, so the endpoints of $\gamma$ are related
by the antipodal relation, which is certainly Lagrangian. So we may
assume $m\geq 1$. In addition, due to the assumptions \eqref{eq:param-Poisson},
the $\scHg^a$ bicharacteristics will only break upon hitting $C_{c_1}$
with $c_1\neq a$ (i.e.\ $a$-tangential bicharacteristics cannot become
normal to $a$ without hitting $C_{a,\sing}$). Similarly, $c_m\neq b$.

The only change from the proof of Proposition~\ref{prop:comp-Lagr}
is that we need to replace the first and last Lagrangians
in the composition \eqref{eq:composite-rel} by the Lagrangians
$\Lambdat_\beta\subset T^*C'_b\times T^*X_{c_m}$ and
$\Lambdat_\alpha\subset T^*X_{c_1}\times T^* C'_a$ given
by the (twisted) graphs of $d(-\sqrt{\lambda-\ep_\beta}\,\omega\cdot w')$
and $d(-\sqrt{\lambda-\ep_\alpha}\,\omega'\cdot w)$ respectively; here we
wrote $(\omega,\zeta,w',\xi')$ and $(w,\xi,\omega',\zeta')$ for
the coordinates on $T^*C'_b\times T^*X_{c_m}$ and $T^*X_{c_1}\times T^* C'_a$
respectively. Thus, on $\Lambdat_\alpha$, $\xi=-\sqrt{\lambda-\ep_\alpha}\,
\pi_{a,c_1}(\omega')$,
so the differential of the projection of $\Lambdat_\alpha$
to $T^*X_{c_1}$ is surjective since $a\neq c_1$. In particular, we
can choose a Lagrangian subspace $V'$ of $T_{(w,\xi)}T^*X_{c_1}$
which is in the
range of the differential of this projection, and which
is a graph over $X_{c_1}$ of the form $(w,Aw)$
with $A$ positive definite (e.g.\ the identity). Then the proof of
Proposition~\ref{prop:comp-Lagr} shows that the composite relation
\begin{equation}\label{eq:composite-rel-16}
\Lambdat_{c_m,\alpha_m,a_m,c_{m-1}}\circ\ldots\circ
\Lambdat_{c_2,\alpha_2,a_2,c_1}\circ
\Lambdat_\alpha
\end{equation}
is a smooth Lagrangian. Since the differential of the projection of
$\Lambdat_\beta$
to $T^*X_{c_m}$ is also surjective as $b\neq c_m$, the composition of
$\Lambdat_\beta$ with \eqref{eq:composite-rel-16}
is transversal, so the result is a smooth Lagrangian submanifold of
$T^*C'_b\times T^* C'_a$. In view of Theorem~\ref{thm:sc-matrix},
this and Remark~\ref{rem:R_-} complete the proof.
\end{proof}

\begin{rem}
The only reason for the individual
composite Lagrangians discussed in the above proof not
to be closed in $T^*C'_b\times T^*C'_a\setminus 0$
is that in the definition of the elementary Lagrangians,
\eqref{eq:Lambdat-def-1}, we excluded $w=w'$. Instead, the boundary points
corresponding to $w=w'$ were included in a composite Lagrangian
given by a different (shorter) string. As discussed after
\eqref{eq:br-strings}, we do not need to assume that the break points
are over the regular part of the $C_{c_j}$. Each of the composite Lagrangians
corresponding to a string \eqref{eq:br-strings}
is smooth through the singular part of the various $C_{c_j}$; they simply
intersect each other there.
\end{rem}

\bibliographystyle{plain}
\bibliography{sm}

\end{document}